\documentclass[10pt]{amsart}

\usepackage{amssymb, amscd}
\usepackage{epsfig, mathtools}
\usepackage[vmargin=1in, hmargin=1.3in]{geometry}
\usepackage[font=small,format=plain,labelfont=bf,up,textfont=it,up]{caption}
\usepackage{microtype}
\usepackage[shortlabels]{enumitem}
\usepackage[backref=page, bookmarks, bookmarksdepth=2, colorlinks=true, linkcolor=blue, citecolor=blue, urlcolor=blue]{hyperref}
\usepackage{etoolbox}
\usepackage{tikz-cd}
\usepackage{mathrsfs}  
\usepackage{float}
\usepackage{import}
\usepackage{comment}

\usepackage{caption}
\newcommand{\nc}{\newcommand}
\nc{\Z}{\mathbb{Z}}
\nc{\N}{\mathbb{N}}
\nc{\R}{\mathbb{R}}
\nc{\Q}{\mathbb{Q}}
\nc{\C}{\mathbb{C}}
\nc{\CP}{\mathbb{C} P}
\nc{\h}{\mathbb{H}}
\nc{\F}{\mathbb{F}}

\usepackage{transparent}
\usepackage[inkscapelatex=true]{svg}

\usepackage{caption}
\apptocmd{\thebibliography}{\raggedright}{}{}

\makeatletter
\patchcmd{\@maketitle}{\global\topskip42\p@\relax}
  {\global\topskip42\p@\relax \vspace*{-38pt}}
  {}{}
\makeatother

\setenumerate[0]{label=(\alph*)}

\renewcommand*{\backref}[1]{}
\renewcommand*{\backrefalt}[4]{%
    \ifcase #1 (Not cited.)%
    \or        (Cited on page~#2.)%
    \else      (Cited on pages~#2.)%
    \fi}

\numberwithin{equation}{section}

\theoremstyle{plain}
\newtheorem{theorem}{Theorem}[section]
\newtheorem{maintheorem}{Theorem}
\newtheorem{maincorollary}[maintheorem]{Corollary}

\newtheorem{proposition}[theorem]{Proposition}
\newtheorem{lemma}[theorem]{Lemma}

\newtheorem{corollary}[theorem]{Corollary}

\theoremstyle{definition}
\newtheorem{asm}[theorem]{Assumption}

\newtheorem{defn}[theorem]{Definition}

\newtheorem{notn}[theorem]{Notation}

\theoremstyle{remark}
\newtheorem{rmk}[theorem]{Remark}
\newenvironment{remark}[1][]{\begin{rmk}[#1] \pushQED{\qed}}{\popQED \end{rmk}}
\newtheorem{eg}[theorem]{Example}

\DeclareMathOperator{\Mod}{Mod}

\title{Flexible algebraic curves and fake vanishing cycles}

\author{Joshua Lehman and Tudur Lewis}
\address{Dept of Mathematics; Univ of Notre Dame; Hurley Hall; Notre Dame, IN 46556; USA \newline School of Mathematics, University of Bristol}
\email{jlehman4@nd.edu \newline et24788@bristol.ac.uk}
\thanks{}

\date{\today}

\begin{document}
\maketitle

\begin{abstract}

%(First as I would say to expert, then make more accesible) \newline 

%A surface in a 4-manifold is called flexible if any mapping class of the surface extends to a mapping class of the ambient 4-manifold, if the surface is characteristic, flexibility is obstructed 
We prove that complex curves in sufficiently ample linear systems on simply connected projective surfaces are topologically (spin) flexible: any (spin) mapping class of the underlying real surface extends to a topological mapping class of the ambient 4-manifold. As a corollary, we prove that smooth curves in the complex projective plane are topologically (spin) flexible. A consequence of topological flexibility is the existence of fake vanishing cycles, that is, simple closed curves on the underlying real surface that are not vanishing cycles, yet still bound a topologically embedded $(\pm1)$-framed disk in the surface exterior. We also give examples of pushed in Seifert surfaces in $D^4$ that are not topologically spin-flexible, but become spin-flexible after a single internal stabilization. Our methods utilize the equivariant intersection form on the surface exterior. \newline

%V1: A surface in a 4-manifold is topologically (spin) flexible if any (spin) mapping class of the surface extends to a topological mapping class of the 4-manifold. We prove that complex curves in sufficiently ample linear systems on simply connected projective surfaces are topologically (spin) flexible. As a corollary, we construct fake vanishing cycles, that is, simple closed curves on the underyling real surface that are not vanishing cycles, yet still bound a topologically embedded $(\pm1)$-framed disk in the surface exterior. We give examples of pushed-in Seifert surfaces in $D^4$ that are not topologically spin-flexible, but become spin-flexible after a single internal stabilization. Our methods utilize the equivariant intersection form on the surface exterior. 

%sufficient conditions under which complex curves in simply connected projective surfaces are topologically (spin) flexible. 

%\textcolor{red}{We also realize certain non-trivial proper subgroups of $\mathrm{Sp}_{2g}(\mathbb{F}_2)$ as the image of topological extendable subgroups under the mod 2 symplectic representation.} 
% We compute the extendable subgroup of smooth curves in sufficiently ample linear systems on simply connected projective surfaces. When the curve is a characteristic surface, the extendable subgroup  we find they are 

%We prove that smooth curves in projective surfaces are topologically flexible surfaces, that is, 
\end{abstract}

\section{Introduction}

Fix a properly embedded oriented surface $\Sigma$ in a compact oriented 4-manifold $M$. Consider the group $\pi_0(\mathrm{Homeo}_{\partial}(M,\Sigma))$ of isotopy classes of orientation preserving homeomorphisms of $M$ which fix $\Sigma$ setwise and $\partial M$ pointwise. Restriction determines a homomorphism 
\[ \mathcal{R}\colon \pi_0(\mathrm{Homeo}_{\partial}(M,\Sigma)) \to \pi_0(\mathrm{Homeo}_{\partial}(\Sigma)). \]
The subgroup $\Mod(\Sigma) %= \pi_0(\mathrm{Homeo}_{\partial}^+(\Sigma))
\leq \pi_0(\mathrm{Homeo}_\partial(\Sigma))$ preserving the orientation of $\Sigma$ is the \textit{mapping class group} of $\Sigma$. The \textit{topological extendable subgroup} $\mathcal{E}^{\mathrm{TOP}}(M,\Sigma)$ of $\Sigma$ in $M$ is the image of $\mathcal{R}$ in $\Mod(\Sigma)$. The surface $\Sigma$ is \textit{topologically flexible} in $M$ if $\mathcal{E}^{\mathrm{TOP}}(M,\Sigma) = \Mod(\Sigma)$. If $\Sigma \subset M$ is characteristic and $H_1(M;\mathbb{Z}) = 0$, Rochlin defined a quadratic form $q_\Sigma\colon H_1(\Sigma;\mathbb{Z}/2\mathbb{Z}) \to \mathbb{Z}/2\mathbb{Z}$ which obstructs flexibility:
\[ \mathcal{E}^{\mathrm{TOP}}(M,\Sigma) \leq  \Mod(\Sigma)[q_\Sigma] \]
where $\Mod(\Sigma)[q_\Sigma]$ denotes the stabilizer of $q_\Sigma$ in $\Mod(\Sigma)$.
The surface $\Sigma$ is \textit{topologically spin-flexible} in $M$ if $\mathcal{E}^{\mathrm{TOP}}(M,\Sigma) = \Mod(\Sigma)[q_\Sigma]$. %There are several calculations of the (smooth) extendable subgroup in the literature.  %The (smooth) flexibility of curves in $\mathbb{C}P^2$ of degree at least six has remained open. By working with the topological extendable subgroup, we obtain a slew of additional techniques for extending mapping classes. We apply this to establish (spin) flexibility of complex curves in $\mathbb{C}P^2$. Our methods apply to curves in linear systems with sufficiently large monodromy groups, and we state our result in this level of generality.   %\textcolor{red}{Transition}. 

We find sufficient conditions under which complex curves in simply connected projective surfaces are topologically (spin) flexible. 

%\begin{question} \label{flexofcurves}
%Fix a complex curve $C$ in a simply connected projective surface $X$. Is $C$ topologically (spin) flexible? %In particular, are complex curves in $\mathbb{C}P^2$ topologically (spin) flexible?
%\end{question}

%We answer this question affirmatively whenever the linear system associated to $C$ has \textit{$r$-spin monodromy}. In particular, we establish that complex curves in $\mathbb{C}P^2$ are topologically (spin) flexible surfaces. Hirose established smooth (spin) flexibility for cubics and quartics, whilst it is a corollary of Salter that quintics are (smoothly) spin-flexible. 

%Hirose proved that cubics and quartics in $\mathbb{C}P^2$ are (smoothly) flexible, whilst work of Salter implies that quintics are (smoothly) spin-flexible.
%We find sufficient conditions under which the answer to Question \ref{flexofcurves} is affirmative. 
Let $X$ be a smooth projective surface and $L \in \mathrm{Pic}(X)$ a line bundle\footnote{Recall that $\mathrm{Pic}(X) = H^1(X,\mathcal{O}_X^*)$ is the group of isomorphism classes of holomorphic line bundles on $X$, and, as $\pi_1(X)=1$, the first Chern class gives an injective homomorphism $c_1 : \mathrm{Pic}(X) \to H^2(X;\mathbb{Z})$ \cite[1.1]{GH}.}. Recall that holomorphic sections of $L$, elements of $H^0(X;L)$, have the same zero locus in $X$ if and only if they are scalar multiples of one another \cite[1.4]{GH}. Therefore elements of $|L| := \mathbb{P}H^0(X;L)$ are identified with complex curves in $X$. 
A source of extendable mapping classes for a smooth curve $C \in |L|$ is the \textit{monodromy} $\Gamma_L \leq \Mod(C)$ of the linear system $|L|$: a homotopy class of a loop in the smooth locus $U_L \subset |L|$ based at $C$ produces a loop of smooth curves in $X$, and in particular, a mapping class of $C$. By isotopy extension, 
\[\Gamma_L \leq \mathcal{E}^{\mathrm{TOP}}(X,C).\] 
For sufficiently ample line bundles, our work establishes a difference between these two groups. Recall that $L$ is \textit{very ample} (or $1$\textit{-jet ample}) if there is a holomorphic embedding $\varphi : X \to \mathbb{C}P^N$ such that $L = \varphi^*H$ where $H$ is the dual of the tautological line bundle on $\mathbb{C}P^N$ \cite[1.4]{GH}. Beltrametti and Sommese \cite{kjetampleref} introduced a quantitative refinement, known as \textit{k-jet ampleness},  with the key property that if $L_1$ is $k_1$-jet ample and $L_2$ is $k_2$-jet ample, $L_1\otimes L_2$ is $(k_1+k_2)$-jet ample (Proposition \ref{ampleadd}). %\textcolor{blue}{Recall definition of very ample}.   %We say a line bundle $L$ is \textit{sufficiently ample} if it is $k$-jet ample (Definition \ref{kjetample}) for $k\geq 6$. % If the line bundle $L$ is sufficiently ample and $\pi_1(X) = 1$, Banerjee-Salter \cite[Theorem A]{banerjee2025monodromyvanishingcyclessufficiently}

%We say $|L|$ has \textit{$r$-spin monodromy} if $\Gamma_L$ is an \textit{$r$-spin mapping class group}: the stabilizer of a mod $r$ cohomology class on the unit tangent bundle of $C$. Here $r$ is necessarily the maximal root of the adjoint line bundle $K_X \otimes L$.  %When $\Gamma_L$ is an $r$-spin mapping class group, we say $|L|$ has \textit{$r$-spin monodromy}.  

\begin{maintheorem} \label{maintheorem}
%Let $C_d \subset \mathbb{C}P^2$ be a smooth degree $d$ plane curve. If $d$ is even, $C_d$ is topologically spin-flexible. If $d$ is odd, $C_d$ is topologically flexible. 
Let $X$ be a projective surface with $\pi_1(X) = 1$. Let $L \in \mathrm{Pic}(X)$. Suppose that $L$ admits a decomposition $L= L_1 \otimes L_2$ for $L_1$ very ample and $L_2$ 6-jet ample. Then smooth curves in the linear system $\mathbb{P}H^0(X;L)$ are topologically (spin) flexible. %if $r$ is odd, and topologically spin-flexible if $r$ is even. 
\end{maintheorem}

Our approach to Theorem \ref{maintheorem} uses work of Conway-Powell \cite[Theorem 4.1]{Conway_Powell} to produce new homeomorphisms of the pair $(X,C)$ and work of Banerjee-Salter \cite[Theorem A]{banerjee2025monodromyvanishingcyclessufficiently} which computes the monodromy of the linear system $|L|$. We highlight an especially concrete application of Theorem \ref{maintheorem}.

\begin{maincorollary} \label{flexcomplexcurvesintro}
Smooth plane curves in $\mathbb{C}P^2$ are topologically (spin) flexible. Smooth curves in $\mathbb{C}P^1 \times \mathbb{C}P^1$ of bidegree $(d_1,d_2)$ with $d_1, d_2\geq 6$ are topologically (spin) flexible. %Moreover, there exist (spin) flexible smooth curves in the Hirzebruch surface $\mathbb{F}_1 \cong \mathbb{C}P^2 \# \overline{\mathbb{C}P^2}$, as well as in complete intersections. 
\end{maincorollary}

%\begin{rmk}[\textbf{Which curves are characteristic?}]
%In Theorem \ref{maintheorem}, a smooth curve $C \in |L|$ is characteristic if and only if $r$ is even \cite{Libgober}. 
%\end{rmk}

Hirose proved that cubics and quartics in $\mathbb{C}P^2$ are (smoothly) flexible \cite{HiroseP2}. Work of Salter implies (smooth) spin-flexibility of quintics \cite{salterquintic}. Whilst Corollary \ref{flexcomplexcurvesintro} establishes topological (spin) flexibility of plane curves in all remaining degrees, their \textit{smooth flexibility} (Definition \ref{extendablegroup}) remains open. We note that after a single \textit{internal stabilization}, plane curves become smoothly (spin) flexible \cite{lehman}. 

Recent work of Baykur and Sunukjian \cite{baykur}, in conjunction with Corollary \ref{flexcomplexcurves}, establishes the existence of surfaces in $\mathbb{C}P^2$ that are topologically isotopic to plane curves, so are topologically flexible, but are \textit{projectively rigid} in the smooth category: the image of their smooth extendable subgroup under the \textit{projective symplectic representation} is trivial \cite[Remark 19]{baykur}.

%\begin{rmk}
%For a smooth curve $C \in |L|$, all homeomorphisms of the pair $(X,C)$ constructed in this paper are topologically isotopic to the identity. %In Theorem \ref{maintheorem}, a smooth curve $C \in |L|$ is characteristic if and only if $r$ is even \cite{Libgober}.
%\end{rmk}
%\begin{rmk}[\textbf{Refinement of Theorem \ref{maintheorem}}] %and systems with $r$-spin monodromy}]
%Theorem \ref{maintheorem} follows from a more technical version, Theorem \ref{thmA}. which utilises work  work of Banerjee-Salter \cite{banerjee2025monodromyvanishingcyclessufficiently} and Conway-Powell \cite{Conway_Powell}.
%Precisely, Theorem \ref{thmA} requires that $L \in \mathrm{Pic}(X)$ is \textit{6-jet ample} and the linear system $|L|$ has \textit{r-spin monodromy} (Definition \ref{rspinmon}). Banerjee-Salter prove that the line bundles considered in Theorem \ref{maintheorem} satisfy the $r$-spin monodromy condition \cite[Theorem A]{banerjee2025monodromyvanishingcyclessufficiently}. 

   %is an application of \cite[Theorem 4.1]{Conway_Powell} to the surface exterior (recorded here as Theorem \ref{extensionhomeo}). We implement this by finding a subsurface $F \subset C$ contained in a 4-ball $D^4 \subset X$ with local knot group $\pi_1(D^4-F)\cong \mathbb{Z}$. The surface $F$ is a \textit{pushed in Seifert surface}. 

\vspace{0.2cm}

\noindent \textbf{Fake vanishing cycles.} A \textit{vanishing cycle} on a smooth complex curve $C \in |L|$ is a simple closed curve on $C$ which contracts to a single point under a \textit{nodal degeneration} of $C$ in $X$. Under this degeneration, the curve traces out a $(-1)$-framed disk in $X$ with interior disjoint from $C$. Moreover, the monodromy of this degeneration can be used to exhibit the \textit{Dehn twist} along the vanishing cycle as an element in $\Gamma_L$, and in particular, as a (smoothly) extendable mapping class. For details, we refer the reader to \cite{SingularitiesDiffMaps}. We call a simple closed curve on $C$ a \textit{fake vanishing cycle} if it is not a vanishing cycle on $C$, but still bounds a topologically embedded $(\pm1)$-framed disk in $X$ with interior disjoint from $C$. By Proposition \ref{HBtrick}, the Dehn twist on any fake vanishing cycle is a topologically extendable mapping class. Let $K_X \in \mathrm{Pic}(X)$ be the \textit{canonical bundle} of $X$ (the top exterior power of the holomorphic cotangent bundle of $X$).

\begin{maincorollary} \label{mainfakevanishing}
Let $X$ be a projective surface with $\pi_1(X) = 1$. Let $L \in \mathrm{Pic}(X)$ satisfy the hypothesis of Theorem \ref{maintheorem}. Fix a smooth curve $C \in |L|$. If the maximal root of the adjoint line bundle $K_X \otimes L$ in $\mathrm{Pic}(X)$ has order at least three, then there exists a fake vanishing cycle $a \subset C$.  %$|L|$ has $r$-spin monodromy (Definition \ref{rspinmon}) and that $r \neq 1,2$. 
%Fix a smooth curve $C \in |L|$. Then there exists a simple closed curve $\gamma \subset C$ which is a fake vanishing cycle.
%Let $X$ be a projective surface with $\pi_1(X) = 1$. Let $L \in \mathrm{Pic}(X)$ be 6-jet ample (Definition \ref{kjetample}). Assume that $|L|$ has $r$-spin monodromy (Definition \ref{rspinmon}) and that $r \neq 1,2$. Fix a smooth curve $C \in |L|$. Then there exists a simple closed curve $\gamma \subset C$ which is a fake vanishing cycle. 
\end{maincorollary}

\begin{rmk}[\textbf{Fake vanishing cycles don't always exist}]
%If $r=1$ or $2$ and the linear system $|L|$ has $r$-spin monodromy (flexibility, resp. spin-flexibility), fake vanishing cycles do not exist. For an explicit example, 
Recall that the \textit{hyperplane bundle} $\mathcal{O}(1) \to \mathbb{C}P^2$ is the dual of the tautological line bundle on $\mathbb{C}P^2$. Let $C_d \subset \mathbb{C}P^2$ be a smooth curve of degree $d$. Then $C_d$ may be identified as the zero locus of a section of $\mathcal{O}(d) = \mathcal{O}(1)^{\otimes d}$ \cite[Chapter 1, Section 3]{GH}. The adjoint line bundle of $\mathcal{O}(d)$ is $K_{\mathbb{C}P^2} \otimes \mathcal{O}(d) \cong \mathcal{O}(d-3) \cong \mathcal{O}(1)^{\otimes (d-3)}$. It follows from Corollary \ref{mainfakevanishing} that complex curves of degree $d \geq 6$ have fake vanishing cycles\footnote{The careful reader will note that smooth curves of degree $d=6$ are not covered by Theorem \ref{maintheorem}. This case is covered by Theorem \ref{thmA} which we use to prove Theorem \ref{maintheorem}.}. On the other hand, complex curves of degree $d \leq 5$ do not. This is due to the fact that the monodromy $\Gamma_{\mathcal{O}(d)}$ (Definition \ref{monodromydefn}) of the associated linear system is either the entire mapping class group $(d=3,4)$, or the stabilizer of Rochlin's form $(d=5)$.

\vspace{0.2cm}

%Fake vanishing cycles do not exist on $C_d$ if $d \leq 5$. This follows from the fact that the monodromy group $\Gamma_{\mathcal{O}_{\mathbb{C}P^2}(d)}$ is either the entire mapping class group $\Mod(C_d)$ ($d=3,4$) or the stabilizer of Rochlin's form $(d=5)$. On the other hand, they do exist if $d \geq 6$, in which case the adjoint line bundle is  monodromy group is the stabilizer of a $(d-3)$-spin structure (a proper subgroup of $\Mod(C_d)$).  % Indeed, for $d = 3,4,5$ we have $\Gamma_{|\mathcal{O}(d)|} = \mathcal{E}^{\mathrm{TOP}}(\mathbb{C}P^2,C_d)$. 
\end{rmk}

%\begin{proof}
%By Theorem \ref{maintheorem}, $C$ is topologically (spin) flexible. On the other hand, homeomorphisms of the pair $(X,C)$ will preserve (up to sign) the relative Euler number of  an embedded disk $D \subset X$ with interior disjoint from $C$ and $\partial D$ a simple closed curve on $C$. 
%\end{proof}

For surfaces with knot group $\mathbb{Z}$, Conway-Powell \cite{Conway_Powell} show that the \textit{equivariant intersection form} on the surface exterior gives a \textit{complete} obstruction to extending mapping classes topologically. In the proof of Theorem \ref{maintheorem}, the results in \cite{Conway_Powell} are implemented by finding a subsurface $F \subset C$ contained in a 4-ball $D^4 \subset X$ with local knot group $\pi_1(D^4-\nu F)\cong \mathbb{Z}$. The surface $F$ is a \textit{pushed in Seifert surface}. Their results allow us to produce homeomorphisms of the pair $(D^4,F)$ by constructing suitable isometries of the equivariant intersection form on $D^4-\nu F$. On the other hand, the equivariant intersection form can also be used to \textit{obstruct} extendable mapping classes. 

\vspace{0.2cm}

\noindent \textbf{Obstructing extendable mapping classes.} If $F \subset S^3 = \partial D^4$ is a Seifert surface, consider the properly embedded surface (again denoted) $F \subset D^4$ obtained by pushing the interior of $F$ into the interior of $D^4$. We refer to such a surface as a \textit{pushed in Seifert surface}. Let $T \subset S^4$ denote an unknotted torus. A \textit{trivial internal stabilization} of $F \subset D^4$ is the relative connected sum $(D^4,F) \# (S^4,T)$. 

\begin{maintheorem} \label{nonflex}
There exist infinitely many non-isotopic pushed in Seifert surfaces in $D^4$ that are not topologically spin-flexible rel boundary, but become topologically spin-flexible after a trivial internal stabilization. 
\end{maintheorem}

The obstruction in Theorem \ref{nonflex} arises from a specialization of the equivariant intersection form on the exterior of $F$ in $D^4$. 

\begin{rmk}
Several approaches to obstructing extendable mapping classes in the literature have exploited the \textit{peripheral structure}\footnote{If $\nu \Sigma$ is a closed tubular neighbourhood of $\Sigma$ in $M$, the \textit{peripheral subgroup} of $\Sigma$ in $M$ is the image of the map $\iota_* : \pi_1(\partial \nu \Sigma) \to \pi_1(M-\Sigma)$ induced by inclusion.} associated to an embedding $\Sigma \subset M$. For example, work of Hirose \cite{HiroseTori} used the peripheral structure to compute the extendable subgroup of (turned) spun tori in $S^4$. Recently, Niu \cite{WeizheTrivial} used the peripheral subgroup and \textit{iterated} rim surgery to produce closed surfaces in any simply connected 4-manifold with \textit{trivial} topological extendable subgroup (a topologically \textit{rigid} surface). Note that \cite{baykur}  constructed a totally geodesic topologically rigid surface in a hyperbolic 4-manifold using Mostow rigidity and the Dehn-Nielsen-Baer Theorem.  The peripheral subgroup of pushed in Seifert surfaces is too simple to obstruct extendable mapping classes, indeed it is  $\mathbb{Z}$, generated by a meridian of the surface.  
\end{rmk}

%Precisely, if $\nu \Sigma \subset M$ denotes a closed tubular neighbourhood of $\Sigma$ in $M$, inclusion induces homomorphism $\iota_* : \pi_1(\partial \nu \Sigma) \to \pi_1(M-\Sigma)$. 

Problem 4.39 on the K3 list \cite{K3} asks which subgroups of $\Mod(\Sigma)$ arise as extendable subgroups. A linearization of this problem asks which subgroups of $\mathrm{Sp}_{2g}(\mathbb{Z}/2\mathbb{Z})$ arise as the image of extendable subgroups under the mod 2 symplectic representation (see \cite{subgroupsofmod2sp}). We show that \textit{mod 2 Weyl groups} of type $A$ and $E$ may be realized.

The \textit{Milnor fiber} $F_\Gamma \subset D^4$ for the link $L$ of a singularity of type ADE with Dynkin diagram $\Gamma$ is, up to isotopy, a pushed in Seifert surface. We equip the surface $F_\Gamma$ with an ordered configuration $\mathcal{C}_\Gamma$ of oriented curves satisfying the following properties: the homology classes of curves in $\mathcal{C}_\Gamma$ generate $H_1(F_\Gamma;\mathbb{Z})$, the \textit{intersection graph} associated to $\mathcal{C}_\Gamma$ is $\Gamma$, and if $S$ is the Seifert form of $F_\Gamma$ relative to $\mathcal{C}_\Gamma$, then $S+S^T$ is the \textit{Cartan matrix} of $\Gamma$.

%Let $F \subset D^4$ be the \textit{Milnor fiber} of an isolated hypersurface singularity of type ADE. Up to isotopy, $F$ is a pushed-in minimal genus Seifert surface whose boundary is the closure of a positive braid. Moreover, the surface $F$ is equipped with a distinguished configuration $\mathcal{C}$ of curves whose intersection graph is the corresponding ADE Dynkin diagram. 

%Given a Coxeter-Dynkin diagram $\Gamma$, there is an associated minimal genus Seifert surface $F$ obtained by plumbing negative Hopf bands according to $\Gamma$. The surface $F$ is equipped with a distinguished configuration of simple closed curves $\{\gamma_i\}$. 

%A fundamental question is which subgroups of the mapping class group arise as  extendable subgroups \cite[Problem 4.39]{K3}. 

%A simpler version of this question asks which subgroups of $\mathrm{Sp}_{2g}(\mathbb{F}_2)$ arise as the image of extendable subgroups under the mod 2 symplectic representation?  

%In this direction we show that mod 2 Weyl groups of type $A$ and $E$ may be realized. \textcolor{red}{Recall basic setup}.  

\begin{maintheorem} \label{realizationofWeyl}
 Let $\Gamma$ be a Coxeter--Dynkin diagram of type $A_n$ with $n$ even, $E_6$ or $E_8$.
    Let $W_2(\Gamma)$ be the subgroup of $\mathrm{Sp}(H_1(F_\Gamma; \Z_2))$ generated by the transvections in the homology classes of the curves in $\mathcal{C}_\Gamma$. This is the mod 2 reduction of the Weyl group of type $\Gamma$. Then 
    \[ \Psi_2(\mathcal{E}^{\mathrm{TOP}}(D^4,F_\Gamma)) = W_2(\Gamma) \]
    where $\Psi_2 : \Mod(F_\Gamma) \to \mathrm{Sp}(H_1(F_\Gamma; \Z_2))$ is the mod 2 symplectic representation. %Suppose $f \in \epsilon^{\mathrm{DIFF}}( D^4, F)$, then $f_* \in W(\Gamma)$.   
\end{maintheorem}

\begin{rmk}
The singularities $D_n$, $A_n$ for $n$ odd and $E_7$ are excluded from Theorem \ref{realizationofWeyl} because our methods rely on the associated Seifert surface having boundary a \textit{knot} in $S^3$, rather than a link of more than one component.  %the associated Seifert Our methods rely on the boundary of the associated Seifert surface being a knot in $S^3$, 
\end{rmk}

 %In the literature, all obstructions in the topological category have used the \textit{peripheral structure} of the surface (which is trivial for the case of pushed-in Seifert surfaces in $B^4$). 

%\subsection{Context and prior work}

\subsection{General overview}

\subsubsection{Theorem \ref{maintheorem}/Theorem \ref{thmA}} \label{overview} We begin by finding a 4-ball $D^4 \subset X$ and a smooth curve $C \in |L|$ such that $F = C \cap D^4$ has knot group $\pi_1(D^4-F) \cong \mathbb{Z}$ and sufficiently large genus (the surface $F$ is a pushed in Seifert surface of genus three). This is where the 6-jet ampleness of $L$ is used (Lemma \ref{localmodel}). We then use results of Conway-Powell \cite{Conway_Powell} to upgrade isometries of the equivariant intersection form on the exterior of $F$ in $D^4$ to homeomorphisms of the pair $(D^4,F)$ rel boundary. In particular, we obtain that the \textit{Torelli subgroup} of $\Mod(F)$ is topologically extendable: \[\mathcal{I}(F) \leq \mathcal{E}^{\mathrm{TOP}}(D^4,F).\] These homeomorphisms may be extended by the identity to the rest of $X$, yielding extendable mapping classes for the pair $(X,C)$. By assumption (Theorem \ref{thmA}), the linear system $|L|$ has $r$-spin monodromy (Definition \ref{rspinmon}). This implies that $\mathcal{E}^{\mathrm{TOP}}(X,C)$ contains an $r$-spin mapping class group (Corollary \ref{rspinisextendable}). We conclude by exploiting how $r$-spin mapping class groups intersect this Torelli subgroup.  

%\begin{rmk}
%By constructing the homeomorphisms locally, we do not need complete control of the knot group $\pi_1(X-C)$. 
%\end{rmk}

\subsubsection{Theorems \ref{nonflex} and \ref{realizationofWeyl}}
The obstruction in Theorem \ref{realizationofWeyl} is a consequence of the fact that homeomorphisms of the pair $(D^4,F)$ can be promoted to isometries of the intersection form of the 2-fold branched cover $\Sigma_2(D^4,F)$ of $D^4$ branched over $F$. We realize this as a specialization of the equivariant intersection form on the exterior of $F$ in $D^4$ (c.f. Lemma \ref{thm:obtruction_for_pushed_in_Seifert_surfaces}). The Seifert surfaces associated to the $A_n$-family for $n$ even are the family used in Theorem \ref{nonflex}. Their flexibility after a single stabilization is a computation using basic properties of $r$-spin mapping class groups.   

%We give a concrete application of this obstruction to Milnor fibers of ADE singularities in To see that this is a genuine obstruction, we prove that $\Psi_2(\mathcal{E}^{\mathrm{TOP}}(D^4,F)) \leq W(\Gamma)$. The Seifert surfaces associated to the $A_n$-family for $n$ even then give the first part of Theorem \ref{nonflex}. Their flexibility after a single stabilization is a computation using properties of $r$-spin mapping class groups.   

\begin{rmk}
Recent work of Conway and Kasprowski \cite[Theorem 1.12]{conwaykasprowski} gives a complete algebraic characterization of the topological extendable subgroup whenever the knot group of the surface is finite cyclic. For a smooth degree $d$ plane curve $C_d \subset \mathbb{C}P^2$, we have $\pi_1(\mathbb{C}P^2-C_d) \cong \mathbb{Z}_d$. Whilst it is unclear to the authors if the computation \cite[Theorem 1.12]{conwaykasprowski} immediately implies the topological (spin) flexibility of $C_d$, we note that it gives a direct approach to constructing homeomorphisms of the pair $(\mathbb{C}P^2,C_d)$. In this paper we use \cite[Theorem 4.1]{Conway_Powell} to extend homeomorphisms, and the analogous criterion for finite cyclic knot groups is \cite[Theorem 1.6]{conwaykasprowski}. 
\end{rmk}

\subsection{Outline of the paper} Section 2 defines the extendable subgroup. Section 3 recalls the necessary background on $r$-spin mapping class groups and discusses how $r$-spin mapping class groups intersect the Torelli subgroup (Lemma \ref{maintechnicallemma}). Section 4 begins by recalling techniques for extending homeomorphisms of the boundary of a 4-manifold with fundamental group $\mathbb{Z}$  due to Conway-Powell \cite{Conway_Powell} (Theorem \ref{extensionhomeo}). From this we extract methods for extending mapping classes (Corollary \ref{torelliextends}). We give an explicit treatment of this procedure for pushed in Seifert surfaces in $D^4$, illustrated by generalizing classical methods (Proposition \ref{HBtrick}) for extending mapping classes in the smooth category (Corollary \ref{GeneralHBtrick}). Section 5 proves Theorems \ref{nonflex} and \ref{realizationofWeyl}. Section 6 proves Theorem \ref{maintheorem}.

%\subsection{Notation and conventions}
%The mapping class group of a surface $\Sigma$ is denoted $\Mod(\Sigma) = \pi_0(\mathrm{Homeo}_\partial^+(\Sigma))$. Throughout, $\Sigma$ denotes a properly embedded surface in a 4-manifold $M$, $F \subset D^4$ a Seifert surface with its interior pushed into $D^4$, and $C \subset X$ a complex curve. Isometries of the equivariant intersection form are denoted 

%\textcolor{red}{Use Mod to denote the mapping class group? Use $\Sigma$ for an arbitrary surface properly embedded in a $4$-manifold $M$? Use $F$ to denote a Seifert surface with its interior pushed into $D^4$? Use $C \subset X$ in the context of a curve in a complex surface? Denote isometries of equivariant forms over $\Lambda$ by $\alpha$? And denote matrices representing them by $A$? Should we use $\mathrm{Aut}(H,\lambda)$ or $\mathrm{Isom}(H,\lambda)$? Agree with all of this! Could we use $\mathrm{Iso}(H,\lambda)$? We can use $\alpha,\beta \in \mathrm{Iso}(H,\lambda)$ and reserve $a,b,c,\ldots$ for simple closed curves on surfaces? Works for me!}

\subsection{Acknowledgements} The authors thank Anthony Conway and Mark Powell for email correspondence regarding applications of their work \cite{Conway_Powell} to computations of extendable subgroups. The authors thank Larry Taylor for useful conversations regarding \cite{KirbyTaylor}. The first author thanks his advisor, Nick Salter, for his continued support and guidance.

\subsection{AI Declaration} No AI tools were used to prove any of the new results in this paper. ChatGPT assisted the authors in searching the literature and proofreading.

\section{Extendable mapping classes and Rochlin's form}

 Let $M$ be a compact, oriented, smooth 4-manifold. Let $\Sigma \subset M$ be a smooth, compact, properly embedded, oriented surface in $M$. 
 
 \subsection{Extendable subgroups} The \textit{mapping class group} of $\Sigma$ is the group \(\Mod(\Sigma) = \pi_0(\mathrm{Homeo}_{\partial}^+(\Sigma))\) of isotopy classes of orientation preserving homeomorphisms of $\Sigma$ that fix $\partial \Sigma$ pointwise. Let $\pi_0(\mathrm{Homeo}_{\partial}(M,\Sigma))$ (resp. $\pi_0(\mathrm{Diff}_{\partial}(M,\Sigma))$) denote the group of isotopy classes of orientation preserving homeomorphisms (resp. diffeomorphisms) of $M$ which fix $\Sigma$ setwise and $\partial M$ pointwise. Restriction determines a homomorphism $\mathcal{R}^{\mathrm{TOP}}\colon \pi_0(\mathrm{Homeo}_{\partial}(M,\Sigma)) \to \pi_0(\mathrm{Homeo}_{\partial}(\Sigma))$ (resp. $\mathcal{R}^{\mathrm{DIFF}}\colon \pi_0(\mathrm{Diff}_{\partial}(M,\Sigma)) \to \pi_0(\mathrm{Homeo}_{\partial}(\Sigma))$).  
 
 %Let $\mathrm{Homeo}^+_{\partial}(\Sigma)$ denote the topological group consisting of orientation preserving homeomorphisms of $\Sigma$ that fix $\partial \Sigma$ pointwise. The \textit{mapping class group} of $\Sigma$ is \[\Mod(\Sigma) = \pi_0(\mathrm{Homeo}^+_{\partial}(\Sigma)).\] Let $\mathrm{Homeo}_{\partial}^{+_\Sigma}(M,\Sigma)$ denote the topological group consisting of homeomorphisms of $M$ which fix $\Sigma$ setwise, preserve the orientation of $\Sigma$, and restrict to the identity on $\partial M$. Restriction determines a homomorphism \[\mathcal{R} : \pi_0(\mathrm{Homeo}^{+_\Sigma}_{\partial}(M,\Sigma)) \to \Mod(\Sigma).\] 

\begin{defn} \label{extendablegroup}
The \textit{topological (resp. smooth) extendable subgroup} $\mathcal{E}^{\mathrm{TOP}}(M,\Sigma)$ (resp. $\mathcal{E}^{\mathrm{DIFF}}(M,\Sigma)$) of $\Sigma$ in $M$ is the image of $\mathcal{R}^{\mathrm{TOP}}$ (resp. $\mathcal{R}^{\mathrm{DIFF}}$) in $\Mod(\Sigma)$. %The \textit{smooth extendable subgroup} of $\Sigma$ in $M$ is the image of $\mathcal{R}^{\mathrm{DIFF}}$ in $\Mod(\Sigma)$.  
%\[\mathcal{E}^{\mathrm{TOP}}(M,\Sigma) = \mathcal{R}(\pi_0(\mathrm{Homeo}_{\partial}(M,\Sigma))) \cap \Mod(\Sigma).\]
We say that $\Sigma$ is \textit{topologically (resp. smoothly) flexible} if $\mathcal{E}^{\mathrm{TOP}}(M,\Sigma) = \Mod(\Sigma)$ (resp. $\mathcal{E}^{\mathrm{DIFF}}(M,\Sigma) = \Mod(\Sigma))$. 
\end{defn}

If the surface $\Sigma$ is \textit{characteristic} in $M$ and $M$ is simply connected, Rochlin defined a quadratic form $q_\Sigma \colon H_1(\Sigma;\Z/2\Z) \to \Z/2\Z$ (see Definition \ref{Rochlinform}) which, in general, obstructs the smooth flexibility of $\Sigma$:

\[ \mathcal{E}^{\mathrm{DIFF}}(M,\Sigma) \leq \mathrm{Stab}_{\Mod(\Sigma)}[q_\Sigma].\]
It is well-known that this is in fact also an obstruction \textit{topologically}. This follows, for example, from work of Kirby-Taylor \cite{KirbyTaylor}. For the benefit of the reader, we recall some of the details.  

\subsection{Rochlin's quadratic form as an obstruction to topological flexibility}

We begin this subsection by recalling the relevant machinery from \cite{KirbyTaylor} (Theorem \ref{KirbyTaylor}). We then use this to observe that Rochlin's quadratic form (Definition \ref{Rochlinform}) is well-defined in the topological category (Proposition \ref{recoveringRochlin}), and moreover, obstructs extendable mapping classes (Corollary \ref{toppreservesroch}). 

In this and the next subsection, we work in the topological category. Let $\Sigma$ be a properly embedded, locally flat, connected, oriented surface in a compact, oriented topological 4-manifold $M$. Assume that $\Sigma$ is \textit{characteristic} in $M$, that is, the Poincar\'{e} dual of $[\Sigma,\partial \Sigma] \in H_2(M,\partial M;\mathbb{Z}/2\mathbb{Z})$ is the second Stiefel-Whitney class $w_2(M) \in H^2(M;\Z/2\Z)$. Note that, whilst $M$ may not have a linear tangent bundle, the Stiefel-Whitney classes of $M$ can still be defined and satisfy all of their usual properties (see, for example, \cite{friedl2024surveyfoundationsfourmanifoldtheory} for details).

\subsubsection{Kirby-Taylor quadratic forms}  \label{kirbytaylorsection}
We recall the definition of a topological spin structure. Let $\mathrm{STOP}(n)$ denote the topological group consisting of orientation preserving homeomorphisms of $\mathbb{R}^n$ that fix the origin. Recall that $\mathrm{Spin}(n)$ is the universal cover of $\mathrm{SO}(n)$ if $n\geq 3$. 

%By (the spin version of) \cite[Theorem 2.4]{KirbyTaylor}, there exists a spin structure $\mathfrak{s}$ on $M\backslash \Sigma$ that does not extend over $\Sigma$. 

\begin{defn} \label{topspinstructures}
Let $n\geq 3$. A topological spin structure on a topological $n$-manifold $N$ is a (homotopy class of) lift of the classifying map for the topological \textit{microbundle}\footnote{By Kister's Theorem \cite{Kister}, the topological microbundle $N \to N\times N \to N$ underlies a(n essentially unique) topological $\mathbb{R}^n$-bundle.} $\tau_N : N \to \mathrm{BTOP}(n)$ to $\mathrm{B TOP Spin(n)}$ where $\mathrm{TOPSpin(n)}$ denotes the universal cover of $\mathrm{STOP}(n)$. 
\end{defn}

If $N$ is smoothable and a smooth structure on $N$ is fixed, a topological spin structure on $N$ determines a unique spin structure on $N$. See, for example, \cite[Proof of Theorem 13.1]{Siebenmann}.

\begin{defn}
Let $\mathrm{Char}_{\mathrm{Spin}}(M,\Sigma)$ denote the set of \textit{topological spin structures} on $M\backslash \Sigma$ which do not extend over $\Sigma$. 
\end{defn}

\begin{proposition} \label{transitiveaction}
If $\Sigma$ is characteristic in $M$, the set $\mathrm{Char}_{\mathrm{Spin}}(M,\Sigma)$ is non-empty, and moreover, $H^1(M;\Z/2\Z)$ acts simply transitively on $\mathrm{Char}_{\mathrm{Spin}}(M,\Sigma)$. 
\end{proposition}

\begin{proof}
This is the spin case of \cite[Theorem 2.4]{KirbyTaylor}. Following the discussion in \cite[Section 9]{KirbyTaylor}, it holds equally well in the topological category. 
\end{proof}

\begin{defn}
A quadratic form on $H_1(\Sigma;\Z/2\Z)$ is a function $q\colon H_1(\Sigma;\Z/2\Z) \to \Z/2\Z$ satisfying \[q(x+y) = q(x) +q(y) + \langle x,y\rangle_2\] where $\langle -,-\rangle_2$ denotes the mod 2 intersection pairing on $H_1(\Sigma;\Z/2\Z)$. Let $\Theta_2(\Sigma)$ denote the set of quadratic forms on $H_1(\Sigma;\Z/2\Z)$.    
\end{defn}

A spin structure $\mathfrak{s} \in \mathrm{Char}_{\mathrm{Spin}}(M,\Sigma)$ determines the \textit{Kirby-Taylor quadratic form} \[q_{\mathfrak{s}} \colon H_1(\Sigma;\Z/2\Z) \to \Z/2\Z.\] Following the terminology of \cite{KirbyTaylor}, the assignment $\mathfrak{s} \to q_{\mathfrak{s}}$ defines a \textit{descent of structures} map 
 \[ \theta : \mathrm{Char}_{\mathrm{Spin}}(M,\Sigma) \to \Theta_2(\Sigma)\]
For the benefit of the reader, we recall a geometric definition of $q_{\mathfrak{s}}$.

Let $\nu_{\Sigma \subset M}$ be an oriented normal vector bundle to $\Sigma$ in $M$ \cite[Theorem 9.3A]{freedman}. We may equip the surface $\Sigma$ with a smooth structure. This (together with the standard smooth structure on $\mathbb{R}^2$) allows us to choose a smooth structure on the total space  $\nu_{\Sigma \subset M}$ turning it into a smooth vector bundle over $\Sigma$. The restriction of $\mathfrak{s}$ to $\nu_{\Sigma \subset M} \backslash \Sigma $ is a topological spin structure on a smooth 4-manifold. It determines a unique spin structure $\mathfrak{s} \in H^1(\mathrm{Fr}^+(\nu_{\Sigma \subset M} \backslash \Sigma);\Z/2\Z)$ where $\mathrm{Fr}^+(\nu_{\Sigma \subset M}\backslash\Sigma)$ denotes the oriented \textit{frame bundle} of $\nu_{\Sigma\subset M}\backslash \Sigma$.

Fix an oriented simple closed curve $\gamma \subset \Sigma$. Framing a push-off of $\gamma$ in $\nu_{\Sigma\subset M}\backslash \Sigma$ provides a loop in $\mathrm{Fr}^+(\nu_{\Sigma\subset M}\backslash \Sigma)$ and the (restricted) spin structure $\mathfrak{s}$ assigns an element of $\mathbb{Z}/2\mathbb{Z}$ to this loop. We now explain how the surface $\Sigma$ frames a push-off.

For convenience, write $D(\nu_{\Sigma \subset M})$ and $S(\nu_{\Sigma \subset M})$ for associated disk and unit sphere bundles in $\nu_{\Sigma\subset M}$. Choose a section $s : \gamma \to S(\nu_{\Sigma\subset M})\rvert_{\gamma}$. Let $\gamma' = s(\gamma)$. The $D^2$-fibers of the solid torus $D(\nu_{\Sigma\subset M})\rvert_{\gamma}$ bounding $S(\nu_{\Sigma\subset M})\rvert_{\gamma}$ are oriented, and their boundary circles have inherited an orientation through the inward-pointing normal convention. Let $\xi$ be a vector field along $\gamma'$ positively tangent to the circle fibers of $S(\nu_{\Sigma \subset M})\rvert_{\gamma}$. Let $n$ be an inward pointing vector field along $\gamma'$  for $S(\nu_{\Sigma \subset M})\rvert_{\gamma} = \partial D(\nu_{\Sigma \subset M})\rvert_{\gamma}$. Finally, consider the thickening of the solid torus $D(\nu_{\Sigma \subset M})\rvert_{\gamma}$ determined by the surface $\Sigma$, that is $D(\nu_{\Sigma \subset M})\rvert_{\nu_{\gamma\subset \Sigma}}$ where $\nu_{\gamma\subset \Sigma} \subset \Sigma$ is a normal bundle to $\gamma$ in $\Sigma$. The orientation of $\gamma$ and the orientation of $\Sigma$ determine a preferred normal vector field $v$ to $D(\nu_{\Sigma \subset M})\rvert_{\gamma} \subset D(\nu_{\Sigma \subset M})\rvert_{\nu_{\gamma\subset \Sigma}}$ along $\gamma'$. Define a loop $\Psi_{\gamma,\Sigma} : \gamma' \to \mathrm{Fr}^+(\nu_{\Sigma\subset M}\backslash \Sigma)$ by $\Psi_{\gamma,\Sigma} = (\overrightarrow{\gamma'},v,\xi,n)$. The quadratic form $q_{\mathfrak{s}}$ is defined by \[q_{\mathfrak{s}}(\gamma) = \mathfrak{s}(\Psi_{\gamma,\Sigma}(\gamma')) +1 \in \Z/2\Z.\]

\begin{rmk}
The normalization term appearing in the definition of the Kirby-Taylor form arises so that $q_{\mathfrak{s}}$ assigns the value zero to a loop bounding a disk on $\Sigma$. %We provide the reader with a more complete justification. Let $UT\Sigma \to \Sigma$ denote the unit tangent bundle of $\Sigma$. The construction actually produces a \textit{spin structure} $\xi_{\mathfrak{s}} \in H^1(UT\Sigma;\Z/2\Z)$ on the surface $\Sigma$, defined by the property that $\xi_{\mathfrak{s}}(\overrightarrow{\gamma}) = \Psi_\Sigma(\gamma')$ where $\overrightarrow{\gamma}$ denotes the forward pointing unit velocity vector field along $\gamma$. There is a canonical bijection between the set of spin structures on $\Sigma$ and the set $\Theta_2(\Sigma)$ \cite{Johnsonspin}. Under this bijection, the quadratic form determined by $\xi_{\mathfrak{s}}$ is defined by the equation $q_{\mathfrak{s}}(\gamma) = \xi_{\mathfrak{s}}(\overrightarrow{\gamma})+1 \in \Z/2\Z$.  
\end{rmk}

%The quadratic form $q_{\mathfrak{s}}$ is defined by $q_{\mathfrak{s}}(\gamma) = \mathfrak{s}(\Psi_\Sigma(\gamma')) +1 \in \Z/2\Z$. Let $\Theta_2(\Sigma)$ denote the set of quadratic forms on $H_1(\Sigma;\Z/2\Z)$. 

\begin{theorem}[Kirby-Taylor] \label{KirbyTaylor}
Let $\Sigma \subset M$ be a locally flat, properly embedded, oriented, characteristic surface in a compact oriented 4-manifold $M$. The descent of structures map  
\[ \theta : \mathrm{Char}_{\mathrm{Spin}}(M,\Sigma) \to  \Theta_2({\Sigma})\]
is well-defined and is equivariant with respect to the map $H^1(M;\Z/2\Z) \to H^1(\Sigma;\Z/2\Z)$ induced by inclusion. 
\end{theorem}

\begin{proof}
See \cite[Lemma 6.2]{KirbyTaylor} for the statement and the discussion in \cite[Section 9]{KirbyTaylor} which asserts that the construction applies in the topological category. 
\end{proof}

\begin{rmk}
For the benefit of the reader who wishes to consult Kirby-Taylor \cite{KirbyTaylor} for details, we explain how to extract Theorem \ref{KirbyTaylor} directly from \cite[Lemma 6.2]{KirbyTaylor}. Let $\mathrm{Char}(M,\Sigma)$ denote the set of \textit{$\mathrm{Pin}^{-}$-structures}\footnote{A $\mathrm{Pin}^-$-structure on a vector bundle $\xi$ is a spin structure on $\xi \oplus \det(\xi)$ where $\det(\xi)$ denotes the determinant line bundle of $\xi$.} on $M\backslash \Sigma$ that do not extend over $\Sigma$. Kirby-Taylor define a \textit{descent of structure} map $\theta : \mathrm{Char}(M,\Sigma)\to \mathrm{Pin}^{-}(\Sigma)$, where $\mathrm{Pin}^-(\Sigma)$ denotes the set of $\mathrm{Pin}^-$-structures on $\Sigma$. The set $\mathrm{Pin}^-(\Sigma)$ is canonically identified with the set of functions $q : H_1(\Sigma;\Z/2\Z) \to \Z/4\Z$ satisfying $q(x+y) = q(x)+q(y) + 2\bullet \langle x,y\rangle_2$ where $2\bullet$ denotes the inclusion $\mathbb{Z}/2\mathbb{Z} \cong 2\mathbb{Z}/4\mathbb{Z} \subset \Z/4\Z$. When $M$ and $\Sigma$ are oriented, these forms are only valued in $2\mathbb{Z}/4\mathbb{Z}$, so may be pulled back and seen as valued in $\Z/2\Z$. Moreover, the orientations provide a bijection between $\mathrm{Char}(M,\Sigma)$ and $\mathrm{Char}_{\mathrm{Spin}}(M,\Sigma)$. This is the perspective we take.  
\end{rmk}

Theorem \ref{KirbyTaylor} says that the choice of normal vector bundle $\nu_{\Sigma \subset M}$ is irrelevant when defining the Kirby-Taylor quadratic form. Essential to our story is the following naturality consequence. Recall that there is a left action of $\Mod(\Sigma)$ on the set of quadratic forms on $H_1(\Sigma;\mathbb{Z}/2\mathbb{Z})$: \[(f\cdot q)(x) = q(f^{-1}_{*}(x))\]for $f\in \Mod(\Sigma)$ and $x \in H_1(\Sigma;\Z/2\Z)$. We introduce some shorthand notation for the stabilizer of a quadratic form in $\Mod(\Sigma)$ (a \textit{2-spin mapping class group}, see Theorem \ref{johnsonspin} and Definition \ref{rspinmcg}). 

\begin{defn}
For a quadratic form $q : H_1(\Sigma;\mathbb{Z}/2\mathbb{Z}) \to \mathbb{Z}/2\mathbb{Z}$, let 
\[ \Mod(\Sigma)[q] = \mathrm{Stab}_{\Mod(\Sigma)}(q).\]
\end{defn}

\begin{proposition} \label{stabilizespinstructure}
% Assume that $H^1(M;\mathbb{Z}/2\mathbb{Z}) = 0$. Then there is a unique spin structure $\mathfrak{s} \in \mathrm{Char}_{\mathrm{Spin}(M,\Sigma)}$ on $M-\Sigma$ that does not extend over $\Sigma$. Let $q_{\mathfrak{s}} ; H_1(\Sigma;\Z/2\Z) \to \Z/2\Z$ be the associated Kirby-Taylor quadratic form. 
Let $f\colon(M,\Sigma) \to (M,\Sigma)$ be a homeomorphism of the pair and fix $\mathfrak{s} \in \mathrm{Char}_{\mathrm{Spin}}(M,\Sigma)$. Then 
\[ q_{f^*\mathfrak{s}} = f^{-1} \cdot q_{\mathfrak{s}}. \]
If $H^1(M;\Z/2\Z) = 0$, there is a unique spin structure $\mathfrak{s} \in \mathrm{Char}_{\mathrm{Spin}}(M,\Sigma)$, and therefore  \[\mathcal{E}^{\mathrm{TOP}}(M,\Sigma) \leq \Mod(\Sigma)[q_{\mathfrak{s}}].\]  %If the pair $(M,\Sigma)$ is smooth, then Rochlin's quadratic form is defined $q_\Sigma : H_1(\Sigma;\Z$
\end{proposition}

\begin{proof}
By the uniqueness of normal vector bundles up to isotopy \cite[Theorem 9.3D]{freedman}, we may assume that $f$ fixes $E:= \nu_{\Sigma \subset M}$ setwise and $f\rvert_{E} : E \to E$ is an isomorphism of vector bundles covering $f\rvert_{\Sigma}$. By pulling back a smooth structure on $E$ along $f$, we may assume that $f\rvert_{E}$ is actually a fiber preserving diffeomorphism of smooth vector bundles, which, for clarity, we write as $f : E' \to E$. 

Restricting the spin structure $\mathfrak{s}$  to $E\backslash \Sigma$ determines a cohomology class $\rho_\mathfrak{s} \in H^1(\mathrm{Fr}^+(E\backslash \Sigma);\Z/2\Z)$. Restricting $f^*\mathfrak{s}$ to $E'\backslash \Sigma$ is the same as pulling back the restriction of $\mathfrak{s}$ to $E\backslash \Sigma$ under $f\rvert_{E'} : E' \to E$. It follows that the cohomology class $\rho_{f^*\mathfrak{s}} \in H^1(\mathrm{Fr}^+(E'\backslash \Sigma);\Z/2\Z)$ determined by restricting $f^*\mathfrak{s}$ to $E'\backslash \Sigma$ is the pullback of $\rho_{\mathfrak{s}}$ under the induced map on frame bundles $\varphi_{f} : \mathrm{Fr}^+(E'\backslash \Sigma) \to \mathrm{Fr}^+(E\backslash \Sigma)$:
\[ \varphi_f^*(\rho_{\mathfrak{s}}) = \rho_{f^*\mathfrak{s}}\]
Fix an oriented simple closed curve $\gamma \subset \Sigma$. The conclusion $q_{f^*\mathfrak{s}}(\gamma) = q_{\mathfrak{s}}(f(\gamma))$ now follows from the fact that the surface framing on a push-off $\gamma' \subset E'\backslash \Sigma$ is carried to the surface framing of $f(\gamma') \subset E \backslash \Sigma$ by $\varphi_f$, and that $f(\gamma')$ is a push-off of $f(\gamma) \subset \Sigma$. The final assertion is now an immediate consequence of Proposition \ref{transitiveaction}. 
\end{proof}

\subsubsection{Rochlin's quadratic form and spin-flexibility} \label{secrochlin} Assume now that $M$ is simply connected and let $\mathfrak{s} \in \mathrm{Char}_{\mathrm{Spin}}(M,\Sigma)$ denote the unique spin structure on $M-\Sigma$ that does not extend over $\Sigma$. Kirby-Taylor showed that $q_{\mathfrak{s}} \colon H_1(\Sigma;\Z/2\Z) \to \Z/2\Z$ recovers  \textit{Rochlin's quadratic form}. We record their calculation here after introducing some machinery. Recall that $\nu_{\Sigma \subset M}$ denotes a normal vector bundle to $\Sigma$ in $M$. 

\begin{defn} \label{membranedef}
A \textit{membrane} for the pair $(M,\Sigma)$ is a connected, orientable, compact, locally flat surface $S \subset M$ with the following properties:
\begin{enumerate}
    \item[1)]$\partial S \subset \Sigma$ is a disjoint collection of simple closed curves on $\Sigma$.
    \item[2)] $S$ is transverse to $\Sigma$ on its interior. 
    \item[3)] $S$ is \textit{transverse to $\Sigma$ along its boundary} $\partial S$: there is a trivialization \[D(\nu_{\Sigma \subset M})\rvert_{\nu_{\partial S \subset \Sigma}} \cong \partial S \times \mathbb{R} \times D^2\] under which $S \cap D(\nu_{\Sigma \subset M})\rvert_{\nu_{\partial S \subset \Sigma}}$ is given by $\partial S \times \{0\} \times J$ where $J = \{(x,0)\,|\, 0\leq x\leq 1\} \subset D^2$.    
\end{enumerate}
\end{defn}

%We provide a brief outline as to why membranes exist. Fix a simple closed curve $\gamma \subset \Sigma$. To ensure transversality near the boundary, choose first the collar $C$ of the membrane in $\nu_{\Sigma \subset M}$ by choosing a push-off $\gamma' \subset S(\nu_{\Sigma\subset M})\rvert_{\gamma}$. A null-homotopy of $\gamma'$ in $M$ allows us to define a continuous map $f : D^2 \to M$ with boundary collar $C$. By a homotopy of $f$ relative to this collar, we may assume that $f$ is a generic immersion, with double points disjoint from $\Sigma$ (see \cite{friedl2024surveyfoundationsfourmanifoldtheory}). By resolving the double points, we may replace our immersed disk with some embedded, orientable, connected surface $S$. By topological transversality, we may assume that $S$ is transverse to $\Sigma$ on its interior. 

\begin{rmk}The existence of a membrane $S$ for $(M,\Sigma)$ with prescribed boundary can be arranged as follows. First, build the relevant collar in $\nu_{\Sigma\subset M}$. Now invoke simple connectivity of $M$ and approximation theorems\footnote{See \cite[Corollary 9.5C]{freedman}. We also recommend the proof of \cite[Theorem 10.6]{friedl2024surveyfoundationsfourmanifoldtheory}, which illustrates the basic tools.} to cap this off with a \textit{generically topologically immersed} disk (care must be taken so that the boundary collar is unchanged when arranging this). Resolving the double points of the immersed disk yields an embedded surface $S \subset M$. We make $S$ transverse to $\Sigma$ relative to its boundary collar (applying \textit{topological transversality}, see \cite[Theorem 9.5A]{freedman}). 
\end{rmk}

Fix a membrane $S$ for the pair $(M,\Sigma)$. Consider the trivialization of $\nu_{\Sigma \subset M}\rvert_{\nu_{\partial S \subset \Sigma}} \cong \partial S \times \mathbb{R} \times \mathbb{R}^2$ supplied by Definition \ref{membranedef}. View the total space as a 3-plane bundle over $\partial S$. Let $V$ denote the 2-plane subbundle of this 3-plane bundle obtained by splitting off the inward-pointing collar direction for $\partial S \subset S$.

%We now define the \textit{relative Euler number} $e_\Sigma(S) \in \mathbb{Z}$ of a membrane $S$ for the pair $(M,\Sigma)$.

\begin{defn}
A normal vector bundle $\nu_{S \subset M}$ for $S$ in $M$ is \textit{$\Sigma$-compatible} if $\nu_{S\subset M}\rvert_{\partial S} = V$. In particular, $\nu_{\partial S \subset \Sigma} \subset \nu_{S\subset M}\rvert_{\partial S}$. 
\end{defn}

\begin{proposition}
If $S$ is a membrane for the pair $(M,\Sigma)$, it has a $\Sigma$-compatible normal vector bundle. Moreover, any two $\Sigma$-compatible normal vector bundles are isotopic through $\Sigma$-compatible normal vector bundles.  
\end{proposition}

\begin{proof}
By \cite[Theorem 9.3A]{freedman}, extend $V$ to a normal vector bundle of $S$ in $M$, denoted $\nu_{S \subset M}$. By construction, \(\nu_{\partial S \subset \Sigma} \subset \nu_{S \subset M}\rvert_{\partial S} =V.\) Uniqueness follows from \cite[Theorem 9.3D]{freedman} relative to (a neighbourhood of) $\partial S$.  
\end{proof}

The bundle $\nu_{\partial S\subset \Sigma}$ is a trivial line bundle. An orientation of $\nu_{\partial S \subset \Sigma}$ makes a choice of (homotopy class of) nowhere vanishing section canonical. An orientation of $\partial S$ determines (using  the orientation of $\Sigma$) an orientation of $\nu_{\partial S\subset \Sigma}$, and so a distinguished (homotopy class of) section  \[s_\Sigma : \partial S \to \nu_{\partial S\subset \Sigma} \subset \nu_{S \subset M}\rvert_{\partial S}.\]

\begin{defn}
Let $S \subset M$ be an oriented membrane for the pair $(M,\Sigma)$ and $\nu_{S \subset M}$ a $\Sigma$-compatible normal vector bundle to $S$ in $M$. The \textit{relative Euler class}  $e_\Sigma(S)$ is the obstruction  $e(\nu_{S \subset M},s_\Sigma) \in H^2(S,\partial S;\Z)$ to extending the section $s_\Sigma : \partial S \to \nu_{S\subset M}$ to a nowhere vanishing section $S \to \nu_{S \subset M}$. The \textit{relative Euler number} is $e_\Sigma(S) = e(\nu_{S\subset M},s_\Sigma) \cap [S,\partial S] \in H_0(S)  \cong \mathbb{Z}$. 
\end{defn}

\begin{rmk}
We recommend the account \cite[Section 2]{conway2024unknottingnonorientablesurfaces} for a discussion of the relative Euler number. In practice, we let $s' : S \to \nu_{S \subset M}$ be a \textit{generic extension} of $s_\Sigma$. The \textit{relative Euler number}  $e_\Sigma(S)$ is the signed count of the zeros of $s'$.  
\end{rmk}

\begin{defn}[Rochlin's quadratic form] \label{Rochlinform}
Let $\Sigma \subset M$ be characteristic\footnote{If $\Sigma$ is not characteristic, tubing a membrane to a surface witnessing the failure of $\Sigma$ to be characteristic shows that $q_\Sigma$ is not well-defined.} and $M$ simply connected. Let $\gamma \subset \Sigma$ be a simple closed curve and $S \subset M$ a membrane for the pair $(M,\Sigma)$ with boundary $\gamma$. Define a function $q_\Sigma \colon H_1(\Sigma;Z/2\Z) \to \Z/2\Z$ by 
\[q_\Sigma(\gamma) = |\mathrm{Int}(S) \cap \Sigma| + e_\Sigma(S) \mod 2.\]
\end{defn}

Kirby-Taylor identify $q_\Sigma$ with the quadratic form $q_{\mathfrak{s}}$ for $\mathfrak{s} \in \mathrm{Char}_{\mathrm{Spin}}(M,\Sigma)$.  

\begin{proposition}[Recovering Rochlin's quadratic form] \label{recoveringRochlin}
Assume that $M$ is simply connected and $\Sigma \subset M$ characteristic. Let $\mathfrak{s}$ denote the unique spin structure on $M-\Sigma$ that does not extend over $\Sigma$, and $q_{\mathfrak{s}} \colon H_1(\Sigma;\Z/2\Z) \to \Z/2\Z$ the associated Kirby-Taylor quadratic form. Let $\gamma \subset \Sigma$ be a simple closed curve. Choose a membrane $S \subset M$ for the pair $(M,\Sigma)$ with $\partial S = \gamma$. Then 
\[ q_{\mathfrak{s}}(\gamma)  = |\mathrm{Int}(S) \cap \Sigma| + e_\Sigma(S) = q_\Sigma(\gamma)\mod 2.\]
\end{proposition}

\begin{proof}
We sketch the argument appearing in \cite[Discussion after Theorem 6.3]{KirbyTaylor}. Choose a $\Sigma$-compatible normal vector bundle $\nu_{S \subset M}$ to $S$ in $M$. Choose a smoothing of $\nu_{S\subset M}$ turning it into a smooth vector bundle over the surface $S$. Let $s' \colon S \to \nu_{S\subset M}$ be a generic extension of the nowhere vanishing section $s \colon \partial S \to \nu_{S \subset M}$ determined by $\Sigma$. Let $T \subset S$ be obtained by removing a small boundary collar of $S$ and small disks about $\mathrm{Int}(S) \cap \Sigma$ and $s'^{-1}(0)$. 

Let $N_T =\nu_{S\subset M}\rvert_T \subset M-\Sigma$. The topological spin structure $\mathfrak{s}$ on $M-\Sigma$ determines a spin structure on $N_T$. The nowhere vanishing section $s' : T \to N_T$ determines a trivialization of $N_T$. This allows us to restrict the spin structure $\mathfrak{s}$ on $N_T$ to a spin structure $\eta$ on the surface $T$. Spin structures satisfy \textit{homological coherence}  of Definition \ref{windingnumberfunction}(b) and Theorem \ref{bijectionbetweenwindingandspin}). If $g(T)$ is the genus of $T$ and $b(T)$ its number of boundary components,
\[ \eta(\partial T) = \chi(T) = 2-2g(T) -b(T) =  |\mathrm{Int}(S)\cap \Sigma| + e_\Sigma(S) +1 \mod2 .\] Kirby-Taylor argue that the spin structure $\eta$ does \textit{not} extend over the components of $\partial T$ coming from $S\cap \Sigma$ and the zeros of $s'$. In particular, it assigns the value zero to these components\footnote{If it did extend, it must assign the value one. This can be seen through the homological coherence property.}. Finally, the boundary component of $T$ coming from the collar of $S$ is a push-off $\gamma'$ of $\gamma$ (relative to $\Sigma)$ and tracing through the definitions shows that $\eta(\overrightarrow{\gamma'}) +1 = q_{\mathfrak{s}}(\gamma)$. 
\end{proof}

\begin{corollary}
 \label{toppreservesroch}
Fix a locally flat, properly embedded, characteristic surface $\Sigma$ in a compact, simply connected topological 4-manifold $M$. Then
\[ \mathcal{E}^{\mathrm{TOP}}(M,\Sigma) \leq \Mod(\Sigma)[q_\Sigma]\]
where $q_\Sigma : H_1(\Sigma;\mathbb{Z}/2\mathbb{Z}) \to \mathbb{Z}/2\mathbb{Z}$ denotes Rochlin's quadratic form. 
\end{corollary}

\begin{proof}
This is Proposition \ref{stabilizespinstructure} and Proposition \ref{recoveringRochlin}. 
\end{proof}

%\begin{proof}
%Fix a homeomorphism of the pair $f : (M,\Sigma) \to (M,\Sigma)$. Fix a simple closed curve $\gamma \subset \Sigma$ and a surface $S \subset M$ with $\partial S = \gamma$ and interior transverse to $\Sigma$. 
%\end{proof}

\begin{defn}
A properly embedded characteristic surface $\Sigma \subset M$ is \textit{topologically spin-flexible} in $M$ if \(\mathcal{E}^{\mathrm{TOP}}(M,\Sigma) = \Mod(\Sigma)[q_\Sigma]\). If $M$ is smooth and $\Sigma$ smoothly embedded, we say $\Sigma$ is \textit{smoothly spin-flexible} if $\mathcal{E}^{\mathrm{DIFF}}(M,\Sigma) = \Mod(\Sigma)[q_\Sigma]$. 
\end{defn}

\begin{rmk}
We conclude this subsection by noting that Proposition \ref{stabilizespinstructure} is all we require to establish spin-flexibility of the characteristic surfaces considered in this paper. We consider smooth pairs $(M,\Sigma)$ with $M$ simply connected, in which case \textit{Rochlin's form} $q_\Sigma \colon H_1(\Sigma;\Z/2\Z) \to \Z/2\Z$ is classically defined. We establish the lower bound 
\[\Mod(\Sigma)[q_\Sigma] \leq  \mathcal{E}^{\mathrm{TOP}}(M,\Sigma).\] The containment $\Mod(\Sigma)[q_\Sigma] \leq \Mod(\Sigma)[q_{\mathfrak{s}}]$ guaranteed by Proposition \ref{stabilizespinstructure} implies, for the examples that we consider, that  $q_\Sigma = q_{\mathfrak{s}}$. 
\end{rmk}

\subsection{Some extendable mapping classes} In the sequel, we use the following methods for extending mapping classes. Montesinos proved that an unknotted torus $\Sigma_1 \subset S^4$ is smoothly spin-flexible \cite{Montesinos}. Hirose proved that unknots $\Sigma_g \subset S^4$ of all genera are smoothly spin-flexible \cite{HiroseS4}. These statements have rel boundary versions. Let $\Sigma_g^1 \subset D^4$ be the result of removing a 4-ball from the pair $(S^4,\Sigma_g)$ meeting $\Sigma_g$ in a closed 2-disk. 
\begin{proposition} \label{spinflexunknot}
For $g\geq 1$, $\mathcal{E}^{\mathrm{DIFF}}(D^4,\Sigma_g^1) = \Mod(\Sigma_g^1)[q_{\Sigma_g^1}]$. 
\end{proposition} 
\begin{proof}
This follows from the statements in the closed case. For details we refer the reader to \cite[Proposition 2.15]{lehman}. 
\end{proof}

For a simple closed curve $\gamma \subset \Sigma$, let $T_{\gamma}$ be the (left-handed) Dehn twist about $\gamma$.

\begin{proposition}[Hopf band trick] \label{HBtrick}
Let $\gamma \subset \Sigma \subset M$ be a simple closed curve. Assume that $\gamma$ bounds a topologically (resp. smoothly) embedded disk $D\subset \mathrm{int}(M)$ with interior disjoint from $\Sigma$ and $e_\Sigma(D) = \pm1$. Then $T_{\gamma} \in \mathcal{E}^{\mathrm{TOP}}(M,\Sigma)$ (resp. $T_{\gamma} \in \mathcal{E}^{\mathrm{DIFF}}(M,\Sigma)$). 
\end{proposition}

\begin{proof}
See the proof of \cite[Lemma 6.2]{Slopegenera}. 
\end{proof}

\begin{rmk}
In the setting of Proposition \ref{HBtrick}, let $\nu_{\gamma \subset \Sigma} \cong S^1 \times [0,1]$ be a closed tubular neighbourhood of $\gamma$ in $\Sigma$. Let $\nu_D \cong D \times D^2$ be a closed tubular neighbourhood of $D$ in $M$. Up to isotopy, $\nu_{\gamma \subset \Sigma}$ is a \textit{Hopf band} in $S^3 = \partial (D\times D^2)$, that is, an annulus with $\pm 1$ full twist (computed as the linking number of the boundary components in $S^3$) and unknotted core. 
\end{rmk}

\section{Background on $r$-spin mapping class groups}

We recall fundamental results about stabilizers of $r$-spin structures on closed surfaces ($r$-spin mapping class groups). We conclude this section by exploiting how these stabilizers intersect the Torelli subgroup if $r > 2$. This section follows \cite{SalterInvent} and \cite{CalderonSalter} closely.

\subsection{Higher spin structures and winding number functions} \label{windingbackground} Let $\Sigma$ be a closed, oriented surface. Let $UT\Sigma \to \Sigma$ denote the unit tangent bundle of $\Sigma$ and $\langle -,-\rangle$ the oriented intersection pairing on $H_1(\Sigma;\mathbb{Z})$. 

\begin{defn}
An \textit{$r$-spin structure} on $\Sigma$ is a cohomology class $\xi \in H^1(UT\Sigma;\mathbb{Z}/r\mathbb{Z})$ where $\xi(\mathbf{z}) =1$ for $\mathbf{z}$ a positively oriented fiber of $UT\Sigma \to \Sigma$.
\end{defn}

Let $\mathcal{S}(\Sigma)$ denote the set of oriented isotopy classes of simple closed curves on $\Sigma$. 

\begin{defn} \label{windingnumberfunction}
A \textit{$\mathbb{Z}/r\mathbb{Z}$-winding number function} on $\Sigma$ is a function $\psi : \mathcal{S}(\Sigma) \to \mathbb{Z}/r\mathbb{Z}$ satisfying 
\begin{enumerate}
    \item (\textit{Twist-linearity}) For any $a,b \in \mathcal{S}(\Sigma)$, 
    $$ \psi(T_a(b)) = \psi(b) + \langle b,a\rangle \psi(a) \mod r.$$
    %where $\langle -,-\rangle : H_1(\Sigma;\mathbb{Z}) \times H_1(\Sigma;\mathbb{Z}) \to \mathbb{Z}$ denotes the oriented intersection pairing. 
    \item (\textit{Homological coherence}) Let $S \subset \Sigma$ be a subsurface with boundary components $c_1,\ldots,c_k$ oriented with $S$ to their left. Then 
    $$\sum_{i=1}^k \psi(c_i) = \chi(S) \mod r$$
    where $\chi(S)$ denotes the Euler characteristic of $S$. 
\end{enumerate}

\begin{rmk}
By homological coherence, $r\,| \chi(\Sigma)$. If $\overline{a}$ denotes $a \in \mathcal{S}(\Sigma)$ with the opposite orientation, then $\psi(\overline{a}) = -\psi(a) \mod r$. 
\end{rmk}

\end{defn}

\begin{eg}[Constructing winding number functions]
\label{windingnumberfromrspin}
 Fix an $r$-spin structure $\xi \in H^1(UT\Sigma;\mathbb{Z}/r\mathbb{Z})$. Given $a \in \mathcal{S}(\Sigma)$, let $\overrightarrow{a} \subset UT\Sigma$ denote the unit forward pointing velocity vector field of $a$. The assignment $a \to \overrightarrow{a}$ is referred to as the \textit{Johnson lift} \cite{Johnsonspin}. Define a function $\phi_\xi : \mathcal{S}(\Sigma) \to \mathbb{Z}/r\mathbb{Z}$ by $\phi_\xi(a) = \xi(\overrightarrow{a})$. 
\end{eg}

\begin{theorem} \label{bijectionbetweenwindingandspin}
The assignment \[\xi \in H^1(UT\Sigma;\mathbb{Z}/r\mathbb{Z})\to \phi_{\xi} : \mathcal{S}(\Sigma) \to \mathbb{Z}/r\mathbb{Z}\] (Example \ref{windingnumberfromrspin}) defines a bijective correspondence between $r$-spin structures on $\Sigma$ and mod $r$ winding number functions on $\Sigma$. In particular, there are $r^{2g}$ winding number functions on $\Sigma$ valued in $\mathbb{Z}/r\mathbb{Z}$. 
\end{theorem}

\begin{proof}
See \cite{HJ} and e.g. \cite[Sections 2, 3 and Theorem 3.2]{salterquintic}. 
\end{proof}

%\begin{corollary}
%A winding number function $\psi : \mathcal{S}(\Sigma) \to \mathbb{Z}/r\mathbb{Z}$ is determined by its values on a collection of simple closed curves whose homology classes generate $H_1(\Sigma;\mathbb{Z})$. 
%\end{corollary}

\subsection{$r$-spin mapping class groups} The mapping class group $\Mod(\Sigma)$ acts on the set of $\mathbb{Z}/r\mathbb{Z}$ valued winding number functions by pullback. Precisely
\[(f\cdot \psi)(a) = \psi(f^{-1}(a))\] for $f \in \Mod(\Sigma)$ and $a \in \mathcal{S}(\Sigma)$.

\begin{defn} \label{rspinmcg}
Let $\psi : \mathcal{S}(\Sigma) \to \mathbb{Z}/r\mathbb{Z}$ be a winding number function on $\Sigma$. The stabilizer of $\psi$ in $\Mod(\Sigma)$ is denoted
\[ \Mod(\Sigma)[\psi] = \{ f \in \Mod(\Sigma)\,|\, f\cdot \psi = \psi\}.\]
We call $\Mod(\Sigma)[\psi]$ an \textit{$r$-spin mapping class group}. 
\end{defn}

\begin{defn}
A simple closed curve $a \subset \Sigma$ is \textit{admissible} under a winding number function $\psi : \mathcal{S}(\Sigma) \to \mathbb{Z}/r\mathbb{Z}$ if it is nonseparating and $\psi(a) = 0$.  
\end{defn}

By twist-linearity, if $a \subset \Sigma$ is $\psi$-admissible, $T_a \in \Mod(\Sigma)[\psi]$. Calderon-Salter \cite{CalderonSalter} established that $\psi$-admissible twists are enough to generate the $r$-spin mapping class group $\Mod(\Sigma)[\psi]$.

\begin{theorem}[{\cite[Proposition 5.1]{CalderonSalter}}] \label{generatingrspinmcg}
Let $\psi : \mathcal{S}(\Sigma_g) \to \mathbb{Z}/r\mathbb{Z}$ be a winding number function on a closed surface $\Sigma_g$ of genus $g \geq 3$. Then
\[ \Mod(\Sigma_g)[\psi] = \langle T_{a}\,|\, a \subset \Sigma_g, a \,\mathrm{ admissible}\rangle.\]
\end{theorem}

We require two instances of the $r$-spin \textit{change of coordinates} principle. We refer the reader to \cite[Section 5.2]{CalderonSalter} for more details. Let $\Sigma_g$ be a closed oriented surface of genus $g$ and fix a winding number function $\psi : \mathcal{S}(\Sigma_g) \to \mathbb{Z}/r\mathbb{Z}$. 

\begin{proposition} \label{rspinchangeofcoordinates}
Let $a,b \subset \Sigma_g$ be nonseparating simple closed curves on $\Sigma_g$, $g \geq 2$. Assume that $\psi(a) = \psi(b)$. Then there exists $f \in \Mod(\Sigma_g)[\psi]$ such that $f(a) = b$.  
\end{proposition}

\begin{proof}
See \cite[Lemma 5.8]{CalderonSalter}.
\end{proof}

\begin{proposition} \label{admissibledual}
Let $a \subset \Sigma_g$ be a nonseparating simple closed curve and $g \geq 2$. Then there exists a $\psi$-admissible curve $b \subset \Sigma_g$ with geometric intersection number $i(a,b) = 1$. 
\end{proposition}

\begin{proof}
If $r$ is even, this follows from \cite[Lemma 5.7]{CalderonSalter}, and if $r$ is odd, this follows from \cite[Lemma 5.8]{CalderonSalter}. 
\end{proof}

\subsection{The Torelli group and $r$-spin mapping class groups} Let $\Psi : \Mod(\Sigma) \to \mathrm{Sp}(H_1(\Sigma;\mathbb{Z}))$ be the \textit{symplectic representation} of $\Mod(\Sigma)$. The kernel $\mathcal{I}(\Sigma)$ of $\Psi$ is the \textit{Torelli subgroup} of $\Mod(\Sigma)$. Recall that two disjoint simple closed curves $a,b \subset \Sigma$ are a genus $h$ \textit{bounding pair} if they cobound some genus $h$ subsurface on $\Sigma$. A bounding pair is nonseparating if both curves are nonseparating. The \textit{bounding pair map} is $T_aT_b^{-1}$ and, as $a$ and $b$ are homologous, $T_aT_b^{-1} \in \mathcal{I}(\Sigma)$. The following result is an essential ingredient in proving Theorem \ref{maintheorem}. It exploits how $r$-spin mapping class groups intersect Torelli. 

\begin{lemma} \label{maintechnicallemma}
Let $\psi\colon\mathcal{S}(\Sigma) \to \mathbb{Z}/r\mathbb{Z}$ be a winding number function on $\Sigma$. Assume that $\Sigma$ has genus $\geq 3$. If $r$ is even, let $\overline{\psi}$ denote the mod 2 reduction of $\psi$. Let $T_{a}T_{b}^{-1} \in \Mod(\Sigma)$ be a nonseparating genus one bounding pair map. Then 
\[ \langle \Mod(\Sigma)[\psi],T_{a}T_{b}^{-1}  \rangle =  \begin{cases}  \Mod(\Sigma) & \mathrm{if  } \;r\; \mathrm{ is \; odd, } \\ \Mod(\Sigma)[\overline{\psi}] & \mathrm{ if }\; r \mathrm{ \;is\; even}. \end{cases} \]
\end{lemma}

\begin{proof}
Let $\Gamma = \langle \Mod(\Sigma)[\psi],T_{a}T_{b}^{-1}  \rangle$. Orient $a$ and $b$ so that they bound a surface of genus one to their left. Consider a simple closed curve $c \subset \Sigma$. By homological coherence, $\psi(b) +\psi(a) = -2$. Compute $\psi(T_a T_{b}^{-1}(c))$ using twist-linearity and homological coherence,
\begin{align*}
\psi(T_a T_{b}^{-1}(c)) &= \psi(c) - \langle c,b \rangle \psi(b) + \langle c,a\rangle \psi(a)  \mod r\\
&= \psi(c) + \langle c,a \rangle(\psi(b)+\psi(a)) \mod r \\
&= \psi(c) - 2\langle c,a \rangle \mod r.
\end{align*}
Fix $c\subset \Sigma$ for which $\langle c,a \rangle = 1$ and $\psi(c) = 0$ (Proposition \ref{admissibledual}, note that here we are using the fact that $a$ is nonseparating). Then $\psi(T_{a}T_b^{-1}(c)) = -2$.  Let $d = T_{a}T_{b}^{-1}(c)$. As $T_c \in \Mod(\Sigma)[\psi]$, it follows that $T_d \in \Gamma$. Let $e \subset \Sigma$ be a simple closed curve with $\langle e,d \rangle = 1$  and $\psi(e) = 0$ (Proposition \ref{admissibledual}). Then \[ \psi(T_{d}^{- k}(e)) = 2k \]
and $T_{T_{d}^{-k}(e)} \in \Gamma$. Assume $r$ is even. Using the family of curves $T^{-k}_d(e)$ and Proposition \ref{rspinchangeofcoordinates}, $\Gamma$ contains the Dehn twist on any nonseparating simple closed curve $\gamma \subset \Sigma$ with $\psi(\gamma) = 0 \mod 2$. By  Theorem \ref{generatingrspinmcg}, 
\[ \Mod(\Sigma)[\overline{\psi}] \leq \Gamma.\]
We see directly that $\Gamma \leq \Mod(\Sigma)[\overline{\psi}]$, and so equality follows. 

Assume $r$ is odd. We have shown that there exists a Dehn twist on a nonseparating curve $\gamma \subset \Sigma$ with $\psi(\gamma) = 1 \mod r$ contained in $\Gamma$. Use again a geometrically dual admissible curve to show that $\Gamma$ contains Dehn twists on nonseparating simple closed curves with any prescribed winding number. By Proposition \ref{rspinchangeofcoordinates}, $\Gamma$ contains the Dehn twist on any nonseparating simple closed curve, and hence $\Gamma = \Mod(\Sigma)$.  
\end{proof}

\subsection{Explicit generating sets for $r$-spin mapping class groups}

We record a result which gives explicit finite generating sets for $r$-spin mapping class groups. It is used in proving Theorem \ref{nonflex}. A \textit{network} $\mathcal{C}$ on $\Sigma$ is a collection of simple closed curves on $\Sigma$ such that any two curves in $\mathcal{C}$ intersect at most once, and there are no triple intersections. We say $\mathcal{C}$ is \textit{filling} if $\Sigma \backslash \bigcup_{c \in \mathcal{C}} c$ is a collection of disks. 

\begin{defn} \label{intersectiongraphofnetwork}
The \textit{intersection graph} $\Lambda_{\mathcal{C}}$ of a network $\mathcal{C}$ on $\Sigma$ is the graph with vertex set $\mathcal{C}$ and an edge between two vertices if they intersect once. 
\end{defn}

\begin{proposition}[{\cite[Corollary 3.11]{strata3}}]
\label{gensetrspin}
Let $\mathcal{C}$ be a filling network of curves on $\Sigma_g$ for $g\geq 5$. Suppose that the intersection graph $\Lambda_C$ is a tree which contains the $E_6$ Dynkin diagram as a subgraph and that $\mathcal{C}$ cuts the surface into $n$ polygons with $4(k_1+1),\ldots,4(k_n+1)$ sides. Set $r = \gcd(k_1,\ldots,k_n)$. Then there is a unique mod $r$ winding number function $\psi : \mathcal{S}(\Sigma_g) \to \mathbb{Z}/r\mathbb{Z}$ for which 
\[ \mathcal{T}_C = \langle T_c\,:\, c \in \mathcal{C}\rangle  =  \Mod(\Sigma_g)[\psi].\]
\end{proposition}

\subsection{Quadratic forms and mod 2 winding number functions} A result of Johnson  together with Theorem \ref{bijectionbetweenwindingandspin} allows us to identify the stabilizer of a quadratic form $q : H_1(\Sigma;\mathbb{Z}/2\mathbb{Z}) \to \mathbb{Z}/2\mathbb{Z}$ in $\Mod(\Sigma)$ with a 2-spin mapping class group. We may then apply results from Subsection \ref{windingbackground}.

\begin{theorem}[] 
A winding number function $\psi : \mathcal{S}(\Sigma) \to \mathbb{Z}/2\mathbb{Z}$ defines a quadratic form $q_{\psi} : H_1(\Sigma;\mathbb{Z}/2\mathbb{Z}) \to \mathbb{Z}/2\mathbb{Z}$ by the requirement that 
\[q_{\psi}(\gamma) = \psi(\gamma) + 1\mod 2\]
for $\gamma \in \mathcal{S}(\Sigma)$. Conversely, a quadratic form on $H_1(\Sigma;\mathbb{Z}/2\mathbb{Z})$ defines a $\mathbb{Z}/2\mathbb{Z}$ winding number function by the same formula. %Moreover, quadratic forms on $H_1(\Sigma;\mathbb{Z}/2\mathbb{Z})$ are determined by their values on a generating set for $H_1(\Sigma;\mathbb{Z}/2\mathbb{Z})$.  
 \label{johnsonspin}
\end{theorem}
\begin{proof}
This follows from \cite[Theorem 3A]{Johnsonspin} in combination with Theorem \ref{bijectionbetweenwindingandspin}. 
\end{proof}

In particular, we have a bijective correspondence between mod 2 winding number functions on $\Sigma$ and quadratic forms on $H_1(\Sigma;\mathbb{Z}/2\mathbb{Z})$. We use the following convention when dealing with admissible curves. For a form $q$, %a quadratic form $q : H_1(\Sigma;\mathbb{Z}/2\mathbb{Z}) \to \mathbb{Z}/2\mathbb{Z}$ defines a mod 2 winding number function $\psi_q : \mathcal{S}(\Sigma) \to \mathbb{Z}/2\mathbb{Z}$ by \[\psi_q(a) = q([a])+1.\] 
 a curve $a \subset \Sigma$ is \textit{admissible} if it is admissible under the corresponding 2-spin structure, in particular, if $q([a]) = 1$. 

\section{Extending homeomorphisms} %using techniques of Conway-Powell}

 Let $\Sigma \subset X$ be a closed surface in a closed 4-manifold $X$. Then $\Sigma$ admits a closed tubular neighbourhood $\nu \Sigma \subset X$ \cite[Theorem 9.3A]{freedman}. A homeomorphism of the pair $h : (X,\Sigma)  \to (X,\Sigma)$ induces a homeomorphism of the exterior $X_\Sigma = X-\nu \Sigma$ (by uniqueness of tubular neighbourhoods up to isotopy \cite[Theorem 9.3D]{freedman}). Lifting to the universal cover of $X_\Sigma$, we obtain a (semilinear) isometry of the \textit{equivariant intersection form} on the surface exterior $X_\Sigma$. Conway-Powell \cite{Conway_Powell} establish a converse for \textit{$\mathbb{Z}$-surfaces}. If $\pi_1 X_\Sigma \cong \mathbb{Z}$, an isometry of the equivariant intersection form on the surface exterior can be promoted to a homeomorphism of the pair $(X,\Sigma)$. We recall the relevant machinery from \cite{Conway_Powell} (Theorem \ref{extensionhomeo}) and applications to extending mapping classes (Corollary \ref{torelliextends}).
\subsection{Equivariant intersection form} Let $M$ be a connected, compact, oriented, topological\footnote{The results in this section are all in the topological category. However, we note that the 4-manifolds we consider are smoothable.} 4-manifold with fundamental group $\pi = \pi_1(M)$ and connected boundary. 

\begin{defn}
 The \textit{equivariant intersection form} \[\lambda_M : H_2(M;\mathbb{Z}[\pi]) \times H_2(M;\mathbb{Z}[\pi]) \to \Z[\pi]\] has right adjoint given by 
% https://q.uiver.app/#q=WzAsNCxbMCwwLCJIXzIoTTtcXG1hdGhiYntafVtcXHBpXzFNXSkiXSxbMSwwLCJIXzIoTSxcXHBhcnRpYWwgTTtcXG1hdGhiYntafVtcXHBpXzFNXSkiXSxbMiwwLCJIXjIoTTtcXG1hdGhiYntafVtcXHBpXzFNXSkgIl0sWzMsMCwiXFxtYXRocm17SG9tfV97XFxtYXRoYmJ7Wn1bXFxwaV8xTV19KEhfMihNO1xcbWF0aGJie1p9W1xccGlfMU1dKSxcXG1hdGhiYntafVtcXHBpXzFNXSkiXSxbMCwxXSxbMSwyLCJcXG1hdGhybXtQRH0iXSxbMiwzLCJcXG1hdGhybXtldn0iXV0=
\[\begin{tikzcd}
	{H_2(M;\mathbb{Z}[\pi])} & {H_2(M,\partial M;\mathbb{Z}[\pi])} & {H^2(M;\mathbb{Z}[\pi]) } & {\mathrm{Hom}_{\mathbb{Z}[\pi]}(H_2(M;\mathbb{Z}[\pi]),\mathbb{Z}[\pi]),}
	\arrow[from=1-1, to=1-2]
	\arrow["{\mathrm{PD}^{-1}}", from=1-2, to=1-3]
	\arrow["{\mathrm{ev}}", from=1-3, to=1-4]
\end{tikzcd}\]
the first map induced by inclusion of pairs, the second  Poincar\'{e}-Lefschetz duality with twisted coefficients, and the third evaluation (see e.g. \cite[Section 3.1]{equivariantintersectionform}). 
\end{defn}

%Two basic properties of $\lambda_M$ are as follows.
In the sequel, we make use of two basic algebraic properties that $\lambda_M$ satisfies. Recall that the group ring $\mathbb{Z}[\pi]$ carries a natural involution $x \mapsto \overline{x}$ determined by extending inversion on group elements linearly.

\begin{proposition}[{c.f. discussion after \cite[Proposition 3.12]{equivariantintersectionform}}]
The pairing $\lambda_M$ is \textit{Hermitian}, that is, it satisfies $\lambda_M(x,y) = \overline{\lambda_M(y,x)}$. It is also \textit{sesquilinear}, that is, $\lambda_M(px,qy) = p\lambda_M(x,y)\overline{q}$ for $p,q \in \mathbb{Z}[\pi]$.   
\end{proposition}

\begin{defn} 
 A homomorphism $\alpha : H_2(M;\mathbb{Z}[\pi]) \to H_2(M;\mathbb{Z}[\pi])$ of $\mathbb{Z}[\pi]$-modules  is an \textit{isometry} of $\lambda_M$ if it is an isomorphism of $\mathbb{Z}[\pi]$-modules and $\lambda_M(\alpha(x),\alpha(y)) = \lambda_M(x,y)$ for all $x,y \in H_2(M;\mathbb{Z}[\pi])$. Let $\mathrm{Isom}(\lambda_M)$ denote the isometry group of $\lambda_M$. 
\end{defn}

%A homomorphism of $\mathbb{Z}[\pi]$-modules $\alpha : H_2(M;\mathbb{Z}[\pi]) \to H_2(M;\mathbb{Z}[\pi]$ is an \textit{isometry} of $\lambda_M$ if it is an isomorphism of $\mathbb{Z}[\pi]$-modules and $\lambda_M(\alpha(x),\alpha(y)) = \lambda_M(x,y)$ for all $x,y \in H_2(M;\mathbb{Z}[\pi])$. Let $\mathrm{Iso}(\lambda_M)$ denote the isometry group of $\lambda_M$. 

\begin{rmk}
For simplicity, assume that $M$ is smooth. The pairing $\lambda_M\colon H_2(M;\mathbb{Z}[\pi]) \times H_2(M;\mathbb{Z}[\pi]) \to \mathbb{Z}[\pi]$ has a geometric description, justifying its title. It is given by
$$\lambda_M(x,y) = \sum_{g\in \pi} \langle y,gx\rangle g^{-1} \in \mathbb{Z}[\pi],$$
where $\langle-,-\rangle$ denotes the ordinary algebraic intersection of cycles in the universal cover of $M$ (see e.g. \cite[Section 3.1]{equivariantintersectionform}). 
\end{rmk}

%Let $\pi = \pi_1M$. Inclusion, Poincar\'{e}-Lefschetz duality with twisted coefficients, and the evaluation homomorphism defines a homomorphism
%\[ H_2(M;\mathbb{Z}[\pi]) \to H_2(M,\partial M;\mathbb{Z}[\pi]) \to H^2(M;\mathbb{Z}[\pi]) \to \mathrm{Hom}_{\mathbb{Z}[\pi]}(H_2(M;\mathbb{Z}[\pi]),\mathbb{Z}\pi).\]

%We now recall the relevant results of Conway-Powell \cite{Conway_Powell}.

\subsection{4-manifolds with $\pi_1 \cong \mathbb{Z}$}

 Assume that $\pi_1M \cong \mathbb{Z}$.  Let $\Lambda = \mathbb{Z}[\pi_1M] \cong \mathbb{Z}[t,t^{-1}]$. Assume further that $M$ has \textit{ribbon boundary}, that is, inclusion induces a surjection $\pi_1\partial M \to \pi_1M$, and that $H_1(\partial M;\Lambda)$ is a torsion $\Lambda$-module. Under these hypotheses, $H_2(M;\Lambda)$ and $H_2(M,\partial M;\Lambda)$ are free $\Lambda$-modules.

%Recall that $H_2(M;\Lambda) \cong H_2(\widetilde{M};\mathbb{Z})$ is equipped with an \textit{equivariant} intersection pairing \cite[Definition 3.2]{Conway_Piccirillo_Powell}

%\textcolor{red}{Define through duality and record remark about explicit formulae}. 

%\textcolor{red}{Conway--Piccirillo--Powell use $\lambda_M(x,y) = \sum_{k \in \Z} \langle y, t^k x\rangle t^{-k}$, perhaps we should stick to this convention since we use Conway--Powell so much? Agree yes, my bad!}

An isometry of $\lambda_M$ acts on $H_1(\partial M;\Lambda)$ as follows.  The evaluation map $\mathrm{ev}\colon H_2(M;\Lambda) \to \mathrm{Hom}_{\Lambda}(H_2(M;\Lambda),\Lambda)$ is injective (by an application of the universal coefficient spectral sequence, see e.g. \cite[Proposition 3.3]{equivariantintersectionform}) and $H_2(\partial M;\Lambda) \to H_2(M;\Lambda)$ is the zero map (proof of \cite[Lemma 3.2]{Conway_Powell}). Let $\overline{H_2(M;\Lambda)^*}$ denote the dual of $H_2(M;\Lambda)$ but with $\Lambda$-module structure $p\cdot x = \overline{p}x$ for $p\in \Lambda$. The right adjoint of $\lambda_M$ determines an injective $\Lambda$-linear map \[\widehat{\lambda}_M: H_2(M;\Lambda) \to \overline{H_2(M;\Lambda)} ^*\] defined by $\widehat{\lambda}_M(y) = \lambda_M(-,y)$. Inverting the dual of an element $\alpha \in \mathrm{Isom}(\lambda_M)$ induces an isomorphism \[[\alpha^{-*}]: \mathrm{coker}(\widehat{\lambda}_M) \to \mathrm{coker}(\widehat{\lambda}_M).\] We have the following commutative diagram with exact rows \cite[Remark 3.3]{Conway_Powell}.

\begin{comment}
\[\begin{tikzcd} \label{diagramdefiningdm}
	0 & {H_2(M;\Lambda)} & {H_2(M;\Lambda)^*} & {\mathrm{coker}(\widehat{\lambda}_M)} & 0 \\
	0 & {H_2(M;\Lambda)} & {H_2(M,\partial M;\Lambda)} & {H_1(\partial M;\Lambda)} & 0.
	\arrow[from=1-1, to=1-2]
	\arrow["{\widehat{\lambda}_M}", from=1-2, to=1-3]
	\arrow[from=1-3, to=1-4]
	\arrow[from=1-4, to=1-5]
	\arrow["{D_M}"', from=1-4, to=2-4]
	\arrow[from=2-1, to=2-2]
	\arrow["{=}", no head, from=2-2, to=1-2]
	\arrow[from=2-2, to=2-3]
	\arrow["{\mathrm{ev}\circ\mathrm{PD}}", from=2-3, to=1-3]
	\arrow["\partial", from=2-3, to=2-4]
	\arrow[from=2-4, to=2-5]
\end{tikzcd}\]
\end{comment}

% https://q.uiver.app/#q=WzAsOCxbMSwxLCJIXzIoTTtcXExhbWJkYSkiXSxbMCwxLCIwIl0sWzIsMCwiXFxvdmVybGluZXtIXzIoTTtcXExhbWJkYSleKn0iXSxbMiwyLCJIXzIoTSxcXHBhcnRpYWwgTTtcXExhbWJkYSkiXSxbMywwLCJcXG1hdGhybXtjb2tlcn0oXFx3aWRlaGF0XFxsYW1iZGFfTSkiXSxbNCwwLCIwIl0sWzMsMiwiSF8xKFxccGFydGlhbCBNO1xcTGFtYmRhKSJdLFs0LDIsIjAiXSxbMSwwXSxbMCwyLCJcXHdpZGVoYXRcXGxhbWJkYV9NIl0sWzAsMywiaiIsMl0sWzIsNF0sWzQsNV0sWzMsNl0sWzYsN10sWzMsMiwiXFxtYXRocm17ZXZ9XFxjaXJjIFxcbWF0aHJte1BEfSIsMl0sWzQsNiwiRF9NIl1d

\begin{equation}
\label{commutativediagram}
\begin{tikzcd} 
	&& {\overline{H_2(M;\Lambda)^*}} & {\mathrm{coker}(\widehat\lambda_M)} & 0 \\
	0 & {H_2(M;\Lambda)} \\
	&& {H_2(M,\partial M;\Lambda)} & {H_1(\partial M;\Lambda)} & 0
	\arrow[from=1-3, to=1-4]
	\arrow[from=1-4, to=1-5]
	\arrow["{D_M}", from=1-4, to=3-4]
	\arrow[from=2-1, to=2-2]
	\arrow["{\widehat\lambda_M}", from=2-2, to=1-3]
	\arrow["j"', from=2-2, to=3-3]
	\arrow["{\mathrm{ev}\circ \mathrm{PD}^{-1}}"', from=3-3, to=1-3]
	\arrow[from=3-3, to=3-4]
	\arrow[from=3-4, to=3-5]
\end{tikzcd}
\end{equation}

The lower sequence arises from the long exact sequence of the pair $(M,\partial M)$. The middle vertical arrow is evaluation after Poincar\'{e}-Lefschetz duality with twisted coefficients. The diagram therefore induces an isomorphism of $\Lambda$-modules $D_M : \mathrm{coker}(\widehat{\lambda}_M) \to H_1(\partial M;\Lambda)$. 

\begin{defn}[{\cite[Definition 2.6]{Conway_Powell}}] \label{boundaryofiso}
The \textit{boundary} of an isometry $\alpha \in \mathrm{Isom}(\lambda_M)$ is \[\Pi_M(\alpha) = D_M \circ [\alpha^{-*}] \circ D_M^{-1}\colon H_1(\partial M;\Lambda) \to H_1(\partial M;\Lambda)\]
This defines a homomorphism 
\begin{align*} \Pi_M :  \mathrm{Isom}(&\lambda_M) \to \mathrm{Aut}_\Lambda(H_1(\partial M;\Lambda)) \\
&\alpha \to D_M\circ[\alpha^{-*}] \circ D_M^{-1}.
\end{align*}
Denote the image of $\Pi_M$ by $\Gamma_M \leq \mathrm{Aut}(H_1(\partial M;\Lambda))$. 
%where $D_M : \mathrm{coker}(\widehat{\lambda}_M) \to H_1(\partial M;\Lambda)$ is defined by $D_M([x]) = \partial (\mathrm{PD}^{-1}(x))$ where $\partial : H_2(M,\partial M) \to H_1(\partial M)$ is the connecting homomorphism and $PD$ is Poincar\'{e}-Lefschetz duality with twisted coefficients. 
\end{defn}

\begin{theorem} \label{extensionhomeo}
Let $M$ be a compact oriented 4-manifold with $\pi_1M \cong \mathbb{Z}$ whose boundary is ribbon and has $H_1(\partial M;\Lambda)$ a torsion $\Lambda$-module. Fix an orientation preserving homeomorphism $h : \partial M \to \partial M$ which intertwines the inclusion map $\pi_1\partial M \to \pi_1M \cong \mathbb{Z}$. Assume that a lift $h_*$ of $h$ to $H_1(\partial M;\Lambda)$ is the boundary (Definition \ref{boundaryofiso}) of an isometry in $\mathrm{Isom}(\lambda_M)$, that is, $h_* \in \Gamma_M$. Then there exists an orientation preserving homeomorphism $H : M \to M $ with $H\rvert_{\partial M} = h$. 
\end{theorem}

\begin{proof}
This is a special case of Theorem 4.1 in \cite{Conway_Powell}. 
\end{proof}

\subsection{The image of $\Pi_M$ and the Blanchfield pairing} The group $\Gamma_M$ is governed by additional structure on $H_1(\partial M;\Lambda)$, namely the \textit{Blanchfield pairing} on $\partial M$. Let $\mathcal{Q}$ denote the field of fractions of $\Lambda$. As the boundary of $M$ is ribbon, inclusion induces a surjection $\pi_1(\partial M) \to \mathbb{Z}$. Moreover, we've assumed that $H_1(\partial M;\Lambda)$ is a torsion $\Lambda$-module. It follows that $H_1(\partial M;\Lambda)$ is equipped with a sesquilinear, Hermitian, non-singular \textit{linking form}
\[ \mathrm{Bl}_{\partial M} : H_1(\partial M;\Lambda) \times H_1(\partial M;\Lambda) \to \mathcal{Q}/\Lambda \]
referred to as the Blanchfield pairing on $\partial M$ (see \cite[Definition 3.4]{Conway_Powell} for a precise definition). Let $\mathrm{Isom}(\mathrm{Bl}_{\partial M})$ denote the isometry group of $\mathrm{Bl}_{\partial M}$. 

\begin{proposition} \label{boundaryofisopreservesblanch}
The boundary of an isometry of $\lambda_M$ is an isometry of the Blanchfield pairing $\mathrm{Bl}_{\partial M}$, that is, $\Gamma_M \leq \mathrm{Isom}(\mathrm{Bl}_{\partial M})$. 
\end{proposition}

\begin{proof}
This is \cite[Proposition 3.5]{Conway_Powell} in combination with \cite[Lemma 2.5]{Conway_Powell}. 
\end{proof}

Recall that a surface $\Sigma$ in a 4-manifold $X$  is a $\mathbb{Z}$\textit{-surface} if $\pi_1(X-\Sigma) \cong \mathbb{Z}$. Conway-Powell identify the group $\mathrm{Isom}(\mathrm{Bl}_{\partial M})$ explicitly when $M$ is a $\mathbb{Z}$-surface exterior (this identification will be used in Section 5). To state the results, we first describe the boundary of a $\mathbb{Z}$-surface exterior (in the rel boundary setting). %If $\Sigma$ is a closed $\mathbb{Z}$-surface in a closed 4-manifold $X$, then $\Sigma$ has trivial normal bundle $\nu\Sigma \subset X$ \cite[Lemma 5.2]{Conway_Powell}. Let $X_\Sigma = X-\nu\Sigma$. 
If $\Sigma \subset D^4$ is a properly embedded $\mathbb{Z}$-surface with boundary a knot $K \subset S^3$, let $D^4_\Sigma = D^4 -\nu \Sigma$. Then $\partial D^4_\Sigma \cong \Sigma \times S^1 \cup_\partial E_K$ where $E_K$ denotes the knot exterior of $K$ in $S^3$. Mayer-Vietoris with $\Lambda$-coefficients shows that 
\begin{equation} H_1(\partial D^4_\Sigma;\Lambda) \cong H_1(\Sigma \times S^1;\Lambda) \oplus H_1(E_K;\Lambda). \label{SplittingofH_1MV}
\end{equation}
Moreover, $H_1(\Sigma \times S^1;\Lambda) \cong (\Lambda/(t-1))^{\oplus 2g}$ where $g$ is the genus of $\Sigma$. 
\begin{proposition} \label{isoofblanchfield}
%For a closed $\mathbb{Z}$-surface $\Sigma \subset X$, a trivialization of $\partial X_\Sigma \cong \Sigma \times S^1$ determines an isomorphism 
%\[ \mathrm{Iso}(\mathrm{Bl}_{\partial X_\Sigma}) \cong \mathrm{Sp}(H_1(\Sigma;\mathbb{Z})).\] 
For a properly embedded $\mathbb{Z}$-surface $\Sigma \subset D^4$ with boundary a knot $K \subset S^3$, the isomorphism \eqref{SplittingofH_1MV} induces a decomposition
\[ \mathrm{Isom}(\mathrm{Bl}_{\partial D^4_\Sigma}) \cong \mathrm{Isom}(\mathrm{Bl}_{\Sigma \times S^1})  \oplus \mathrm{Isom}(\mathrm{Bl}_{E_K})\]
and $\mathrm{Isom}(\mathrm{Bl}_{\Sigma \times S^1}) \cong \mathrm{Sp}(H_1(\Sigma;\mathbb{Z}))$. 
\end{proposition}
\begin{proof}
This is \cite[Propositions 5.6 and 5.7]{Conway_Powell}. 
\end{proof}

\subsection{Extending mapping classes}

We now record an application of Theorem \ref{extensionhomeo} in the language of extendable mapping classes. For a compact surface $\Sigma$, let  \[\Psi : \Mod(\Sigma) \to \textrm{Sp}(H_1(\Sigma;\mathbb{Z}))\] be the \textit{symplectic representation} of $\Mod(\Sigma)$. The kernel of $\Psi$ is the \textit{Torelli subgroup} of $\Mod(\Sigma)$, denoted $\mathcal{I}(\Sigma)$.

%\begin{corollary}
%Let $\Sigma$ be a closed $\mathbb{Z}$-surface in a closed 4-manifold $X$. Let $X_\Sigma = X-\nu \Sigma$. A trivialization of the $S^1$-bundle $\partial X_\Sigma \cong \Sigma \times S^1$ determines an isomorphism \[ \mathrm{Iso}(\mathrm{Bl}_{\partial X_\Sigma}) \cong \mathrm{Sp}(H_1(\Sigma;\mathbb{Z})).\] 
%View the image $\Gamma_{X_\Sigma}$ of \(\Pi_{X_\Sigma}\colon\mathrm{Iso}(\lambda_{X_\Sigma}) \to \mathrm{Iso}(\mathrm{Bl}_{\partial X_\Sigma})\) as a subgroup of the symplectic group $\mathrm{Sp}(H_1(\Sigma;\mathbb{Z}))$. Then 
%\[\mathcal{E}^{\textrm{TOP}}(X,\Sigma) = \Psi^{-1}(\Gamma_{X_\Sigma}).\]
%It follows that $\mathcal{I}(\Sigma) \leq \mathcal{E}^{\textrm{TOP}}(X,\Sigma)$. 
%\end{corollary}

%\begin{proof}
%This follows from Theorem \ref{extensionhomeo} and Proposition \ref{isoofblanchfield}. 
%\end{proof}

%We now record an application in the presence of boundary. Corollary \ref{torelliextends} is essential in proving Theorem \ref{maintheorem}. 

\begin{corollary} \label{torelliextends}
Let $\Sigma \subset D^4$ be a properly embedded $\mathbb{Z}$-surface with a single boundary component. Then \[\mathcal{I}(\Sigma) \leq \mathcal{E}^{\textrm{TOP}}(D^4,\Sigma).\]
\end{corollary}

\begin{proof}
This follows from Theorem \ref{extensionhomeo} and the methods in the proof of \cite[Theorem 5.8]{Conway_Powell}. Fix a homeomorphism $f : \Sigma \to \Sigma$ which is the identity on $\partial \Sigma$ and acts trivially on $H_1(\Sigma;\mathbb{Z})$. Write $\partial \Sigma = K \subset \partial D^4$ and $E_K$ for the exterior of the knot $K$. Let $D^4_\Sigma = D^4-\nu\Sigma$ be the surface exterior. The boundary of $D^4_\Sigma$ is obtained from the knot exterior $E_K$ by gluing in $\Sigma \times S^1$ along $\partial E_K$. Define a homeomorphism $h : \partial D^4_\Sigma \to \partial D^4_\Sigma$ by 
\[ h = (f \times \mathrm{Id}_{S^1}) \cup_{\partial} \mathrm{Id}_{E_K}.\] The homeomorphism $h$ intertwines $\pi_1 \partial D^4_\Sigma \to \mathbb{Z}$ as it is the identity on $\partial E_K$. In particular, we may consider lifts to automorphisms of $H_1(\partial D^4_\Sigma;\Lambda)$. Recall that Mayer-Vietoris with $\Lambda$-coefficients yields an isomorphism \cite[Lemma 5.5]{Conway_Powell}
\[ H_1(\partial D^4_\Sigma ;\Lambda) \cong  H_1(\Sigma\times S^1;\Lambda) \oplus H_1(E_K;\Lambda).\] In particular, $h_* : H_1(\partial D^4_\Sigma;\Lambda) \to H_1(\partial D^4_\Sigma;\Lambda)$ is the identity. By Theorem \ref{extensionhomeo}, there exists a homeomorphism $H : D^4_\Sigma \to D^4_\Sigma$ extending $h$. Extend $H$ over $D^4$ using $(f \times {\mathrm{Id}_{D^2}})$ on $\nu \Sigma \cong \Sigma \times D^2$. By construction, the resulting homeomorphism of $D^4$ is the identity on $\partial D^4$. It follows that $[f] \in \mathcal{E}^{\mathrm{TOP}}(D^4,\Sigma)$. 
\end{proof}

\begin{rmk}
By the Alexander trick, the homeomorphism $H : (D^4,\Sigma) \to (D^4,\Sigma)$ constructed in the proof of Corollary \ref{torelliextends} is topologically isotopic to the identity relative to the boundary $\partial D^4$. 
\end{rmk}

When the equivariant intersection form can be computed explicitly, isometries may be  constructed directly. We carry this out for pushed in Seifert surfaces in the 4-ball.

\subsection{Pushed in Seifert surfaces in $D^4$}
 Let $F' \subset S^3$ be a Seifert surface for a knot $K$. Let $\iota \colon F'\times [0,1] \to D^4$ be an isotopy which pushes the interior of $F'$ into the interior of $D^4$. Our pushed in Seifert surface is $F = \iota(F'\times \{1\})$. Let $V_F = D^4-\nu F$ denote the exterior of $F$. 
 
 \begin{proposition} \label{pushedinareZsurfaces}
The pushed in Seifert surface $F \subset D^4$ is a $\mathbb{Z}$-surface, that is, $\pi_1(V_F) \cong \mathbb{Z}$. 
 \end{proposition}

 \begin{proof}
By \cite[Proposition 6.2.1]{gompfstipsicz}, the exterior $V_F$ admits a handle decomposition with only a single 1-handle (recall that $F$ is, by definition, connected). Therefore $\pi_1(D^4-\nu F)$ is a quotient of $\mathbb{Z}$. The normal bundle $\nu F$ to $F$ in $D^4$ is the trivial disk bundle $F\times D^2$ over $F$ (it is an orientable disk bundle over $F$, and $F$ is homotopy equivalent to a wedge of circles). The long exact sequence of the pair $(D^4-F,\partial \nu F) \simeq (D^4-F, \nu F-F)$ and excision of $F$ applied to the pair $(D^4,\nu F)$ show that $H_1(V_F;\mathbb{Z}) \cong \mathbb{Z}$. 
 \end{proof}
 Let $\widetilde{V_F}$ denote the universal cover of $V_F$. We begin by describing the standard model for $\widetilde{V_F}$. For details, we refer the reader to \cite[Page 537]{Ko}, \cite[proof of Lemma 5.4]{Cochran_Orr_Teichner} as well as \cite{AkbulutKirby}. Note that this procedure works for any pushed in $F^{n-2}$ in an $n$-ball $D^n$.  %the equivariant intersection form $\lambda_{V_F}$ on $H := H_2(V_F;\Lambda)$. 

Cut $D^4$ along the trace $S = \iota(F'\times [0,1])$ of the isotopy $\iota$ which pushed the interior of $F'$ into the interior of $D^4$. The result $D^4-\nu S$ is homeomorphic to $D^4$ (cf. \cite[Proposition 6.2.1]{gompfstipsicz}). This identification provides a thickening $\bar{F'} \subset \partial (D^4-\nu S) \cong \partial D^4$ of the interior of the Seifert surface $F'$ in the boundary of $D^4$. Precisely,  there is an embedding of $\bar{F'} = F' \times [-1,1]$ in $\partial (D^4-\nu S) \cong \partial D^4$. Take countably many copies $t^k D^4,\; k \in \Z$ of $D^4$. Fix $\epsilon >0$ small. Set $$\bar{F'}_{+}(k) = \{ (x,t) \in \bar{F'} \; | \; x \in F', t \in [\epsilon,1]\} \subset t^k D^4,$$
and  $$\bar{F'}_-(k) = \{ (x,t) \in \bar {F'} \; | \; x \in F', t \in [-1,-\epsilon]\} \subset t^k D^4.$$
Let 
\[\widetilde{V_F} = \bigsqcup_{k \in\mathbb{Z}} t^k D^4/ \sim \]
where we identify $\bar{F'}_+(k)$ with $\bar {F'}_-(k-1)$ via the map $h: \bar{F'}_+(k) \rightarrow \bar{F'}_-(k-1)$ that sends $(x,t)$ to $(x,-t)$.

%The boundary of $\widetilde{V_F}$ can computed as follows. Let $\bar F(k) = \{ (x,t) \in \bar F \; | \; t \in [-\epsilon, \epsilon]\} \subset t^kD^4$. Then  $$\partial \widetilde{V_F} = (\bigcup_{k \in \Z} \bar F(k)) \cup \widetilde{E_K},$$
%where $\cup_{k \in \Z} \bar F(k) \cong F \times \R$ and $\widetilde{E_K}$ denotes the universal abelian cover of $E_K = S^3  - \nu(K)$ (constructed the usual manner [Rolfsen, Chapter 5c, Remark 5 and Exercise 6]). 

The connecting homomorphism \[\delta: H_2(\widetilde{V_F};\Z) \rightarrow H_1(F';\Z) \otimes_{\Z} \Lambda\] for the Mayer-Vietoris sequence of the decomposition $A = \cup_{k \; \mathrm{even}} t^kD^4,\; B = \cup_{k \; \mathrm{odd}}t^kD^4,\; A \cap B = \cup_{k \in \Z} \bar{F'}_+(k)$ gives an isomorphism of $\Lambda$-modules \cite[Page 538]{Ko}. Therefore $H_2(\widetilde{V_F};\Z)$ is a free $\Lambda$-module. %Henceforth, we no longer distinguish between the Seifert surface $F' \subset S^3$ and its push in $F \subset D^4$.  

\subsection{Basing $H_2(V_F;\Lambda)$ and computing $\lambda_{V_F}$} To obtain a $\Lambda$-basis for $H_2(V_F;\Lambda)$, fix a collection $\{ a_i\}_{i=1}^{2g}$ of oriented simple closed curves on $F'$ that give a $\Z$-basis for $H_1(F';\Z)$. Take immersed $2$-disks $(\Sigma_i^+, \partial \Sigma_i^+) \subset (tD^4, \bar{F'}_+(1))$ and $(\Sigma_i^-,\partial \Sigma _i^-) \subset (t^0 D^4, \bar{F'}_-(0))$ such that $\partial \Sigma_i^+ = (a_i,\epsilon) \in \bar{F'}_+(1)$ and $\partial \Sigma_i^- = (-a_i,-\epsilon)$; we identify $(a_i,\epsilon)$ with $a_i^+$, a pushoff of $a_i$ in the positive normal direction to $F'$ in $S^3$. Similarly, we identify $(-a_i,-\epsilon)$ with $-a_i^-$, which denotes $a_i$ with opposite orientation pushed off in the negative normal direction to $F'$. Set $\Sigma_i = \Sigma_i^+ \cup \Sigma_i^-$. Then the spheres $\{ \Sigma_i\}_{i=1}^{2g}$ represent a $\Lambda$-basis for $H_2(\widetilde{V_F};\Z)$. Let $S$ denote the Seifert form of $F'$ relative to the basis $\{a_i\}_{i=1}^{2g}$ for $H_1(F';\mathbb{Z})$. In the constructed basis, the equivariant intersection form $\lambda_{V_F}$ is a \textit{deformation} of $S$. %(The computation of the equivariant intersection form is an application of [Rolfsen, Chapter 5.D, Exercise 9] to $\Sigma_i^+ \subset tD^4$ and $\Sigma_i^- \subset t^0D^4$.)

\begin{proposition} \label{basisandform}
The spheres $\Sigma_1,\ldots,\Sigma_{2g}$ constructed by doubling along cycles above give a basis $\mathcal{C}$ for the free $\Lambda$-module $H_2(\widetilde{V_F};\Z)$. In the basis $\mathcal{C}$ the equivariant intersection form 
\begin{align*}
        \lambda_{V_F}: H_2(\widetilde{V_F};\Z) \times H_2(\widetilde{V_F};\Z) \rightarrow \Lambda \\
        (\Sigma_j, \Sigma_i) \mapsto \sum_k (\Sigma_i \cdot t^k \Sigma_j)t^{-k}
\end{align*}
is given by the matrix $Q \in M_{2g}(\Lambda)$ where
$$Q_{ij} = \lambda_{V_F}(\Sigma_i, \Sigma_j) = (1-t)S_{ji} + (1-t^{-1})S_{ij}.$$ 
\end{proposition}

\begin{proof}
A reference is \cite{Ko}. See also e.g. the proof of \cite[Lemma 5.4]{Cochran_Orr_Teichner} (note our conventions for $\lambda_{V_F}$). %The computation of the equivariant intersection form is an application of \cite[Chapter 5.D, Exercise 9]{Rolfsen} to $\Sigma_i^+ \subset tD^4$ and $\Sigma_i^- \subset t^0D^4$. 
\end{proof}

Henceforth, we identify $F \subset D^4$ and $F' \subset S^3$. Let \(\delta \colon H_2(V_F,\partial V_F;\Lambda) \to H_1(\partial V_F;\Lambda) \) denote the connecting homomorphism in the long exact sequence of the pair $(V_F,\partial V_F)$ with $\Lambda$ coefficients. Recall that $H_1(\partial V_F;\Lambda) \cong H_1(F\times S^1;\Lambda) \oplus H_1(E_K;\Lambda)$ \eqref{SplittingofH_1MV} and $H_1(F \times S^1;\Lambda) = H_1(F\times\mathbb{R};\mathbb{Z}) \cong H_1(F;\mathbb{Z})$. 

\begin{lemma} \label{connectinghomomorphism}
%Let \[\delta \colon H_2(V_F,\partial V_F;\Lambda) \to H_1(\partial V_F;\Lambda) \cong H_1(F\times S^1;\Lambda) \oplus H_1(E_K;\Lambda) \] denote the connecting homomorphism in the long exact sequence of the pair $(V_F,\partial V_F)$ with $\Lambda$ coefficients. 
Fix $i \in \{1,\ldots,2g\}$. Then  $\delta(\Sigma_i^+) = [a_i] \in H_1(F;\mathbb{Z}) \subset H_1(\partial V_F;\Lambda). $ 
%The boundary $\partial \widetilde{V_F} \cong F'\times \mathbb{R} \cup \widetilde{E_K}$ where $\widetilde{E_K}$ denotes the universal abelian covering of the knot exterior $E_K$. The Mayer-Vietoris sequence associated to this union yields an isomorphism 
%\[ H_1(\partial \widetilde{V_F};\mathbb{Z}) \cong H_1(F'\times \mathbb{R};\mathbb{Z}) \oplus H_1(\widetilde{E_K};\mathbb{Z}).\] In particular, for the connecting homomorphism \[\delta:H_2(V_F, \partial V_F;\Lambda) \rightarrow H_1(\partial V_F;\Lambda)\] in the long exact sequence of the pair $(V_F, \partial V_F)$, we have $\delta(\Sigma_i^+) = (a_i,\epsilon) \in H_1(F' \times \R;\Z) \subset H_1(\partial \widetilde{V_F};\Z)$.
\end{lemma}

\begin{proof} 
The connecting homomorphism is \[\delta : H_2(\widetilde{V_F},\partial \widetilde{V_F};\mathbb{Z}) \to H_1(\partial \widetilde{V_F};\mathbb{Z}).\] Set $\bar{F'}(k) = \{ (x,t) \in \bar{F'} \; | \; t \in [-\epsilon, \epsilon]\} \subset t^kD^4$. Then  $$\partial \widetilde{V_F} = \bigcup_{k \in \Z} \bar{F'}(k)\cup \widetilde{E_K},$$
where $\bigcup_{k \in \Z} \bar{F'}(k) \cong F \times \R$ and $\widetilde{E_K}$ denotes the universal abelian cover of $E_K$. Here we implicitly use the standard description of $\widetilde{E_K}$ (c.f. \cite{Rolfsen}). The decomposition \eqref{SplittingofH_1MV} is $H_1(\partial \widetilde{V_F} ;\mathbb{Z}) \cong H_1(F\times \mathbb{R};\mathbb{Z}) \oplus H_1(\widetilde{E_K};\mathbb{Z})$. Now $\delta(\Sigma_i^+) = [\partial \Sigma_i^+] = [(a_i,\epsilon)] \in H_1(F\times \mathbb{R};\mathbb{Z})$.   
\end{proof}

\subsection{Computing the boundary homomorphism $\Pi_{V_F}$} 
%\textcolor{red}{Set up dual basis. Recall boundary of isometry as acting on the cokernel. Recall decomposition of this cokernel. We determine the explicit action on the $H_1(F\times S^1;\Lambda)$ factor.}. 

We now compute the restriction of the boundary of an isometry of $\lambda_{V_F}$ to $H_1(F\times S^1;\Lambda)$. Precisely, for $\alpha \in \mathrm{Isom}(\lambda_{V_F})$ and 
\[ \Pi_{V_F} \colon \mathrm{Isom}(\lambda_{V_F}) \to \mathrm{Isom}(\mathrm{Bl}_{\partial V_F}) \]%\cong \mathrm{Sp}(H_1(F;\mathbb{Z})) \oplus \mathrm{Iso}(\mathrm{Bl}_{E_K})\]
we compute $\Pi_{V_F}(\alpha)\rvert_{H_1(F\times S^1;\Lambda)}$ (see Proposition \ref{isoofblanchfield}).

Let $H^* = \overline{H_2(V_F;\Lambda)^*}$ (the dual of $H_2(V_F;\Lambda)$ but with $\Lambda$-module structure $p\cdot x = \overline{p}x$). Write $\lambda = \lambda_{V_F}$. 
%along the \textit{surface factor} $F\times S^1$ of $\partial V_F $. Recall that inclusion induces a decomposition \[H_1(\partial V_F;\Lambda) \cong H_1(F\times S^1;\Lambda) \oplus H_1(E_K;\Lambda),\] and that this induces a splitting
%$\[ \mathrm{Iso}(\mathrm{Bl}_{\partial V_F}) \cong \mathrm{Sp}(H_1(F;\mathbb{Z})) \oplus \mathrm{Iso}(\mathrm{Bl}_{E_K}).\]
Recall that we have an identification $D_{V_F} \colon \mathrm{coker}(\widehat{\lambda}_{V_F}) \to H_1(\partial V_F;\Lambda)$ (Definition \ref{boundaryofiso}), and that the map
$$H_2(V_F,\partial V_F;\Lambda) \xrightarrow{\mathrm{PD}^{-1}} H^2(V_F;\Lambda) \xrightarrow{\mathrm{ev}} H^*$$
is identified with the currying map of the equivariant intersection form 
$$\lambda^{\delta}: H_2(V_F;\Lambda) \times H_2(V_F,\partial V_F;\Lambda) \rightarrow \Lambda.$$
A basis for the $H_1(F \times S^1 ; \Lambda)$ factor in $\mathrm{coker}(\widehat \lambda_{V_F})$ is given by $[ \lambda^{\delta}(-,\Sigma_i^+)] \in \mathrm{coker}(\widehat{\lambda}_{V_F})$
for $i=1,\ldots,2g$. The $\Lambda$-basis $\{\Sigma_i\}_{i=1}^{2g}$ for $H_2(V_F;\Lambda)$ gives the dual basis $\{ \Sigma_i^*\}_{i=1}^{2g}$ for $H^*$. Recall that $S$ denotes the Seifert form of $F$ in the basis $\{a_i\}_{i=1}^{2g}$ for $H_1(F;\mathbb{Z}) \cong H_1(F\times S^1;\Lambda)$,  and that $j : H_2(V_F;\Lambda) \to H_2(V_F,\partial V_F;\Lambda)$ is the homomorphism induced by inclusion of pairs \eqref{commutativediagram}.

\begin{lemma} \label{coker_basis_for_surface_factor}
    We have $$\lambda^{\delta}(-,\Sigma_i^+) = \sum_{j=1}^{2g}(S_{ij} - S_{ji}t^{-1}) \Sigma_j^*,$$ and $$j(\Sigma_i) = (1-t^{-1})\Sigma_i^+.$$ 
    %Fix the $\Lambda$-basis $\mathcal{C} = \{ \Sigma_i\}_{i=1}^{2g}$ for $H_2(V_F;\Lambda)$ considered above. 
    Therefore $(1-t)\lambda^{\delta}(-,\Sigma_i^+)= \lambda(-,\Sigma_i) = 0 \in \mathrm{coker}(\widehat\lambda_{V_F})$. Let $\alpha \in \mathrm{Isom}(\lambda_{V_F})$ and let $A$ denote the matrix of $\alpha$ with respect to the basis $\mathcal{C}$ (Proposition \ref{basisandform}). Consider the $\Z$-submodule $U \subset \mathrm{coker}(\widehat \lambda_{V_F})$ spanned by the $\Z$-basis $\{ \lambda^{\delta}(-,\Sigma_i^+)\}_{i=1}^{2g} = \mathcal{B}$. Then the matrix representing the action of $\alpha^*$ on $U$ is given by $ A^{-1} \pmod{t-1}$ in the basis $\mathcal{B}$.
\end{lemma}
\begin{proof}
    We have $\lambda^{\delta}(-,\Sigma_i^+) = \sum_{j=1}^{2g}a_j \Sigma_j^*$ where $a_j = \lambda^{\delta}(\Sigma_j,\Sigma_i^+) = \sum_{k \in \Z} (\Sigma_i^+ \cdot t^k \Sigma_j)t^{-k}$. We have $\Sigma_i^+ \cdot t^k \Sigma_j = 0 $ for $k \neq 0,1$. For $k=0$, we have $\Sigma_i^+ \cdot \Sigma_j = \Sigma_i^+ \cdot \Sigma_j^+ = lk(a_i^+, a_j^+) = lk(a_i, a_j^+) = S_{ij}$. For $k=1$ we have $\Sigma_i^+ \cdot t\Sigma_j = \Sigma_i^+ \cdot t\Sigma_j^- = lk(a_i^+, -a_j^-) = - S_{ji}$. This shows the first claim.

    For the second claim, we have a commutative diagram 
    \[\begin{tikzcd}
	{H_2(V_F;\Lambda)} & {H^*} \\
	{H_2(V_F,\partial V_F;\Lambda)}
	\arrow["{\hat \lambda}", from=1-1, to=1-2]
	\arrow["j"', from=1-1, to=2-1]
	\arrow["{\hat \lambda^{\delta}}"', from=2-1, to=1-2]
    \end{tikzcd}\]
where $j$ comes from the exact sequence 
    $$0 \rightarrow H_2(V_F;\Lambda) \xrightarrow{j} H_2(V_F,\partial V_F;\Lambda) \xrightarrow{\delta} H_1(\partial V_F;\Lambda) \rightarrow 0.$$
The basis $\{\Sigma_i\}$ for $H_2(V_F;\Lambda)$ has the form $\Sigma_i = \Sigma_i^+ \cup \Sigma_i^-$ and $j(\Sigma_i) = \Sigma_i^+ + \Sigma_i^-$. We have (using Proposition \ref{basisandform} which computes $\lambda(\Sigma_j,\Sigma_i$)), 
    \begin{align*}
       (1-t)S_{ij} + (1-t^{-1})S_{ji} &=  \lambda(\Sigma_j, \Sigma_i) \\
       &= \lambda^{\delta}(\Sigma_j, j(\Sigma_i)) \\
       &= \lambda^{\delta}(\Sigma_j, \Sigma_i^+ + \Sigma_i^-) \\
        &= S_{ij} - S_{ji}t^{-1} + \lambda^{\delta}(\Sigma_j, \Sigma_i^-). %\\
        %= (1-t)S_{ij} + (1-t^{-1})S_{ji}.
    \end{align*}
Hence we have $$\lambda^{\delta}(\Sigma_j, \Sigma_i^-) = S_{ji} - tS_{ij},$$
which implies that 
    $$\lambda^{\delta}(-, \Sigma_i^-) = \sum_{j=1}^{2g}\lambda^{\delta}(\Sigma_j, \Sigma_i^-) \Sigma_j^* = \sum_{j=1}^{2g}(S_{ji}-tS_{ij})\Sigma_j^* = -t\lambda^{\delta}(-,\Sigma_i^+).$$
Since $\hat \lambda^{\delta}$ is an isomorphism, we get that $j(\Sigma_i) = \Sigma_i^+ + \Sigma_i^- = (1-t^{-1})\Sigma_i^+$.

For the last claim, write $\alpha^{-1}(\Sigma_i) = \sum_j b_{ji} \Sigma_j$, then $\alpha^*\lambda^{\delta}(-, \Sigma_i^+) = \lambda^{\delta}(\alpha(-), \Sigma_i^+)$. We have 
\begin{align*}
    (1-t)\lambda^{\delta}(\alpha(-), \Sigma_i^+) = \lambda^{\delta}(\alpha(-), j(\Sigma_i)) \\
    = \lambda(\alpha(-), \Sigma_i) \\
    = \lambda(-, \alpha^{-1}(\Sigma_i)) \\
    = \sum_j \bar b_{ji} \lambda(-, \Sigma_j) \\
    = (1-t) \sum_j \bar b_{ji} \lambda^{\delta}(-, \Sigma_j^+).
\end{align*}
Since $\Lambda$ is an integral domain, we obtain 
$$\alpha^* \lambda^{\delta}(-,\Sigma_i^+) = \sum_j \bar b_{ji} \lambda^{\delta}(-,\Sigma_j^+).$$
Note that $H^* = \mathrm{Hom}_{\Lambda}(H_2(V_F;\Lambda),\Lambda)$ is equipped with the $\Lambda$-module structure $p \cdot h = \bar p h$ so that the maps $\hat \lambda,\; \hat \lambda^{\delta}$ become $\Lambda$-module maps, so the final claim follows.
\end{proof}

%We have that $U \subset \mathrm{coker}(\hat \lambda)$ is the $H_1(F \times S^1;\Lambda)$ factor, and we have that $(1-t) \lambda^{\delta}(-,\Sigma_i^+) = 0 \in \mathrm{coker}(\hat \lambda)$. 

Relative to our basis $\mathcal{C}$ for $H_2(V_F;\Lambda)$ we have $$\mathrm{Isom}(H_2(V_F;\Lambda), \lambda) = \{ A \in GL_{2g}(\Lambda) \; | \; A^TQ \bar A = Q\}.$$ We have $Q= (1-t)P$, where $P = S^T - t^{-1}S$. Any matrix $A \in \mathrm{GL}_{2g}(\Lambda)$ with $A^TQ \bar A = Q$ must satisfy \begin{equation}
    A^T P \bar A = P . \label{PmatrixAeqn}
\end{equation} Note that $P \pmod{t-1} = S^T-S$ is the  intersection form on $H_1(F;\Z)$ in our basis of curves $\{a_i\}_{i=1}^{2g}$ \cite[Chapter 8, Exercise 7]{Rolfsen}. Therefore $A \pmod{t-1}$ has image in $\mathrm{Sp}_{2g}\Z$. 

\begin{corollary} \label{t=1}
The map 
\begin{align*}
    \Pi_{V_F}: \mathrm{Isom}(\lambda_{V_F}) \rightarrow \mathrm{Isom}(\mathrm{Bl}_{\partial V_F}) \\
    \alpha \mapsto [\alpha^{-*}]
\end{align*}
of Definition \ref{boundaryofiso} is given by $A \mapsto A \pmod{t-1}$ on the $H_1(F \times S^1;\Lambda)$ factor of $H_1(\partial V_F;\Lambda)$ where $A$ is the matrix of $\alpha$ with respect to the basis $\mathcal{C}$ of $H_2(V_F;\Lambda)$ considered in Lemma \ref{coker_basis_for_surface_factor}. \end{corollary}

\begin{proof}
This follows from Lemma \ref{coker_basis_for_surface_factor} and the definition of $\Pi_{V_F}$.
\end{proof}

%Note that Corollary \ref{t=1} gives another proof of Corollary \ref{torelliextends} in the setting of pushed-in Seifert surfaces. %We now use Corollary \ref{t=1} to exhibit explicit Dehn twists within $\mathcal{E}^{\mathrm{TOP}}(D^4,F)$.  %This is done by constructing explicit unitary transvections in $\mathrm{Iso}(\lambda_{V_F})$. 
%We give some example computations for Seifert surfaces $F \subset S^3$ below.

\subsection{Extending mapping classes with explicit transvections}

We now use the computation of $\lambda_{V_F}$ (Proposition \ref{basisandform}) to write down isometries of $\lambda_{V_F}$. We then use these isometries to exhibit certain Dehn twists as topologically extendable mapping classes. 

\begin{defn}
Let $H$ be a free $\Lambda$-module equipped with a sesquilinear Hermitian pairing $\lambda : H \times H \to \Lambda$. We say $\lambda$ is divisible by $r \in \Lambda$ if for any $x,y\in H$, there exists some $s \in \Lambda$ (depending on $x$ and $y$) such that $\lambda(x,y) = rs$. 
\end{defn}

\begin{rmk} \label{equivariantformisdivisible}
By Proposition \ref{basisandform}, $ \lambda_{V_F}$ is divisible by $1-t$. Indeed, in the basis $\mathcal{C} = \{\Sigma_i\}_{i=1}^{2g}$, we have $\lambda = (1-t)P$ where $P = S^T-t^{-1}S$. 
\end{rmk}

%Relative to our basis $\mathcal{C}$ for $H_2(V_F;\Lambda)$,  $$\mathrm{Isom}(\lambda_{V_F}) = \{ A \in GL_{2g}(\Lambda) \; | \; A^TQ \bar A = Q\}.$$ We have $Q= (1-t)P$ where $P = S^T - t^{-1}S$. Any matrix $A \in \mathrm{GL}_{2g}(\Lambda)$ with $A^TQ \bar A = Q$ must satisfy \begin{equation}
    %A^T P \bar A = P . \label{PmatrixAeqn}
%\end{equation} %Note that $P \pmod{t-1}$ is the Gram matrix of the intersection form on $H_1(F;\Z)$ in our basis of curves $a_i \subset F$, hence $A$ has image in $\mathrm{Sp}_{2g}\Z$ after setting $t=1$.

%By Proposition \ref{basisandform}, the equivariant intersection form $\lambda = \lambda_{V_F}$ on the exterior $V_F$ of a pushed-in Seifert surface $F \subset D^4$ is divisible by $1-t$. Indeed, in the basis $\mathcal{C} = \{\Sigma_i\}_{i=1}^{2g}$, we have $\lambda = (1-t)P$ where $P = S^T-t^{-1}S$ where $S$ is the Seifert form of $F$ in the basis $\{a_i\}_{i=1}^{2g}$ for $H_1(F;\mathbb{Z})$. 

\begin{proposition} \label{transvectionformulae}
Let $H$ be a free $\Lambda$-module equipped with a sesquilinear Hermitian pairing $\lambda : H \times H \to \Lambda$. Fix $v \in H$, $p \in \Lambda$. Consider the transvection $T_{v,p} : H \to H$ defined by 
\[ T_{v,p}(x) = x + p \lambda(x,v)v.\]
%If $p + \overline{p} =0$ and $\lambda(v,v) = 0$, then $T_{v,p} \in \mathrm{Iso}(H,\lambda)$. 
We have the following properties. 
\begin{enumerate}
    \item If $p +\overline{p}+p\overline{p}\lambda(v,v) = 0$, then $\lambda(T_{v,p}(x),T_{v,p}(y)) = \lambda(x,y)$ for all $x,y\in H$. If $\lambda(v,v) =0$ and $p + \overline{p} = 0$, then $T_{v,p} \in  \mathrm{Isom}(H,\lambda)$, with inverse $T_{v,-p}$. 
    
    \item If $\lambda$ is divisible by $1-t$, we may allow $p = s/(1-t)$ with $s\in \Lambda$ in the definition of $T_{v,p}$ and in the item above (part (a)). 

    \item Assume that $\lambda$ is divisible by $(1-t)$ and $\lambda(v,v) = \pm(1-t)(1-t^{-1})$. Let $p = \pm t/(1-t)$. Then $p + \overline{p} + p\overline{p}\lambda(v,v) = 0$ and $T_{v,\frac{\mp1}{1-t}}$ is the inverse of $T_{v,p}$. Therefore $T_{v,p} \in \mathrm{Isom}(H,\lambda)$. 
    
\end{enumerate}

%If $p +\overline{p}+p\overline{p}\lambda(v,v) = 0$, then $T_{v,p} \in \mathrm{Iso}(H,\lambda)$. 
\end{proposition}

  %Ind In particular, for $\lambda_{V_F}$, we may allow $ p = s/(1-t)$ for $s \in \Lambda$ in Proposition \ref{transvectionformulae}. 

\begin{proposition} \label{explicittransvections}
Let $F \subset D^4$ be a pushed in Seifert surface with exterior $V_F = D^4-\nu F$. Fix a simple closed curve $a \subset F$. Assume that $\mathrm{lk}(a,a^+) = 0$ or $\pm1$. Let $\Sigma_a \subset \widetilde{V_F}$ denote a 2-cycle obtained by doubling along push-offs of $a\subset F$ as in Proposition \ref{basisandform}. Let 
\[ p_a =  \begin{cases} \displaystyle\frac{t+1}{1-t} & \mathrm{if   }\; \mathrm{lk}(a,a^+) = 0,\\ 
\displaystyle\frac{\pm t}{1-t} & \mathrm{if }\; \mathrm{lk}(a,a^+)  =\pm 1. \end{cases} \]
Then $T_{\Sigma_a,p_a} \in \mathrm{Isom}(\lambda_{V_F})$. %Recall the boundary (Definition \ref{boundaryofiso}) of $T_{\Sigma_a,p_a}$ \[\Pi_{V_F}(T_{\Sigma_a,p_a}) : H_1(\partial V_F;\Lambda) \to H_1(\partial V_F;\Lambda)\]  gives an isometry of the Blanchfield pairing on $\partial V_F$, that is $  \Pi_{V_F}(T_{\Sigma_a,p_a})\in \mathrm{Iso}(\mathrm{Bl}_{\partial V_F})$ (Proposition \ref{boundaryiso}) and we have a decomposition of the isometry group \[\mathrm{Iso}(\mathrm{Bl}_{\partial V_F}) \cong \mathrm{Sp}(H_1(F;\mathbb{Z})) \oplus \mathrm{Iso}(\mathrm{Bl}_{E_K})\]  Under this identification 
Moreover, the restriction $\Pi_{V_F}(T_{\Sigma_a,p_a})\rvert_{H_1(F\times S^1;\Lambda)}$ is the ordinary transvection $T_{[a]}^{\mathrm{lk}(a,a^+)} \in \mathrm{Sp}(H_1(F;\mathbb{Z}))$ along $[a] \in H_1(F;\mathbb{Z})$ if $\mathrm{lk}(a,a^+) = \pm1$ and the square $T_{[a]}^2$ if $\mathrm{lk}(a,a^+) = 0$. On the other hand, $\Pi_{V_F}(T_{\Sigma_a,p_a})\rvert_{H_1(E_K;\Lambda)}$ is the identity in either case. 
%the\textcolor{red}{Compute $\Pi_{V_F}(T_{\Sigma_a,p_a})\rvert_{H_1(E_K;\Lambda)}$ by identifying the $H_1(E_K;\Lambda)$ as a further quotient of $\mathrm{coker}(\lambda_{V_F})$. In the notation of Lemma \ref{coker_basis_for_surface_factor} this is $\mathrm{coker}(\lambda_{V_F})/M$? Conclude extendability as in proof of Corollary \ref{torelliextends}}. 
\end{proposition}

\begin{proof}
That $T_{\Sigma_a,p_a} \in \mathrm{Isom}(\lambda_{V_F})$ follows from Proposition \ref{transvectionformulae}. Indeed, if $\mathrm{lk}(a,a^+) =0$, then $\lambda_{V_F}(\Sigma_a,\Sigma_a) =0$ and $p_a + \overline{p_a} = 0$. If $\mathrm{lk}(a,a^+) = \pm 1$, then $\lambda(\Sigma_a,\Sigma_a) = \pm (1-t)(1-t^{-1})$ and the claim follows from part (c) of Proposition \ref{transvectionformulae}. 

Recall that $U \subset \mathrm{coker}(\widehat\lambda_{V_F})$ is the $\mathbb{Z}$-submodule spanned by $\lambda^{\delta}(-,\Sigma_i^+)$ for $i=1,\ldots,2g$. Under the isomorphism \[D_{V_F}\colon\mathrm{coker}(\widehat\lambda_{V_F}) \to H_1(\partial V_F;\Lambda)\] the submodule $U$ corresponds to the $H_1(F\times S^1;\Lambda)\cong H_1(F;\mathbb{Z})$ factor. By definition of $\Pi_{V_F}$, we must compute the action $T_{\Sigma_a,p_a}^{-*}$ on $U$. 
In the notation of Lemma \ref{coker_basis_for_surface_factor}, we have, first working on the level of the dual module $H^*$,  
\begin{align*}
    (1-t)T_{\Sigma_a, p_a}^{-*} ( \lambda^{\delta}(-,\Sigma_i^+)) &= \lambda(T_{\Sigma_a, p_a}^{-1}(-), \Sigma_i) \\
    &= \lambda(-, T_{\Sigma_a, p_a}(\Sigma_i)) \\
    &= \lambda(-, \Sigma_i + p_a\lambda(\Sigma_i, \Sigma_a) \Sigma_a) \\
    &= \lambda(-,\Sigma_i) + \overline{p_a \; \lambda(\Sigma_i, \Sigma_a)} \lambda (-, \Sigma_a) \\
    &= (1-t)(\lambda^{\delta}(-,\Sigma_i^+) + \overline{p_a \; \lambda(\Sigma_i, \Sigma_a)} \lambda^{\delta}(-, \Sigma_a^+)),
\end{align*} and hence $$T_{\Sigma_a, p_a}^{-*} ( \lambda^{\delta}(-,\Sigma_i^+)) = \lambda^{\delta}(-,\Sigma_i^+) + \overline{p_a \; \lambda(\Sigma_i, \Sigma_a)} \lambda^{\delta}(-, \Sigma_a^+).$$ The connecting map $\delta \colon H_2(V_F, \partial V_F;\Lambda) \rightarrow H_1(\partial V_F ; \Lambda)$ sends $\Sigma_i^+$  to $[a_i]$, and $\Sigma_a^+$ to $[a]$ (Lemma  \ref{connectinghomomorphism}). It follows that 
\[ \Pi_{V_F}(T_{\Sigma_a,p_a})[a_i] = [a_i] +\overline{p_a\lambda(\Sigma_i,\Sigma_a)}[a].\]

By (the proof of) Proposition \ref{basisandform}, 
\[\lambda(\Sigma_i,\Sigma_a) = (1-t)(\mathrm{lk}(a,a_i^+)-t^{-1}\mathrm{lk}(a_i,a^+)). \]

Therefore 
\[ p_a \lambda(\Sigma_i,\Sigma_a) = \begin{cases} 
(t+1) (\mathrm{lk}(a,a_i^+)-t^{-1}\mathrm{lk}(a_i,a^+)) & \mathrm{if   }\; \mathrm{lk}(a,a^+) = 0,\\ 
\pm \,t (\mathrm{lk}(a,a_i^+)-t^{-1}\mathrm{lk}(a_i,a^+)) & \mathrm{if }\; \mathrm{lk}(a,a^+)  =\pm 1. \end{cases} 
\]

 Recall that $t = 1$ in the $H_1(F\times S^1;\Lambda)$ factor, and that, for the Seifert form $S$, the difference $S^T-S$ is the oriented intersection matrix of the curves $\{a_i\}_{i=1}^{2g}$ \cite[Chapter 8, Exercise 7]{Rolfsen}.  Therefore 
 \[ \mathrm{lk}(a,a_i^+) - \mathrm{lk}(a_i,a^+) = \langle a_i,a\rangle\] 
 where $\langle-,-\rangle$ denotes the oriented intersection pairing on $H_1(F;\mathbb{Z})$. It follows that
 
 \[ \Pi_{V_F}(T_{\Sigma_a,p_a})[a_i] = \begin{cases} [a_i] + 2\langle a_i,a\rangle [a] &  \mathrm{if   }\; \mathrm{lk}(a,a^+) = 0,\\ 
 
 [a_i] \pm \langle a_i,a\rangle[a] &  \mathrm{if }\; \mathrm{lk}(a,a^+)  =\pm 1.  \end{cases}
 \]
 
 This is the desired (power of a) transvection in either case. 

%\[\frac{1}{1-t}\lambda(\Sigma_a, \Sigma_i) \equiv \langle a, a_i\rangle \pmod{t-1}\] for all $j=1,..,2g$ where $\langle -,-\rangle \colon H_1(F;\Z) \times H_1(F;\Z) \rightarrow \Z$ denotes the intersection form of $F$ over $\Z$. 

%By Corollary \ref{t=1}, the restriction $\Pi_{V_F}(T_{\Sigma_a,p_a})\rvert_{H_1(F\times S^1;\Lambda)}$ is obtained by computing $T_{\Sigma_a,p_a}$ in the basis $\mathcal{C}$ and setting $t=1$. It follows that 

 % \[ \Pi_{V_F}(T_{\Sigma_a,p_a})([a_i])  = \begin{cases} [a_i] + (t+1)\langle a_i,a\rangle[a] & \mathrm{if   }\; \mathrm{lk}(a,a^+) = 0,\\ 
%[a_i] \pm t \langle a_i,a\rangle[a] & \mathrm{if }\; \mathrm{lk}(a,a^+)  =\pm 1.\end{cases}\] 

%\[ \Pi_{V_F}(T_{\Sigma_a,p_a})([a_i])  = \begin{cases} [a_i] + (t+1)\langle a_i,a\rangle[a] & \mathrm{if   }\; \mathrm{lk}(a,a^+) = 0,\\ 
%[a_i] \pm t \langle a_i,a\rangle[a] & \mathrm{if }\; \mathrm{lk}(a,a^+)  =\pm 1.\end{cases}\]

%At $t=1$, we have the desired (power of a) transvection in either case.

To see that $\Pi_{V_F}(T_{\Sigma_a,p_a})$ acts trivially on $H_1(E_K;\Lambda)$, use 
\[ \mathrm{coker}(\widehat\lambda_{V_F}) \cong H_1(\partial V_F;\Lambda) \cong H_1(F\times S^1;\Lambda) \oplus H_1(E_K;\Lambda)\]
to exhibit $H_1(E_K;\Lambda)$ as the quotient $\mathrm{coker}(\widehat\lambda_{V_F})/U$. By naturality, it suffices to show that $T_{\Sigma_a,p_a}^*$ acts trivially on this quotient. For $\varphi \in H^*$, 
\[ T_{\Sigma_a,p_a}^*(\varphi) = \varphi + p_a \lambda(-,\Sigma_a)\varphi(\Sigma_a) = \varphi + p_a (1-t)\lambda^{\delta}(-,\Sigma_a^+)\varphi(\Sigma_a).\]
Now $[\lambda^{\delta}(-,\Sigma_a^+)] = 0$ in $\mathrm{coker}(\widehat \lambda_{V_F})/U$, so $[T_{\Sigma_a,p_a}^*(\varphi)] = [\varphi] \in \mathrm{coker}(\widehat\lambda_{V_F})/U$.

\end{proof}

We now generalize the Hopf band trick (Proposition \ref{HBtrick}) and the \textit{tube isotopy}. The latter realizes the square of a Dehn twist along a curve bounding a zero-framed disk in the surface exterior as extendable (see \cite[Figure 1]{HiroseP2}). Note that both of these techniques require the presence of an embedded disk in the surface exterior.

\begin{corollary}  \label{GeneralHBtrick}
Let $F \subset D^4$ be a pushed in Seifert surface for a knot. If $\gamma \subset F$ is a simple closed curve with $\mathrm{lk}(\gamma,\gamma^+) = \pm1$, then $T_{\gamma} \in \mathcal{E}^{\mathrm{TOP}}(D^4,F)$. If $\gamma \subset F$ is a simple closed curve with $\mathrm{lk}(\gamma,\gamma^+) = 0$, then $T_{\gamma}^2 \in \mathcal{E}^{\mathrm{TOP}}(D^4,F)$.
\end{corollary}

\begin{proof}
Proceed as in the proof of Corollary \ref{torelliextends} using the explicit transvections constructed in Proposition \ref{explicittransvections}.
\end{proof}

Corollary \ref{GeneralHBtrick} does not address curves $\gamma \subset F$ with Seifert framing $|\mathrm{lk}(\gamma,\gamma^+)| > 1$. In general, no claim can be made about such curves. Indeed, the spin-flexibility of unknots (Proposition \ref{spinflexunknot}) implies that there exist curves with any prescribed odd Seifert framing whose associated twists extend. On the other hand, we show in Proposition \ref{non_extendable_twist} that there exist curves with odd Seifert framing for which the associated twists do not extend. As Rochlin's form $q_F : H_1(F;\mathbb{Z}/2\mathbb{Z}) \to \mathbb{Z}/2\mathbb{Z}$ is the mod 2 Seifert form, this is an obstruction to spin-flexibility.

\section{Obstructing extendable mapping classes via the equivariant intersection form}
\label{obstruct}
%\textcolor{red}{Simplifying some of this section}

%\textcolor{red}{I'm trying to rephrase some of the earlier bits of this section to be a little more streamlined. }

In this section, we demonstrate how the equivariant intersection form on the surface exterior can be used to obstruct extendable mapping classes.  %results of Conway--Powell to obtain obstructions to extending mapping classes of a surface to ambient homeomorphisms of the $4$-manifold. 
We apply our obstruction to  %study the extendable subgroup for 
Milnor fibers of isolated hypersurface singularities of type ADE. % case of smooth curves in $\C^2$ intersected with a $4$-ball, where the curves are obtained from the ADE classification of isolated singularities. 
%We highlight this case since 
In this setting, the image of the extendable subgroup under the mod 2 symplectic representation is precisely the associated mod 2 Weyl group.

%\subsection{Extendable mapping classes and isometries of the Blanchfield pairing}

\subsection{Obstructing through the equivariant intersection form}

We begin with a general obstruction to extending a mapping class of a pushed in Seifert surface in $D^4$. We use notation introduced in Section 4. Let $F \subset D^4$ be a genus $g$ Seifert surface for a knot $K \subset S^3$ with its interior pushed into the interior of $D^4$. Let $S$ be the Seifert matrix of $F$ with respect to a fixed basis of curves $\{a_i\}_{i=1}^{2g} \subset F$ for $H_1(F;\Z)$. Let $V_{F} = D^4 - \mathrm{Int}\; \nu(F)$. Recall that the matrix $Q$ of the equivariant intersection form on $H_2(V_{F};\Lambda)$ in the $\Lambda$-basis $\mathcal{C}$ obtained from $\{a_i\}_{i=1}^{2g}$ (Proposition \ref{basisandform}) is $Q = (1-t)P$, where $P = S^T - t^{-1}S$. For any matrix $A$ with entries in $\Lambda$, write $A(\pm 1)$ for the integer matrices $A \pmod{t \mp 1}$. The basis $\{a_i\}_{i=1}^{2g}$ determines a mod 2 symplectic representation $\Psi_2 : \Mod(F) \to \mathrm{Sp}_{2g}(\Z_2)$.

\begin{lemma} \label{thm:obtruction_for_pushed_in_Seifert_surfaces}
     Let $$G_P = \{ w \in \mathrm{GL}_{2g}(\Z) \;  | \; w^T P(-1) w = P(-1) \},$$ where $P(-1) = S^T + S$. If $f \in \mathcal{E}^{\mathrm{TOP}}(D^4,F)$, then $\Psi_2(f)$ lies in the mod 2 reduction of $G_P$.
\end{lemma}
\begin{proof}
    Suppose $f:D^4 \rightarrow D^4$ is an orientation-preserving homeomorphism fixing $F$ setwise and $\partial D^4$ pointwise. %Then by Lemmas \ref{isotopy_tofix_tube}, \ref{isotoped_map_intertwines}, and \ref{isotoped_action_Blanchfield}, 
    Up to isotopy $\mathrm{rel} \; \partial D^4 \cup F$, $f$ restricts to a bundle map on a tubular neighbourhood $\nu(F)$ of $F$ in $D^4$, and the restriction $f\rvert_{V_F}$ lifts to a homeomorphism $\widetilde{f}$ of the universal cover $\widetilde{V_F}$. This % Any such map lifts to a homeomorphism $\hat \varphi: \widetilde{V_F} \rightarrow \widetilde{V_F}$, which 
    induces an isometry $\alpha_f=\widetilde{f}_*$ of the equivariant intersection form $(H_2(V_F;\Lambda), Q)$. The restriction $f\rvert_{\partial V_F} : \partial V_F \to \partial V_F$ intertwines the inclusion map $i : \pi_1(\partial V_F) \to \mathbb{Z}$. In particular, it induces an automorphism $(f\rvert_{\partial V_F})_* : H_1(\partial V_F;\Lambda) \to H_1(\partial V_F;\Lambda)$. By (the proof of) \cite[Proposition 3.11]{Conway_Powell}, the following diagram is commutative
    \[\begin{tikzcd}
	{\mathrm{coker}(\widehat\lambda_{V_F}) } & {\mathrm{coker}(\widehat\lambda_{V_F})} \\
	{H_1(\partial V_F;\Lambda)} & {H_1(\partial V_F;\Lambda)}
	\arrow["{[\alpha_f^{-*}]}", from=1-1, to=1-2]
	\arrow["{D_{V_F}}"', from=1-1, to=2-1]
	\arrow["{D_{V_F}}", from=1-2, to=2-2]
	\arrow["{(f\rvert_{\partial V_F})_*}"', from=2-1, to=2-2]
\end{tikzcd}\]

    Therefore $\Pi_{V_F}(\alpha_f) = (f\rvert_{\partial V_F})_*$. By Proposition \ref{boundaryofisopreservesblanch}, $(f\rvert_{\partial V_F})_* \in \mathrm{Isom}(\mathrm{Bl}_{\partial V_F})$. Recall that $\partial V_F = F \times S^1 \cup_\partial E_K$ and Mayer-Vietoris gives an isomorphism $H_1(\partial V_F;\Lambda) \cong H_1(E_K;\Lambda) \oplus H_1(F \times S^1;\Lambda)$. By Proposition \ref{isoofblanchfield}, 
    \[ \mathrm{Isom}(\mathrm{Bl}_{\partial V_F}) \cong \mathrm{Isom}(\mathrm{Bl}_{F \times S^1})  \oplus \mathrm{Isom}(\mathrm{Bl}_{E_K})\]
and $\mathrm{Isom}(\mathrm{Bl}_{F\times S^1}) \cong \mathrm{Sp}(H_1(F;\mathbb{Z}))$. Write $(f\rvert_{\partial V_F})_* = h_F \oplus h_K$ for $h_F \in \mathrm{Isom}(\mathrm{Bl}_{F\times S^1})$ and $h_K \in \mathrm{Isom}(\mathrm{Bl}_{E_K})$. Recall that the restriction of $\widetilde{V_F}$ to the $F \times S^1$ part of $\partial V_F$ is the covering $p: F \times \mathbb{R} \to F \times S^1$. It follows that the automorphism $h_F$ fits into the following commutative diagram, 
% https://q.uiver.app/#q=WzAsNixbMCwwLCIgSF8xKEZcXHRpbWVzIFxcbWF0aGJie1J9O1xcbWF0aGJie1p9KSJdLFsxLDAsIkhfMShGXFx0aW1lcyBTXjE7XFxtYXRoYmJ7Wn0pIl0sWzIsMCwiSF8xKEZcXHRpbWVzIEReMjtcXG1hdGhiYntafSkiXSxbMCwxLCJIXzEoRlxcdGltZXMgXFxtYXRoYmJ7Un07XFxtYXRoYmJ7Wn0pIl0sWzEsMSwiSF8xKEZcXHRpbWVzIFNeMTtcXG1hdGhiYntafSkiXSxbMiwxLCJIXzEoRiBcXHRpbWVzIEReMjtcXG1hdGhiYntafSkiXSxbMCwxLCJwXyoiXSxbMSwyLCJ7aV9GfV8qIl0sWzAsMywiaF9GIiwyXSxbMyw0LCJwXyoiLDJdLFs0LDUsIntpX0Z9XyoiLDJdLFsyLDUsIihmXFxydmVydF97XFxudSBGfSlfKiJdXQ==
\[\begin{tikzcd}
	{ H_1(F\times \mathbb{R};\mathbb{Z})} & {H_1(F\times S^1;\mathbb{Z})} & {H_1(F\times D^2;\mathbb{Z})} \\
	{H_1(F\times \mathbb{R};\mathbb{Z})} & {H_1(F\times S^1;\mathbb{Z})} & {H_1(F \times D^2;\mathbb{Z})}
	\arrow["{p_*}", from=1-1, to=1-2]
	\arrow["{h_F}"', from=1-1, to=2-1]
	\arrow["{{i_F}_*}", from=1-2, to=1-3]
	\arrow["{(f\rvert_{\nu F})_*}", from=1-3, to=2-3]
	\arrow["{p_*}"', from=2-1, to=2-2]
	\arrow["{{i_F}_*}"', from=2-2, to=2-3]
\end{tikzcd}\]
where $i_F : F \times S^1 \to F \times D^2$ is inclusion into the closed tubular neighbourhood $\nu F\cong F \times D^2$. As ${i_F}_* \circ p_*$ is an isomorphism, the induced action of $(f\rvert_{\partial V_F})_*$ on the $H_1(F\times S^1;\Lambda)$ factor is $(f\rvert_F)_*$.

   % Recall that $\partial V_F = F \times S^1 \cup E_K$. We consider the induced action on $H_1(F \times S^1;\Lambda)$ in the decomposition $H_1(\partial V_F;\Lambda) \cong H_1(E_K;\Lambda) \oplus H_1(F \times S^1;\Lambda)$ induced by inclusion. By Lemma \ref{isotoped_action_Blanchfield}, the induced action of $\hat f_*$ on $H_1(F\times S^1;\Lambda)$ is $(f|_F)_*$. 
   
  On the other hand, we have computed the action of isometries of $\lambda_{V_F}$ on the $H_1(F\times S^1;\Lambda)$ factor explicitly. In terms of our bases, the map \[\Pi_{V_F}: \mathrm{Isom}(H_2(V_F;\Lambda),Q) \rightarrow \mathrm{Isom}(\mathrm{Bl}_{\partial V_F})\] is given by $A \mapsto A \pmod{t-1} =: A(1)$ on the $H_1(F \times S^1;\Lambda)$ factor (Corollary \ref{t=1}). It follows that if a mapping class $f\colon F \to F$ is extendable, its image in $\mathrm{Sp}_{2g}(\mathbb{Z})$ determined by $f_* :H_1(F;\mathbb{Z}) \to H_1(F;\mathbb{Z})$ and $\{a_i\}_{i=1}^{2g}$ is the specialization of a matrix $A  \in \mathrm{Isom}(H_2(V_F;\Lambda), Q) \leq \mathrm{GL}_{2g}(\Lambda)$ to $t=1$. 
  
  The claim follows from specializing $A$ to $t = -1$. Indeed, if $A \in \mathrm{Isom}(H_2(V_F;\Lambda), Q)$, then $A^TP \bar A = P$, so we have $A(-1)^TP(-1) A(-1) = P(-1)$, hence $A(-1) \in G_P$. Now $A(-1) \equiv A(1) \pmod{2}$ and $\Psi_2(f) = A(1) \mod 2$. 
\end{proof}
\subsection{Application to Milnor fibers of ADE singularities}
Let $f \in \C[x,y]$ be the germ of a simple isolated singularity at $0 \in \C^2$. Let $D^4(\epsilon)$ be a ball of radius $\epsilon$ centered at $0 \in \C^2$, and $D_{\delta} = \{z \in \C \,|\, |z|\leq \delta\}$. For $\epsilon >0$ sufficiently small, the intersection $L = f^{-1}(0) \cap \partial D^4(\epsilon)$ is transverse and defines the \textit{link} of the singularity \cite[Corollary 2.9]{milnor_complex_hypersurfaces}. Moreover, for $0 < \delta \ll \epsilon$, the restriction 
\[ f\rvert_{f^{-1}(\partial D_\delta) \cap D^4(\epsilon)} : f^{-1}(\partial D_\delta) \cap D^4(\epsilon) \to \partial D_\delta\]
is a fiber bundle \cite{milnor_complex_hypersurfaces}. The fiber of  $f\rvert_{f^{-1}(\partial D_\delta) \cap D^4(\epsilon)}$ above a point of $\partial D_\delta$ is the \textit{Milnor fiber} $\Sigma_f$ of $f$. The Milnor fiber of $f$ is a compact, oriented surface with non-empty boundary.  

\begin{proposition} \label{MilnorfiberisZsurface}
The Milnor fiber of an isolated hypersurface singularity $f$ is a $\mathbb{Z}$-surface in a 4-ball. In fact, it is a pushed in Seifert surface. 
\end{proposition}

\begin{proof}
Milnor's fibration theorem \cite[Theorem 4.8]{milnor_complex_hypersurfaces} asserts that the link $L \subset \partial D^4(\epsilon)$ is fibered by the map \[\frac{f}{|f|}\colon \partial D^4(\epsilon) -L \rightarrow S^1.\] The fiber $F_\theta$ above $\theta \in S^1$ is the interior of a (minimal genus) Seifert surface for $L$ \cite[Lemma 6.1]{milnor_complex_hypersurfaces}. The proof of \cite[Theorem  5.11]{milnor_complex_hypersurfaces} shows that $\Sigma_f = f^{-1}(c) \cap D^4(\epsilon)$ for $c \in \partial D_\delta$ is obtained by pushing  $F_\theta$ (for $\theta = c/|c|$) into the interior of $D^4(\epsilon)$. Finally, pushed in Seifert surfaces have knot group $\mathbb{Z}$ (Proposition \ref{pushedinareZsurfaces}). 
\end{proof}

%such that $L = f^{-1}(0) \cap \partial D^4(\epsilon)$ defines a transverse intersection. Moreover, the link $L \subset \partial D^4(\epsilon)$ is fibered by the map \[\frac{f}{|f|}\colon \partial D^4(\epsilon) -L \rightarrow S^1,\] hence its fiber above $\theta \in S^1$ is a (minimal genus) Seifert surface for $L$. Moreover, for $\epsilon >0$ sufficiently small and $\delta >0$ small with respect to $\epsilon$, the intersection \(f^{-1}(\delta) \cap D^4(\epsilon)\) defines a properly embedded surface in $D^4_\epsilon$. 

%\begin{defn}
%The \textit{Milnor fiber} of $f$ is \(f^{-1}(\theta) \cap D^4(\epsilon)\). 
%\end{defn}

%The surface \(f^{-1}(\theta) \cap D^4(\epsilon)\) is obtained from the fiber surface by pushing its interior into the interior of $D^4(\epsilon)$ \cite[Proof of Lemma 5.10]{milnor_complex_hypersurfaces}. We refer to such a surface as a \textit{Milnor fiber} for the link $L$ of the isolated singularity.

%\begin{corollary}[c.f. {\cite[Lemma 5.10]{milnor_complex_hypersurfaces}}] \label{MilnorfiberisZsurface}
%The Milnor fiber of an isolated hypersurface singularity is a $\mathbb{Z}$-surface in a 4-ball.   
%\end{corollary}

%we note that the Milnor fiber of an isolated hypersurface singularity is a $\mathbb{Z}$-surface in a 4-ball.  

%\begin{corollary}
%The Milnor fiber of an isolated hypersurface singularity is a $\mathbb{Z}$-surface in a 4-ball.   
%\end{corollary}

%\begin{proof}
%%\end{proof}

Simple isolated singularities are classified into ADE type, and are given by $$(A_n) \; f(x,y) = x^{2} + y^{n+1}, \qquad (D_n) \; f(x,y) = xy^2 + x^{n-1},$$
and the exceptional types $$(E_6) \; f(x,y) = x^3+y^4, \qquad (E_7) \; f(x,y) = x^3 + xy^3, \qquad (E_8) \; f(x,y) = x^3 + y^5.$$
For each family, the following table describes the link of the singularity as the closure of a positive braid (see \cite[Example 2.5]{Casals_lagrangian_skeleta}).

\begin{table}[h] \label{positivebraids}
\begin{tabular}{|l|l|}
\hline
Type & Positive Braid                                                                                                                                  \\ \hline
$A_n$ & $\sigma_1^{n+1}$ \\ \hline
$D_n$ & $\sigma_1^{n-2} \sigma_2 \sigma_1^2 \sigma_2$ \\ \hline
$E_6$ & $\sigma_1^3 \sigma_2 \sigma_1^3 \sigma_2$ \\ \hline
$E_7$ & $\sigma_1^4 \sigma_2 \sigma_1^3 \sigma_2$    \\ \hline
$E_8$ & $\sigma_1^5 \sigma_2 \sigma_1^3 \sigma_2$   \\ \hline
\end{tabular}
\caption{The links of singularities of type ADE as the closure of positive braids.}
\end{table}

For each singularity $\Gamma = A_n,D_n, E_6,E_7$ or $E_8$, the Milnor fiber $F_\Gamma$ is, up to isotopy, the \textit{standard Seifert surface} associated with the braid word (see \cite[Section 2]{rudolph} for a precise definition of the standard Seifert surface associated to a braid word). Depicted in Figure \ref{fig:non_extendable_curves} is $F_\Gamma$ for $\Gamma = A_8$. The Seifert surface $F_\Gamma$ is equipped with a distinguished ordered configuration of curves $\mathcal{C}_\Gamma$ that give a basis for $H_1(F_\Gamma;\Z)$ (c.f. \cite[Figure 1]{ferretti}) and have \textit{intersection graph} (Definition \ref{intersectiongraphofnetwork}) given by the corresponding ADE Dynkin diagram. We abuse notation and denote the Dynkin diagram associated to the singularity by $\Gamma$.  Let $S$ denote the Seifert matrix for $F_\Gamma$ in the basis $C_\Gamma$. Then $P(-1) = S + S^T$ is the \textit{Cartan matrix} associated to the Dynkin diagram (see \cite[Section 2]{Baader_pos_braids}). 

\begin{rmk} \label{diskshb}
The explicit construction of $F_\Gamma$ and $C_\Gamma$ we use exhibits the curves in $C_\Gamma$ as unknots in an $S^3$ slice of a collar about $\partial D^4$ with Seifert framing $-1$. In particular, with our conventions $S+S^T$ has diagonal entries $-2$, and so we differ from the usual Cartan matrix associated to $\Gamma$ by a minus sign. 
\end{rmk}

Lemma \ref{thm:obtruction_for_pushed_in_Seifert_surfaces} gives the following.

\setcounter{maintheorem}{4}
\begin{maintheorem}
 Let $\Gamma$ be a Coxeter--Dynkin diagram of type $A_n$ with $n$ even, $E_6$ or $E_8$.
    Let $W_2(\Gamma)$ be the subgroup of $\mathrm{Sp}(H_1(F_\Gamma; \Z_2))$ generated by the transvections in the homology classes of the curves in $\mathcal{C}_\Gamma$. This is the mod 2 reduction of the Weyl group of type $\Gamma$. Then 
    \[ \Psi_2(\mathcal{E}^{\mathrm{TOP}}(D^4,F_\Gamma)) = W_2(\Gamma) \]
    where $\Psi_2 : \Mod(F_\Gamma) \to \mathrm{Sp}(H_1(F_\Gamma; \Z_2))$ is the mod 2 symplectic representation. %Suppose $f \in \epsilon^{\mathrm{DIFF}}( D^4, F)$, then $f_* \in W(\Gamma)$.   
\end{maintheorem}

\begin{proof}
   We begin with the containment  $\Psi_2(\mathcal{E}^{\mathrm{TOP}}(D^4,F_\Gamma)) \leq W_2(\Gamma)$.  By Lemma \ref{thm:obtruction_for_pushed_in_Seifert_surfaces}, $\Psi_2(\mathcal{E}^{\mathrm{TOP}}(D^4,F_\Gamma))$ is contained in the mod 2 reduction of 
   \[ G_P = \{ w \in \mathrm{GL}_{2g}(\Z) \;  | \; w^T P(-1) w = P(-1) \},\]
   where $P(-1) = S + S^T$ is the Cartan matrix associated to the Coxeter-Dynkin diagram $\Gamma$. The group $G_P$ is the isometry group of the root lattice of type $\Gamma$. We now show the mod 2 reduction of $G_P$ is precisely the mod 2 Weyl group $W_2(\Gamma)$. 
%We have that $P(-1) = S + S^T$ is the Cartan matrix associated to the Dynkin diagram of type $\Gamma$. For a matrix $A$ representing an isometry of the equivariant intersection form $\lambda_{V_F}$ in the basis $\mathcal{C}$,  we have the specialization of \eqref{PmatrixAeqn} to $t = -1$,  \[A(-1)^T (S + S^T) A(-1) = S+S^T.\] Therefore $A(-1)$ is an isometry of the root lattice $L'$ of type $\Gamma$. By Lemma \ref{thm:obtruction_for_pushed_in_Seifert_surfaces}, $\Psi_2(\mathcal{E}^{\mathrm{TOP}}(D^4,F)$ is contained in the mod 2 reduction of the isometry group of the root lattice $L'$. 
   Recall that \[G_P \cong W(\Gamma) \rtimes \mathrm{Aut}(\Gamma)\] where $\mathrm{Aut}(\Gamma)$ denotes the automorphism group of the graph $\Gamma$ and $W(\Gamma)$ the Weyl group of $\Gamma$. For $\Gamma$ of type $A_n$ or $E_6$, we have $\mathrm{Aut}(\Gamma) = \langle \iota \rangle$, where $\iota$ is an involution; in these cases, for the longest element $\omega_0 \in W(\Gamma)$ with respect to the simple reflections determined by $\Gamma$, we have that $\iota = - \omega_0$ as an isometry of the lattice. For $\Gamma = E_8$, we have that $\mathrm{Aut}(\Gamma) = 1$. Hence the mod $2$ reduction of the isometry group of the lattice agrees with the mod $2$ reduction of the Weyl group $W(\Gamma)$.  %It now follows from Lemma \ref{thm:obtruction_for_pushed_in_Seifert_surfaces} that $\Psi_2(\mathcal{E}^{\mathrm{TOP}}(D^4,F)) \leq W(\Gamma)$. 

    We now prove that $W_2(\Gamma) \leq \Psi_2(\mathcal{E}^{\mathrm{TOP}}(D^4,F_\Gamma))$. The curves in $C_\Gamma$ all bound disks in the exterior of $F_\Gamma \subset D^4$ with relative Euler number $-1$ (Remark \ref{diskshb}). %; see \cite[Example 2.5]{Casals_lagrangian_skeleta}. 
    The Hopf band trick (Proposition \ref{HBtrick}) implies that \[ \mathcal{T}_{C_\Gamma} = \langle T_c \,:\, c\in C_\Gamma \rangle \leq  \mathcal{E}^{\mathrm{DIFF}}(D^4,F_\Gamma) \leq \mathcal{E}^{\mathrm{TOP}}(D^4, F_{\Gamma}).\] 
    Therefore $W_2(\Gamma) \leq \Psi_2(\mathcal{E}^{\mathrm{TOP}}(D^4,F_{\Gamma}))$.  
\end{proof}

 Let $q_{F_\Gamma}:H_1(F_\Gamma;\Z/2\Z) \rightarrow \Z/2\Z$ denote the mod 2 Seifert form on $F_\Gamma$ (equivalently, Rochlin's form). Symplectic automorphisms of $H_1(F_\Gamma;\Z/2\Z)$ act on the set of quadratic forms on $H_1(F_\Gamma;\Z/2\Z)$ by pullback. 
 
 \begin{defn}
Denote the stabilizer of $q_{F_\Gamma}$ in $\mathrm{Sp}(H_1(F_\Gamma;\Z/2\Z))$ by $\mathrm{Sp}(H_1(F_\Gamma;\Z/2\Z))[q_{F_\Gamma}]$.
 \end{defn}

 For any $c \in C_\Gamma$, we have $q_{F_\Gamma}(c) = 1$. Direct computation shows that if $q_{F_\Gamma}(c) =1$, then  $T_{[c]} \in \mathrm{Sp}(H_1(F_\Gamma;\Z/2\Z))[q_{F_\Gamma}]$.  Let $g$ be the genus of $F_\Gamma$. Consider reflections $\mathrm{Ref}_{\alpha}: \Z^{2g} \rightarrow \Z^{2g}$ given by \[\mathrm{Ref}_{\alpha}(\beta) = \beta - \frac{2 \beta \cdot \alpha}{ \alpha \cdot \alpha }\alpha,\] where $-\cdot-$ denotes the bilinear form represented by $S^T +S$ in the basis $C_\Gamma$. %Since $S^T + S \pmod{2}$ represents the intersection form on $H_1(F_\Gamma;\Z_2)$,
 Then $\mathrm{Ref}_{c} \equiv T_{[c]} \pmod{2}$. Since $W_2(\Gamma)$ is the subgroup generated by the $T_{[c]} : H_1(F_\Gamma;\Z/2\Z) \to H_1(F_\Gamma;\Z/2\Z)$ for $c \in C_\Gamma$, we obtain $W_2(\Gamma) \subset \mathrm{Sp}(H_1(F_\Gamma;\Z/2\Z))[q_{F_\Gamma}]$. 
 
 \begin{rmk}
By order considerations, for $n\geq 3$, $W_2(\Gamma) \subset \mathrm{Sp}(H_1(F_\Gamma;\Z/2\Z))$ is a  proper subgroup if $\Gamma = A_n, \,E_6$ or $E_8$ (recall that $W(A_n) \cong S_{n+1}$, the symmetric group on $n+1$ letters). 
 \end{rmk} 

 We now give an example of a transvection in $\mathrm{Sp}(H_1(F_\Gamma;\Z/2\Z))[q_{F_\Gamma}]$ that is not in $W_2(\Gamma)$ for $\Gamma = A_n$.

 %is in the stabilizer  of the form $q_{F_\Gamma}$ in $\mathrm{Sp}(H_1(F_\Gamma;\Z_2))$. We denote the stabilizer of $q_{F_\Gamma}$ by $\mathrm{Sp}(H_1(F;\Z_2))[q]$. Since $W(\Gamma)$ is the subgroup generated by the $T_{[c]}$ for $c \in C_\Gamma$, we obtain $W(\Gamma) \subset \mathrm{Sp}(H_1(F;\Z_2))[q]$. By order considerations, this is a proper subgroup if $\Gamma = A_n, E_6$ or $E_8$.    

%\textcolor{red}{For $\Gamma = E_6, E_8$, looking at \cite[Fig.7]{Casals_lagrangian_skeleta} gives a configuration of extendable Dehn twist curves with $E_6$ subgraph. These Dehn twists generate a framed mapping class group, and the subgroup generated by them and the Torelli group is the stabilizer of Rochlin's form, so they are spin flexible.}

%\textcolor{red}{The methods in this section only use Proposition 5.7 of Conway--Powell, do the methods here just go through for Seifert surfaces of links. If so we have obstructions for types $D_n$ too.}

\subsection{Obstructing explicit twists for Type $A_n$}
We illustrate the obstruction in Lemma \ref{thm:obtruction_for_pushed_in_Seifert_surfaces} with the following result. %For a quadratic form $q : H_1(F;\Z_2) \to Z_2$, we let $\mathrm{Sp}(H_1(F;\Z_2))[q]$ denote the symplectic automorphisms of $H_1(F;\Z_2)$ which stabilize $q$ under pullback. 

%\begin{figure}[h]
   % \centering
   % \includegraphics[width=0.5\linewidth]{Screenshot 2026-06-05 142521.png}
   % \caption{The Dehn twist about the curve $b$ does not extend}
    %\label{fig:non_extendable_curves}
%\end{figure}

\begin{figure}[h]
  \centering
  \def\svgwidth{0.95\columnwidth}
  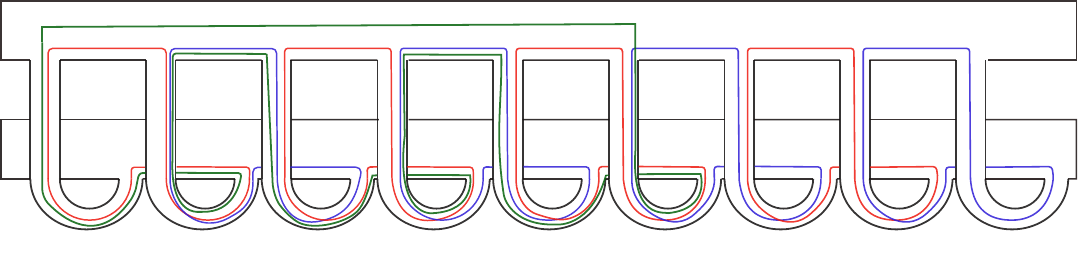
  \caption{\textnormal{The Seifert surface $F_\Gamma$ for $\Gamma = A_8$, together with its distinguished filling configuration $\mathcal{C}_\Gamma = \{\alpha_1,\ldots,\alpha_8\}$ of curves. The curve $b$ is disjoint from $\mathcal{C}_\Gamma\backslash \{\alpha_6\}$ and intersects $\alpha_6$ once. A regular neighbourhood of $b$ in $F_\Gamma$ is an unknotted annulus with three full twists.} }
  \label{fig:non_extendable_curves}
\end{figure}

\begin{proposition} \label{non_extendable_twist} %\textcolor{red}{Make update - works on any with sufficiently large genus}
    Fix $n$ even and $n\geq 8$, and let $F = F_n \subset D^4$ be the minimal genus Seifert surface for the $(2,n+1)$-torus knot with its interior pushed into $D^4$, so $F = F_\Gamma$ for $\Gamma = A_{n}$.  Consider the curve $b \subset F$ depicted in Figure \ref{fig:non_extendable_curves} (for $n=8$). Then $T_{[b]} \in \mathrm{Sp}(H_1(F;\Z/2\Z))[q_{F}]$ but $T_b \notin \mathcal{E}^{\mathrm{TOP}}(D^4,F)$.
\end{proposition}
\begin{proof}
    By construction $\mathrm{lk}(b,b^+) = - 3$, hence $q_{F}([b]) = 1$, so $T_{[b]} \in \mathrm{Sp}(H_1(F;\Z/2\Z))[q_{F}]$. Let $\{[\alpha_1],\ldots,[\alpha_n]\}$ be the homology classes of the chain of curves depicted in Figure \ref{fig:non_extendable_curves} (the homology classes of the curves in $C_\Gamma$ for $\Gamma = A_n$). 
    We have $$T_{[b]}([\alpha_i]) = \begin{cases}
        [\alpha_i] ,\; i \neq 6; \\
        [\alpha_1] + [\alpha_3] + [\alpha_5] + [\alpha_6],\; i=6.
    \end{cases}$$

    Consider the free $\Z$-module with basis $\Delta = \{ [\alpha_i]\}_{i=1}^n$ equipped with the Cartan matrix $P(-1) = S^T+S$. Let $W = \langle \mathrm{Ref}_{[\alpha_i]}\rangle$ be the Weyl group of type $A_n$. Form the root system $\Phi = \bigcup_{w \in W,\; \alpha \in \Delta} \{ w(\alpha)\}$. This consists of all $\alpha \in \Z^n$ with $\alpha \cdot \alpha = -2$, where $\cdot $ is the bilinear form obtained from $P(-1)$. The $-2$-vectors in the $A_{n}$-root lattice are exactly the 
    \begin{equation} \label{eq:consecutive_interval_vector}
        \pm( \alpha_i + \alpha_{i+1} + \cdots + \alpha_j)
    \end{equation} 
    for some $j \geq i$. Any integral matrix $w$ with $w^T P(-1)w = P(-1)$ sends roots to roots, hence the columns of the matrix $w$ take the form (\ref{eq:consecutive_interval_vector}); we call such columns \textit{consecutive interval vectors}.
    
    If $T_b$ extends, then, by Theorem \ref{realizationofWeyl}, $T_{[b]} \equiv w \in W \pmod{2}$ and the column $T_{[b]} \alpha_6$ must be a consecutive interval vector. But $T_{[b]} \alpha_6 \equiv \alpha_1 + \alpha_3 + \alpha_5 + \alpha_6 \pmod{2}$, giving a contradiction.
\end{proof}

\begin{rmk}
We comment on similarities between the obstructions arising in the proof of Proposition \ref{non_extendable_twist} (and Theorem \ref{realizationofWeyl}) and work of A'Campo \cite{ACampo}. A'Campo realized the \textit{monodromy of the singularity} $A_{n}$ as the image of a representation $p\colon B_{n+1} \to \mathrm{Sp}(H_1(F_n;\Z_2)) \cong \mathrm{Sp}_n(\Z_2)$ of the \textit{braid group} on $(n+1)$-strands. A'Campo computed $p(B_{n+1}) = W_2(A_{n})$ \cite[Theorem 1 and Lemma 6]{ACampo}. Another perspective on this obstruction is that it comes from the unitarity of the (reduced) \textit{Burau representation} $B_{n+1} \to \mathrm{GL}_n(\mathbb{Z}[t,t^{-1}])$, reducing at $t = 1$ and $t=-1$. Our obstructions are obtained in a similar fashion.

%coming from this result is that it may be proved in a similar fashion to Theorem \ref{realizationofWeyl} by specializing the (reduced) \textit{Burau representation} \( B_{n+1} \to \mathrm{GL}_n(\mathbb{Z}[t,t^{-1}]) \) to $t = -1$ and $t

%as the mod 2 Weyl group $W_2(A_{n+1})$. In modern language, A'Campos result may be obtained by specializaing 

%principle of specializing an obstruction valued in $\Lambda$ to determine an obstruction  

%here in the $A_n$ case obtained by specializing similarity in techniques between the proof of  Proposition \ref{non_extendable_twist} and A'Campo's  mod 2 homological characterization of \textit{vanishing cycles} on $F_n$ as consecutive interval vectors \cite[Theorem 2]{ACampo}. 
\end{rmk}

\subsection{Flexibility after stabilization} If $h \subset M$ is an embedded 3-dimensional 1-handle intersecting a surface $\Sigma \subset M$ only in its attaching region, an \textit{internal stabilization} of $\Sigma$ replaces $\Sigma \cap h$ with the remaining portion of $\partial h$ (see \cite{InancNathan} for more details). If the knot group of $\Sigma$ is cyclic, stabilization is unique up to isotopy \cite[Lemma 3]{InancNathan}. In particular, it suffices to add a handle supported in a small 4-ball. We call such a handle attachment \textit{trivial}. Up to isotopy, it is equivalent to taking the relative connected sum $(M,\Sigma) \# (S^4,T)$ where $T \subset S^4$ is an unknotted torus.

\setcounter{maintheorem}{3}
\begin{maintheorem} \label{theoremD}
There exist infinitely many non-isotopic pushed in Seifert surfaces in $D^4$ that are not topologically spin-flexible rel boundary, but become topologically spin-flexible after a trivial internal stabilization.
%There exist infinitely many non-isotopic pushed in Seifert surfaces in $B^4$ that are not topologically spin-flexible rel boundary, but become topologically spin-flexible after a single internal stabilization.
\end{maintheorem}

\begin{proof}
Consider the family of Seifert surfaces $F_n$ arising from the $A_{n}$-singularity for $n$ even and $n \geq 8$ (so $g(F_n) = n/2 \geq 4)$. By Proposition \ref{non_extendable_twist}, they are not topologically spin-flexible. Consider the connected sum $F_n \# T$ where $T \subset S^4$ is an unknotted torus. Let $\Delta$ denote the surgery 2-disk for this connected sum. We view the interior of $F_n \#T$ as sitting inside an $S^3$ slice of a collar about $\partial D^4$. Let $\Sigma = F_n \# T \cup D$ be the (abstract) genus $g = n/2+1$ surface obtained by capping off the boundary of $F_n\# T$ with a disk $D$. Extending mapping classes of $F_n\#T$ by the identity yields a homomorphism \[p_D: \Mod(F_n\# T) \to \Mod(\Sigma).\]
Let $UT\Sigma \to \Sigma$ be the unit tangent bundle of $\Sigma$. The kernel of $p_D$ is the \textit{disk-pushing subgroup} $\pi_1(UT\Sigma) \leq \Mod(F_n\#T)$ (see \cite[Section 4.2.5]{farb_margalit}). We have the \textit{disk-pushing sequence}
\[ 1 \to \pi_1(UT\Sigma) \to \Mod(F_n\#T) \overset{p_D}{\to} \Mod(\Sigma) \to 1.\]
Rochlin's form $q_{F_n \# T}$ extends uniquely to a quadratic form $q : H_1(\Sigma;\mathbb{Z}/2\mathbb{Z}) \to \mathbb{Z}/2\mathbb{Z}$. The disk-pushing subgroup acts trivially on $H_1(F_n\# T;\Z) \cong H_1(\Sigma;\Z)$, and so restricting the disk-pushing sequence yields 
\[ 1\to \pi_1(UT\Sigma) \to \Mod(F_n\#T)[q_{F_n\#T}] \overset{p_D}{\to} \Mod(\Sigma)[q]\to 1.\]
 We show \[p_D(\mathcal{E}^{\mathrm{TOP}}(D^4,F_n\#T)) = \Mod(\Sigma)[q].\] By Corollary \ref{torelliextends}, $\mathcal{I}(F_n\#T) \leq \mathcal{E}^{\mathrm{TOP}}(D^4,F_n\#T)$. As $\pi_1(UT\Sigma) \leq \mathcal{I}(F_n\#T)$, it will then follow that
 \[ \mathcal{E}^{\mathrm{TOP}}(D^4,F_n\#T) = \Mod(F_n\#T)[q_{F_n\#T}].\]
We use Lemma \ref{maintechnicallemma} to prove that $p_D(\mathcal{E}^{\mathrm{TOP}}(D^4,F_n\#T)) = \Mod(\Sigma)[q]$. We must show that $p_D(\mathcal{E}^{\mathrm{TOP}}(D^4,F_n\#T))$ contains an $r$-spin mapping class group and a nonseparating genus one bounding pair map. The latter is a consequence of Corollary \ref{torelliextends}.  For the former, we produce the $r$-spin lower bound directly, using the finite generating sets provided by Proposition \ref{gensetrspin}. Let $\mathcal{C} = \{\gamma_1,\ldots,\gamma_n\}$ be the distinguished filling configuration of curves on $F_n$. We start by expanding the configuration $\mathcal{C}$ to fill $F_n \# T$ and to ensure the existence of an $E_6$ subgraph. Choose disjoint nonseparating simple closed curves $\alpha'$ and $\alpha''$ in $F_n$ such that 
\begin{enumerate}
    \item $\alpha'$ is disjoint from $\mathcal{C}\backslash \{\gamma_3\}$ and geometrically dual to $\gamma_3$; %, and $\alpha' \cup \gamma_4 \cup \cdots \cup \gamma_{n-2} \cup \gamma_n$ bounds a $n/2$ times punctured $S^2$ in $F_n$. 
    \item $\alpha''$ is disjoint from $\mathcal{C}\backslash \{\gamma_5\}$, geometrically dual to $\gamma_5$, and $\alpha'\cup \gamma_4 \cup \alpha'' $ bounds a pair of pants in $F_n$. 
\end{enumerate}
For homological reasons, $q_{F_n}(\alpha')$ or $q_{F_n}(\alpha'') = 1$. Let $\alpha \in \{\alpha',\alpha''\}$ be the curve with $q_{F_n}$-value one. The quadratic form $q_{F_n}$ is the mod 2 reduction of the Seifert form on $F_n$ (as is true for any pushed in Seifert surface). Therefore \[\mathrm{lk}(\alpha,\alpha^+) = 2k +1 \] for some $k \in \mathbb{Z}$. Choose a simple closed curve $\delta \subset T\backslash \Delta$ with $\mathrm{lk}(\delta,\delta^+) = -2k \in \mathbb{Z}$ and which intersects the standard meridional loop $m$ of $T\backslash \Delta$ once. Let $\gamma = \alpha +_{\varepsilon}\delta$ be the connected sum of $\alpha$ and $\delta$ along an arc $\varepsilon$ disjoint from the configuration $\mathcal{C}$. As $\delta$ is contained in a ball disjoint from $\alpha$, \[\mathrm{lk}(\gamma,\gamma^+) = \mathrm{lk}(\alpha,\alpha^+) +\mathrm{lk}(\delta,\delta^+) = 2k+1 - 2k = 1.\] Let $\beta = T_{\delta}(m) \subset T\backslash \Delta$. The curve $\beta$ is geometrically dual to $\delta$ and has Seifert framing $\mathrm{lk}(\beta,\beta^+) = 1 \mod 2$, 
since $q_T(\beta) = q_T(m) +q_T(\delta)+1 = 1 \mod 2$ and the restriction of $q_T$ to $H_1(T\backslash \Delta;\mathbb{Z}/2\mathbb{Z})$ is the mod 2 Seifert form.  

Let $\mathcal{D} = \mathcal{C} \cup \{\gamma,\beta\}$. By the Hopf band trick (Proposition \ref{HBtrick}), Dehn twists along curves in $\mathcal{C}$ are (smoothly) extendable. The spin-flexibility of unknots rel boundary (Proposition \ref{spinflexunknot}) implies that $T_{\beta} \in \mathcal{E}^{\mathrm{TOP}}(D^4,F_n\# T)$. Corollary \ref{GeneralHBtrick} applied to the curve $\gamma$ yields \[\Gamma_{\mathcal{D}} := \langle T_{c} : c \in \mathcal{D} \rangle  \leq \mathcal{E}^{\mathrm{TOP}}(D^4,F_n\#T).\] 
By Proposition \ref{gensetrspin},  \[p_D(\Gamma_{\mathcal{D}}) = \Mod(\Sigma)[\psi], \] where $\Mod(\Sigma)[\psi]$ is the stabilizer of an $n = 2g-2$ winding number function $\psi$ on $\Sigma$. Recall that the winding number function $\psi$ is uniquely determined by specifying that $\psi(\mathcal{D}) = 0$, that is, each curve in $\mathcal{D}$ is admissible. Therefore
\[ \Mod(\Sigma)[\psi] \leq p_D(\mathcal{E}^{\mathrm{TOP}}(D^4,F_n\#T)).\]
By Lemma \ref{maintechnicallemma} and Corollary \ref{torelliextends}, 
\[ p_{D}(\mathcal{E}^{\mathrm{TOP}}(D^4,F_n \# T)) = \Mod(\Sigma)[\overline{\psi}],\]
where $\overline{\psi}$ denotes the mod 2 reduction of $\psi$. By construction, the mod 2 reduction $\overline{\psi}$ of $\psi$ is identified with the form $q$ (they agree on the basis $\mathcal{D})$, that is, $\overline{\psi}(c) +1 = q(c)$ for any curve $c \subset \Sigma$. It follows that 
\[ p_D(\mathcal{E}^{\mathrm{TOP}}(D^4,F_n\#T)) = \Mod(\Sigma)[q].\]
\end{proof}

\section{Flexibility within $r$-spin linear systems and fake vanishing cycles}

Let $X$ be a smooth projective surface and $L \in \mathrm{Pic}(X)$ a holomorphic line bundle on $X$.

\subsection{Monodromy of linear systems and higher spin structures} Let $|L| := \mathbb{P}H^0(X,L)$. Recall that elements of $|L|$ may be identified with their zero loci in $X$ (that is, two sections have the same zero locus if and only if they are scalar multiples of one another). Let $U_L \subset |L|$ denote the open subset of $|L|$ consisting of complex curves in $X$ that are smooth in the algebro-geometric sense (the complement of the \textit{discriminant locus}). The tautological family $\mathscr{C}_L$ over $U_L$ is defined by 
\[ \mathscr{C}_L = \{(C,x) \,|\, C \in U_L, x\in C\}.\]
The natural projection $\pi_L : \mathscr{C}_L \to U_L$ equips $\mathscr{C}_L$ with the structure of a smooth fiber bundle over $U_L$. The fibers of $\pi_L$ are the complex curves in $U_L$. In particular, after fixing a $C \in U_L$ there is a \textit{monodromy representation} 
\[ \rho_L : \pi_1(U_L,C) \to  \Mod(C).\]
Geometrically, a loop $\gamma$ in $U_L$ based at $C$ gives a one-parameter family of embedded (real) surfaces in the 4-manifold $X$. By isotopy extension, this determines a diffeomorphism of the pair $(X,C)$ whose restriction to $C$ determines the mapping class $\rho_L(\gamma)$. 
\begin{defn} \label{monodromydefn}
Let $\Gamma_L = \mathrm{im}(\rho_L)$. We refer to $\Gamma_L$ as the \textit{monodromy of the linear system} $|L|$. 
\end{defn}

\begin{corollary} \label{rspinisextendable}
For a smooth curve $C \in |L|$,
\[ \Gamma_L \leq \mathcal{E}^{\mathrm{DIFF}}(X,C).\]
\end{corollary}

In general, the monodromy of a linear system is governed by a \textit{higher spin structure} on $C$. Fix an $r$-th root $E \in \mathrm{Pic}(X)$ of the adjoint line bundle of $L$, that is, $E^{\otimes r} \cong K_X \otimes L$. As $X$ is simply connected, the first Chern class gives a monomorphism $c_1 : \mathrm{Pic}(X) \to H^2(X;\mathbb{Z})$ \cite{GH} and $H^2(X;\mathbb{Z})$ is torsion free. In particular, the $r$-th root $E$ of $K_X \otimes L$ is unique. Adjunction gives an isomorphism \((K_X \otimes L)\rvert_C \cong K_C.\) In particular, we have a fiberwise branched covering of the holomorphic tangent bundle $T_C$ of $C$, 
\[E^{-1}\rvert_C \to (E^{-1})\rvert_C^{\otimes r} \cong T_C,\] on fibers given by $z \to z^r$. This determines a cohomology class $\xi_E \in H^1(UTC;\mathbb{Z}/r\mathbb{Z})$ on the unit tangent bundle $UTC$ of $C$. By Theorem \ref{bijectionbetweenwindingandspin}, $\xi_E$ determines a winding number function \[\psi_E : \mathcal{S}(C) \to \mathbb{Z}/r\mathbb{Z}.\]
The monodromy of the linear system $|L|$ must stabilize this winding number function. 
\begin{proposition}[{\cite[Lemma 3.3]{banerjee2025monodromyvanishingcyclessufficiently}}]
An $r$-th root $E$ of the adjoint line bundle $K_X \otimes L$ determines a winding number function $\psi_E : \mathcal{S}(C) \to \mathbb{Z}/r\mathbb{Z}$, and 
\[ \Gamma_L \leq \Mod(C)[\psi_E].\]
\end{proposition}

\begin{defn} \label{rspinmon}
A linear system $|L|$ has \textit{$r$-spin monodromy} if, for a maximal root $E$ of the adjoint line bundle $K_X \otimes L$, we have \[\Gamma_L = \Mod(C)[\psi_E].\]
\end{defn}

Banerjee and Salter proved that if $\pi_1X = 1$ and $L$ is  \textit{sufficiently ample}, the system $|L|$ has $r$-spin monodromy \cite[Theorem A]{banerjee2025monodromyvanishingcyclessufficiently}. Their result is essential in proving Theorem \ref{maintheorem}. 

\subsection{$k$-jet ample line bundles}

A line bundle $L$ is $0$-jet ample if and only if it is \textit{basepoint free} and 1-jet ample if and only if it is \textit{very ample}. For a precise definition of $k$-\textit{jet ampleness}, we refer the reader to \cite{kjetampleref}. For our purposes, $k$-jet ampleness of $L$ allows us to find a section of $L$ with any prescribed $k$-jet at a point $x\in X$. We make use of the following property when establishing that a bundle is $k$-jet ample.

\begin{proposition} \label{ampleadd}
%The map in Definition \ref{kjetample} is obtained by tensoring an ideal sheaf sequence with $L$. 
 If $L_1$ is $k_1$-jet ample, and $L_2$ is $k_2$-jet ample, then $L_1 \otimes L_2$ is $(k_1+k_2)$-jet ample.
 \end{proposition}

 \begin{proof}
 See \cite[Lemma 2.2]{kjetampleref}.
 \end{proof}

\begin{lemma} \label{localmodel}
Assume that $L \in \mathrm{Pic}(X)$ is 6-jet ample. Then there exists a smooth curve $C \in |L|$ and a closed ball $B \subset X$ such that $C \cap B$ is a $\mathbb{Z}$-surface in $B$ of genus three.  %Let $V_{\varphi}$ denote a versal deformation space of $\varphi$. There is an embedding 
%\[ \iota_p : V_{\varphi} \to H^0(X;L)\]
%and a closed ball $B \subset X$ about $p$ such that $Z(f_0+\iota_p(g))$ is smooth outside $B$, transverse to $\partial B$, and the family of topological surfaces $Z(f_0+\iota_p(g))\backslash B$ is trivial. 
%Moreover, there exists a closed ball $B \subset X$ about $p$ and smooth curve $C \in |L|$ such that the pair $(B,C \cap B)$ is a diffeomorphic to $(D^4_{\varepsilon}(0),D^4_{\varepsilon}(0) \cap \varphi^{-1}(0))$ where $D^4_\varepsilon$ is a ball of radius $\varepsilon$ about $0 \in \mathbb{C}^2$ and $D^4_{\varepsilon}(0) \cap \varphi^{-1}(0)$ is the Milnor fiber of $\varphi$. In particular, $C \cap B$ is a $\mathbb{Z}$-surface in $B$. 
\end{lemma}

\begin{proof}
Fix a point $p \in X$ with coordinate neighbourhood $U \cong \mathbb{C}^2$ and local (analytic) coordinates $x,y$. Let $\varphi(x,y) = x^3+y^4 \in \mathbb{C}[x,y]$. Then there is a section $f_0 \in H^0(X;L)$ such that $Z(f_0)$ has the singularity type of $\varphi$ at $p$ and is smooth everywhere else \cite[Lemma 5.7]{banerjee2025monodromyvanishingcyclessufficiently}. A well-chosen smoothing of $f_0$ does not change the  isotopy type of $Z(f_0)$ outside of a small closed ball $B$ about $p\in X$, but replaces $Z(f_0) \cap B$ with a Milnor fiber $\Sigma_\varphi$ of the singularity $\varphi$. By Proposition \ref{MilnorfiberisZsurface}, the Milnor fiber of the singularity $\varphi$ is a $\mathbb{Z}$-surface in $B$. Moreover, up to isotopy, it is obtained by pushing the interior of a minimal genus Seifert surface for the $(3,4)$-torus knot into the interior of the ball $B$, so $\Sigma_\varphi$ has genus three. 
%We provide some brief commentary. We may assume that $L\rvert_U \cong U \times \mathbb{C}$ and $p \in U \mapsto 0 \in \mathbb{C}^2$. This allows us to identify sections of $L$ over $U$ with holomorphic functions on $U \cong \mathbb{C}^2$ and compute their \textit{jets}. As $L$ is 6-jet ample, it is in particular 4-jet ample.  This guarantees the existence of a global section $f_0 \in H^0(X;L)$ whose 4-jet $j_{f_0}^4(0) = \varphi$. The additional ampleness of $L$ allows us to assume further that this is the only singularity of $Z(f_0)$. As a holomorphic function $\varphi$ is \textit{4-determined}, that is, any holomorphic function $g : \mathbb{C}^2 \to \mathbb{C}$ with an isolated singularity at the origin and whose 4-jet $j_g^4(0) = \varphi$ is \textit{analytically equivalent} to $\varphi$. This follows from \textit{finite determinacy} \cite[Theorem 2.23]{singularities} which states that $\varphi$ is $k$-determined if $\mathfrak{m}^{k+1} \subset \mathfrak{m}^2(\partial_x\varphi,\partial_y\varphi)$. Therefore $Z(f_0)$ has the same singularity type of $\varphi$ at $p$. It follows that a suitable smoothing $C$ of $Z(f_0)$ intersects a closed ball $B \subset U$ in a properly embedded surface such that, for well chosen $\varepsilon,\delta >0$, the pair $(B, C\cap B)$ is diffeomorphic to the pair $(B_{\varepsilon},\varphi^{-1}(\delta)\cap B_{\varepsilon})$ where $B_\varepsilon$ denotes the closed ball of radius $\varepsilon$ in $\mathbb{C}^2$. 
\end{proof}

\begin{rmk}
 Lemma \ref{localmodel} can be extracted from the literature as follows. Replace \cite[Proposition 11]{Banerjee} in the proof of \cite[Lemma 1]{Banerjee} with \cite[Lemma 5.7]{banerjee2025monodromyvanishingcyclessufficiently}. This follows the program suggested at the end of the first section in \cite{Banerjee}. One replaces working with sufficiently large powers of an ample line bundle with sufficient ampleness conditions.  %To see that a deformation is in fact a Milnor fiber, see page 68 of 'Singularities of Differentiable..' Vol 2. 
\end{rmk}

\subsection{Topological flexibility within linear systems} We now prove the main result of this section.

\begin{theorem} \label{thmA}
Let $X$ be a projective surface with $\pi_1(X) = 1$. Let $L \in \mathrm{Pic}(X)$ be 6-jet ample. Assume that $|L|$ has $r$-spin monodromy (Definition \ref{rspinmon}). Then smooth curves in the linear system $\mathbb{P}H^0(X;L)$ are topologically flexible if $r$ is odd, and topologically spin-flexible if $r$ is even.  
\end{theorem}

\begin{proof}
As any two smooth curves in $|L|= \mathbb{P}H^0(X;L)$ are  isotopic in $X$, it suffices to prove the claim for a single smooth curve in $|L|$. By Lemma \ref{localmodel} we can find a smooth curve $C \in |L|$ and an embedded 4-ball $B \subset X$ such that $\Sigma := C\cap B$ is a properly embedded $\mathbb{Z}$-surface in $B$ of genus three. By Corollary \ref{torelliextends}, \[\mathcal{I}(\Sigma) \leq \mathcal{E}^{\mathrm{TOP}}(B,\Sigma).\] Let $T_{a}T_{b}^{-1} \in \mathcal{I}(\Sigma)$ be a nonseparating genus one bounding pair map on $\Sigma$. Let $\psi_E : \mathcal{S}(C) \to \mathbb{Z}/r\mathbb{Z}$ be the mod $r$ winding number function determined by the maximal root $E$ of the adjoint line bundle $K_X \otimes L$. By assumption, \(\Mod(C)[\psi_E] \leq \mathcal{E}^{\mathrm{DIFF}}(X,C).\) 
Therefore \[\langle \Mod(C)[\psi_E],T_a T_b^{-1} \rangle \leq \mathcal{E}^{\mathrm{TOP}}(X,C).\] 
Assume $r$ is odd. By Lemma \ref{maintechnicallemma}, 
\[ \langle \Mod(C)[\psi_E],T_a T_b^{-1} \rangle  = \Mod(C)\]
and so $C$ is topologically flexible. 

Assume $r$ is even. Libgober \cite[Proposition 2.1]{Libgober} identifies the mod 2 reduction $\overline{\psi_E}$ of $\psi_E$ with Rochlin's quadratic form $q_C : H_1(C;\mathbb{Z}/2\mathbb{Z}) \to \mathbb{Z}/2\mathbb{Z}$ ($q_C +1  = \psi_E)$. By Lemma \ref{maintechnicallemma},
\[ \Mod(C)[q_C] \leq \mathcal{E}^{\mathrm{TOP}}(X,C).\]
On the other hand, by Corollary \ref{toppreservesroch}, $\mathcal{E}^{\mathrm{TOP}}(X,C) \leq \Mod(C)[q_C]$, and hence $C$ is topologically spin-flexible. 
\end{proof}

\begin{rmk}[\textbf{Recovering Theorem \ref{maintheorem}}]
If $L = L_1 \otimes L_2$ for $L_1$ very ample (1-jet ample) and $L_2$ 6-jet ample, then \cite[Theorem A]{banerjee2025monodromyvanishingcyclessufficiently} asserts that $|L|$ has $r$-spin monodromy. The bundle $L$ is 7-jet ample, so 6-jet ample. Now apply Theorem \ref{thmA}. 
\end{rmk}

\setcounter{maintheorem}{1}
\begin{maincorollary} \label{flexcomplexcurves}
Smooth plane curves in $\mathbb{C}P^2$ are topologically (spin) flexible. Smooth curves in $\mathbb{C}P^1 \times \mathbb{C}P^1$ of bidegree $(d_1,d_2)$ with $d_1,d_2\geq 6$ are topologically (spin) flexible. % Moreover, there exist (spin) flexible smooth curves in the Hirzebruch surface $\mathbb{F}_1 \cong \mathbb{C}P^2 \# \overline{\mathbb{C}P^2}$, as well as in complete intersections. 
\end{maincorollary}

\begin{proof}
Smooth flexibility for cubics and quartics in $\mathbb{C}P^2$ is due to Hirose \cite{HiroseP2}. Smooth spin-flexibility for quintics in $\mathbb{C}P^2$ follows from work of Salter \cite{salterquintic}. In the remaining cases, we apply Theorem \ref{thmA}. 

For $d \geq6$, smooth degree $d$ plane curves are the zero loci of sections of $\mathcal{O}_{\mathbb{C}P^2}(d) = \mathcal{O}_{\mathbb{C}P^2}(1)^{\otimes d}$, and the line bundle $\mathcal{O}_{\mathbb{C}P^2}(d)$ is $d$-jet ample. It follows from \cite[Theorem 1.3]{SalterInvent} and \cite[Corollary 1.3]{CalderonSalter} that $\mathcal{O}_{\mathbb{C}P^2}(d)$ has $r$-spin monodromy.

%We use \cite[Theorem A]{banerjee2025monodromyvanishingcyclessufficiently}.  %This follows from work of Banerjee-Salter \cite{banerjee2025monodromyvanishingcyclessufficiently}.
%A smooth plane curve $C_d \subset \mathbb{C}P^2$ is the zero locus of a smooth section of $\mathcal{O}_{\mathbb{C}P^2}(d)$. The line bundle $\mathcal{O}_{\mathbb{C}P^2}(d) \cong \mathcal{O}_{\mathbb{C}P^2}(1)^{\otimes d}$ is $d$-jet ample.  The adjoint line bundle of $\mathcal{O}_{\C P^2}(d)$ is $K_{\mathbb{C}P^2}\otimes \mathcal{O}_{\C P^2}(d) \cong \mathcal{O}_{\C P^2}(d-3)$ with maximal root $d-3$. Salter (and remaining cases completed by Calderon-Salter) proved that the linear system $|\mathcal{O}_{\C P^2}(d)|$ has $(d-3)$-spin monodromy. 

A smooth curve $C_{d_1,d_2} \subset \C P^1 \times \C P^1$ of bidegree $(d_1,d_2)$ is the zero locus of a smooth section of $L_{d_1,d_2} = \mathrm{pr}_1^*\mathcal{O}_{\C P^1}(d_1) \otimes \mathrm{pr}_2^*\mathcal{O}_{\C P^1}(d_2)$. The line bundle $L_{d_1,d_2}$ is $k$-jet ample if and only if $d_1 \geq k$ and $d_2\geq k$  \cite[Proposition 5.1]{kjetampleref}. In particular, if $d_1,d_2\geq 6$, the bundle $L_{d_1,d_2}$ is 6-jet ample. By \cite[Theorem A]{SalterInvent} and \cite[Corollary 1.3]{CalderonSalter}, $L_{d_1,d_2}$ has $r$-spin monodromy. The claim now follows from Theorem \ref{thmA}.  
\end{proof}

\subsection{Fake vanishing cycles}

Let $\mathbb{D} \subset \mathbb{C}$ be the closed unit disk.  Recall that a \textit{nodal degeneration} of a smooth curve $C \in |L|$ consists of a family of curves $\{C_t\}_{t\in \mathbb{D}} \subset X$ parameterized by $\mathbb{D}$ such that $C_0$ has a single nodal singularity at $p \in C_0$, each $C_t$ is smooth for $t \neq 0$, and $C_1  =C$. A simple closed curve $a \subset C$ is called a \textit{vanishing cycle} if under parallel transport to $C_0$ within the family  $\{C_t\}_{t\in \mathbb{D}}$ it maps to the point $p$. Under such an isotopy, the curve $a \subset C$ traces out a smoothly embedded disk $D_a \subset X$ with relative Euler number $e_C(D_a) = -1$ (see \cite{gompfstipsicz}). We make the following  definition in the topological category.

%Let $\mathbb{D} \subset \mathbb{C}$ be the closed unit disk.    % complex curve $C \subset X$ has a \textit{nodal singularity} at $p\in X$ if there exists a coordinate neighbourhood $(U,p) \cong (\mathbb{C}^2,0)$ of $p$ such that within $U$ the curve $C$ is identified with the zero locus of $f : \mathbb{C}^2 \to \mathbb{C}$ defined by $f(z,w) = z^2 + w^2$. 

%\begin{defn}[Nodal degeneration]
%Let $f \in H^0(X;L)$ with smooth zero locus $C = Z(f)$. Fix an embedding $\varphi : \mathbb{D} \to H^0(X;L)$. Let $C_t = Z(\varphi(t))$. We say that $\varphi$ is a \textit{nodal degeneration} of $C$ if it satisfies
%\begin{enumerate}
   % \item[(1)] $C_1  = C$,
   %\item[(2)] $C_t$ is smooth for all $t \neq 0$,
   % \item[(3)] The curve $C_0$ has a single node at some point $p \in C_0$, 
  % \item[(4)] $\partial\varphi/\partial\overline{t}\rvert_{t=0}(p) =0$ and $\partial \varphi/\partial t\rvert_{t=0}(p) \neq 0$. 
%\end{enumerate}
%\end{defn}

%For $t \neq 0$, the curve $C_t$ intersects a small ball about $p$ in a cylinder. As $t \to 0$, a core curve of this cylinder collapses along an embedded disk to the point $p \in X$. We call such curves \textit{vanishing cycles}. 

%\begin{theorem}
%Vanishing cycle if and only if admissible. 
%\end{theorem}
%A simple closed curve $a \subset C$ is called a \textit{vanishing cycle} if there exists a nodal degeneration $\{C_t\}_{t \in \mathbb{D}}$ of $C$ such that along some arc in $\mathbb{D}$ with endpoints at 1 and 0 the curve $a$ 

\begin{defn} \label{fakecycle}
Let $\Sigma \subset M$ be a properly embedded surface in a topological 4-manifold $M$. A simple closed curve $a \subset \Sigma$ is a \textit{topological vanishing cycle} if it bounds a topologically embedded, locally flat disk $D \subset M$ with interior disjoint from $\Sigma$ and relative Euler number $e_\Sigma(D) = \pm1$.
\end{defn}

The terminology in Definition \ref{fakecycle} is motivated by Proposition \ref{HBtrick}. Indeed, if $\gamma \subset C$ is a vanishing cycle, the associated nodal degeneration can be used to realize the (right-handed) Dehn twist about $\gamma$ as an extendable mapping class. %Fix a linear system $|L|$ on a smooth surface $X$ and a smooth curve $C \in |L|$. 

\begin{defn}
A simple closed curve $a \subset C \in |L|$ is called a \textit{fake vanishing cycle} if it is a topological vanishing cycle but not a vanishing cycle on $C$. 
\end{defn}

We now use Theorem \ref{maintheorem} to establish the existence of fake vanishing cycles. We need the following result which, for $r$-spin linear systems, characterizes  vanishing cycles as admissible curves. 

\begin{proposition}[{\cite[Corollary 1.6]{banerjee2025monodromyvanishingcyclessufficiently}}]
\label{vanishingiffadmissible}
Let $L \in \mathrm{Pic}(X)$ and $C \in |L|$ a smooth curve. Let $\psi_E$ be the winding number function determined by a maximal root of the adjoint line bundle. Assume that $|L|$ has $r$-spin monodromy, that is, $\Gamma_L = \Mod(C)[\psi_E]$. Then a nonseparating simple closed curve $a \subset C$ is a vanishing cycle if and only if $\psi_E(a) = 0$. 
\end{proposition}

\setcounter{maintheorem}{2}
\begin{maincorollary} \label{fakevanishingcycle}
Let $X$ be a projective surface with $\pi_1(X) = 1$. Let $L \in \mathrm{Pic}(X)$ satisfy the hypothesis of Theorem \ref{maintheorem}. Fix a smooth curve $C \in |L|$. If the maximal root of the adjoint line bundle $K_X \otimes L$ in $\mathrm{Pic}(X)$ has order at least three, then there exists a fake vanishing cycle $a \subset C$. 
\end{maincorollary}

\begin{proof}
By \cite[Theorem A]{banerjee2025monodromyvanishingcyclessufficiently}, the linear system $|L|$ has $r$-spin monodromy, where $r$ is the order of the maximal root of the adjoint line bundle $K_X \otimes L$. By assumption, $r \neq 1,2$. Let $E$ be a maximal root of $K_X \otimes L$ and fix  $C \in |L|$ a smooth curve. Let $\psi_E \colon \mathcal{S}(C) \to \mathbb{Z}/r\mathbb{Z}$ be the associated winding number function. As $L$ is 6-jet ample, the genus of $C$ is at least three (Lemma \ref{localmodel}). By the $r$-spin change of coordinates principle (requiring genus at least two, see \cite[Section 5.2]{CalderonSalter}), there exists a nonseparating simple closed curve $a \subset C$ with $\psi_E(a)\neq 0 \in \mathbb{Z}/r\mathbb{Z}$. If $r$ is even, we can arrange for $\psi_E(a) = 0 \mod 2$. Let $c \subset C$ be an admissible curve, that is, $\psi_E(c) = 0$. By Proposition \ref{vanishingiffadmissible}, $c$ is a vanishing cycle. By Theorem \ref{thmA}, $C$ is topologically (spin) flexible. In particular, there exists an orientation preserving homeomorphism of pairs $f : (X,C)\to (X,C)$ such that $f(c)=a$. This realizes the curve $a \subset C$ as a fake vanishing cycle.  
\end{proof}

\bibliographystyle{alpha} 
\bibliography{template}
\end{document}